\documentclass[11pt,twoside]{amsart}
\usepackage{geometry}
\usepackage{amsmath,amssymb,esint,amsthm,bbm}
\usepackage[active]{srcltx}
\usepackage[dvips]{graphicx}
\catcode`@=11 \@addtoreset{equation}{section} \catcode`@=12

\email{dario.monticelli@gmail.com (Dario Daniele Monticelli)}
\email{scott.rodney@gmail.com (Scott Rodney)}
\email{wheeden@math.rutgers.edu (Richard L. Wheeden)}

\DeclareMathOperator*{\esup}{\text{ess sup }}
\DeclareMathOperator*{\einf}{\text{ess inf }}
\DeclareMathOperator*{\tesup}{\emph{ess sup }}
\DeclareMathOperator*{\teinf}{\emph{ess inf }}

\newcommand{\X}{\mathcal{X}}
\newcommand{\e}{\varepsilon}

\newcommand{\field}[1]{\mathbf{#1}}
\newcommand{\vphi}{\varphi}
\newcommand{\F}{\mathcal{F}}
\newcommand{\Z}{\field{Z}}
\newcommand{\R}{\field{R}}

\newcommand{\N}{\field{N}}

\newcommand{\E}{\mathcal{E}}

\newcommand{\wt}{\widetilde}

\newcommand{\bea}{\begin{eqnarray}}
\newcommand{\cend} {\end{center}}

\newcommand{\eea}{\end{eqnarray}}

\newcommand{\fracc}{\displaystyle\frac}

\newcommand{\disp}{\displaystyle}

\newcommand{\ra}{\rightarrow}

\newcommand{\beginc} {\begin{center}}

\newcommand{\ep} {\varepsilon}

\newcommand{\avnormb}[3]{\| {#1} \|_{{#2} ,{#3};\overline{dx}}}

\newcommand{\gb}{\mathfrak{b}}

\newtheorem{thm}{\textbf{Theorem}}[section]
\newtheorem{lem}[thm]{\textbf{Lemma}}
\newtheorem{pro}[thm]{\textbf{Proposition}}
\newtheorem{rem}[thm]{\textbf{Remark}}
\newtheorem{cor}[thm]{\textbf{Corollary}}
\newtheorem{defn}[thm]{\textbf{Definition}}
\theoremstyle{remark}

\theoremstyle{definition}
\newtheoremstyle{Claim}{}{}{\itshape}{}{\itshape\bfseries}{:}{ }{#1}
\theoremstyle{Claim}

\begin{document}

\title [Harnack's Inequality and H\"older Continuity]{Harnack's
inequality and H\"older continuity for weak solutions of degenerate
quasilinear equations with rough coefficients} 
\vspace{-0,3cm}



\maketitle

\begin{center}
\textsc{\textmd{D. D. Monticelli\footnote{Dipartimento di Matematica
"F. Enriques", Universit\`a degli Studi, Via C. Saldini 50,
20133--Milano, Italy. Member of the Gruppo Nazionale per l'Analisi
Matematica, la Probabilit\`a e le loro Applicazioni (GNAMPA) of
the Istituto Nazionale di Alta Matematica (INdAM). Partially supported by
GNAMPA, projects ``Equazioni differenziali con invarianze in analisi
globale'' and ``Analisi globale ed operatori degeneri''.},
      S. Rodney\footnote{Department of Mathematics, Physics and
      Geology, Cape Breton University, P.O. Box 5300, Sydney, Nova
      Scotia B1P 6L2, Canada.} and R. L. Wheeden\footnote{Department of
      Mathematics, Rutgers University, Piscataway, NJ 08854, USA.}}}
\end{center}

\begin{abstract}
We continue to study regularity results for weak solutions of the
large class of second order degenerate quasilinear equations of the form
\bea \text{div}\big(A(x,u,\nabla u)\big) = B(x,u,\nabla u)\text{ for }x\in\Omega\nonumber
\eea
as considered in our paper \cite{MRW}. There we proved only local boundedness of weak
solutions. Here we derive a version of Harnack's inequality as well
as local H\"older continuity for weak solutions. The possible
degeneracy of an equation in the class is expressed in terms of a
nonnegative definite quadratic form associated with its principal
part. No smoothness is required of either the quadratic form or the
coefficients of the equation. Our results extend ones obtained
by J. Serrin \cite{S} and N. Trudinger \cite{T} for quasilinear
equations, as well as ones for subelliptic linear equations obtained in
\cite[2]{SW1}.
\end{abstract}

\maketitle

\section{\textbf{Introduction}}\label{intro}

\subsection{General Comments} Our main goal is to prove Harnack's
inequality and local H\"older continuity for weak solutions $u$ of
quasilinear equations of the form
\begin{equation} \label{eqdiff}
\text{div}\big(A(x,u,\nabla u)\big) = B(x,u, \nabla u)
\end{equation}
in an open set $\Omega \subset \R^n$. The vector-valued function $A$
and the scalar function $B$ will be assumed to satisfy the same
structural conditions as in our earlier paper \cite{MRW}, where we
proved that weak solutions are locally bounded.  The possible
degeneracy of equation (\ref{eqdiff}) is expressed in terms of a
matrix $Q(x)$, that may vanish or become singular, associated with the
functions $A,B$.  More precisely, given $p$ with $1<p<\infty$ and an
$n\times n$ nonnegative definite symmetric matrix $Q(x)$ satisfying
$|Q| \in L_{loc}^{p/2}(\Omega)$, we assume the following structural
conditions: For $(x,z,\xi)\in \Omega\times\R \times\R^n$, there is a
vector ${\tilde A}(x,z,\xi)$ with values in $\R^n$ such that for
a.e. $x\in \Omega$ and all $(z,\xi)\in \R\times\R^n$,
\begin{equation}\label{struct}
\begin{cases}
(i)\quad A(x,z,\xi)=\sqrt{Q(x)}{\tilde A}(x,z,\xi),\\
(ii)\quad\xi \cdot A(x,z,\xi) \geq a^{-1}\Big|\sqrt{Q(x)}\,
\xi\Big|^p
- h(x)|z|^\gamma - g(x),\\
(iii)\quad\Big|\tilde{A}(x,z,\xi)\Big| \leq a\Big|\sqrt{Q(x)}\,
\xi\Big|^{p-1} + b(x)|z|^{\gamma -1} + e(x),\\

(iv)\quad\Big|B(x,z,\xi)\Big|\leq c(x)\Big|\sqrt{Q(x)}\,
\xi\Big|^{\psi-1} +d(x)|z|^{\delta-1}+f(x),
\end{cases}
\end{equation}
where $a, \gamma, \psi, \delta >1$ are constants, and $b, c, d, e, f,
g, h$ are nonnegative measurable functions of $x\in\Omega$.

The sizes of the exponents are restricted to the ranges
\begin{equation}\label{ranges}
\gamma\in (1,\sigma(p-1)+1),\quad \psi\in (1,p+1-\sigma^{-1}),\quad
\delta \in (1,p\sigma),
\end{equation}
where $\sigma>1$ is a constant that measures the gain in
integrability in a naturally associated Sobolev estimate (see
(\ref{sobolev}) below). For the classical Euclidean metric $|x-y|$,
nondegenerate $Q$ and $1<p<n$, the Sobolev gain factor $\sigma$ is
$n/(n-p)$. Furthermore, the functions $b, c, d, e, f, g, h$ will be
assumed to lie in certain Lebesgue or Morrey spaces, and to satisfy
the minimal integrability conditions
\begin{equation}\label{integrability}
 \begin{array}{c}
  \begin{array}{ccc}
  c\in L^\frac{\sigma p}{\sigma
   p-1-\sigma(\psi-1)}_\text{loc}(\Omega),\,\,\,&\,\,\,e\in
   L^{p^\prime}_\text{loc}(\Omega),\,\,\,&\,\,\,f\in L^{(\sigma
    p)^\prime}_\text{loc}(\Omega),
  \end{array}\\
  \begin{array}{cc}
  b\in L^\frac{\sigma
  p}{\sigma\!(p-1)-\gamma+1}_\text{loc}(\Omega),\,\,\,&\,\,\,d\in
  L^\frac{p\sigma}{p\sigma-\delta}_\text{loc}(\Omega).
  \end{array}
 \end{array}
\end{equation}
Here and elsewhere we use a prime to denote the dual exponent, for
example, $1/p + 1/p' =1$ when $1\le p \le \infty$, with the standard
convention that $1$ and $\infty$ are dual exponents.

The quadratic form associated with $Q(x)$ will be denoted
\begin{equation}\label{quadform}
Q(x, \xi) = \langle Q(x)\xi, \xi \rangle, \quad (x,\xi) \in
\Omega\times\R^n,
\end{equation}
and we note that $Q(x,\xi)$ may vanish when $\xi \neq 0$, i.e., $Q(x)$
may be singular (degenerate). As in \cite{MRW} and
following \cite[2]{SW1}, our weak solutions are pairs $(u, \nabla
u)$ which belong to an appropriate Banach space
$\mathcal{W}^{1,p}_Q(\Omega)$ obtained by isomorphism from the
degenerate Sobolev space $W^{1,p}_Q(\Omega)$, defined as the
completion with respect to the norm
\begin{equation}\label{Wnorm}
 ||u||_{W_Q^{1,p}(\Omega)} = \left(\int_\Omega |u|^p \,dx +
   \int_\Omega Q(x, \nabla u)^{\frac{p}{2}} dx\right)^{\frac{1}{p}}
\end{equation}
of the class of functions in $\text{Lip}_{\text{loc}}(\Omega)$ with
finite $W^{1,p}_Q(\Omega)$ norm. Technical facts about these Banach
spaces are given in \cite{MRW}, \cite[2]{SW1} with weighted versions
in \cite{CRW}, and some of them will be recalled below. For now, we
mention only that when $Q$ is degenerate, it is important to think of
an element of the Banach space $\mathcal{W}^{1,p}_Q(\Omega)$
as a pair $(u,\nabla u)$ rather than as just the first component $u$,
due to the possibility that $\nabla u$ may not be uniquely determined
by $u$.  Nonuniqueness of $\nabla u$ causes us little difficulty since our
primary regularity results concern estimates of $u$ rather than
$\nabla u$. Except for the need to consider a pair, the notions of
weak solution, weak supersolution and weak subsolution that we will
use are standard, namely, we say that a pair $(u,\nabla u)\in
W^{1,p}_Q(\Omega)$ satisfies
\bea \label{subsupsol} \mbox{div}(A(x,u,\nabla u)) &=(\leq\,,\, \geq)&
B(x,u,\nabla u) \quad \text{for $x\in\Omega$}
\eea
in the weak sense if for every nonnegative test function $\vphi\in
Lip_0(\Omega)$, the corresponding integral expression
\begin{equation}\label{weaksolution}
\int_\Omega \big[\nabla\varphi\cdot A(x,u,\nabla u) +\varphi
B(x,u,\nabla u)\big]\,dx =(\geq\,,\, \leq)\, 0
\end{equation}
holds. The integrals in (\ref{weaksolution}) converge absolutely due
to (\ref{struct})--(\ref{integrability}); see \cite[Proposition 2.5,
Corollary 2.6, Proposition 2.7]{MRW}.

Our results and analysis are carried out in the context of a
quasimetric $\rho$ on $\Omega$, that is, $\rho:\Omega\times\Omega\ra
[0,\infty)$ and satisfies the following for all $x,y,z\in\Omega$:
\bea \label{tri}
\bullet &&\rho(x,y) = \rho(y,x) \;(\mbox{symmetry}),\nonumber\\
\bullet &&\rho(x,y) = 0 \iff x=y \;(\mbox{positivity}),\nonumber\\
\bullet &&\rho(x,y) \leq \kappa[\rho(x,z) + \rho(y,z)]\;
(\mbox{triangle inequality}),
\eea
where $\kappa\geq 1$ is independent of $x,y,z\in\Omega$.  In particular, we
will assume that appropriate Sobolev-Poincar\'e estimates hold and that
Lipschitz cutoff functions exist for the class of quasimetric
$\rho$-balls defined for $x\in \Omega$ and $r>0$ by
\bea \label{rhoball}
B(x,r)&=& \;\;\{y\in\Omega : \rho(x,y) <r\}.
\eea
We will refer to $B(x,r)$ as the $\rho$-ball of radius $r>0$ and center
$x$. All $\rho$-balls lie in $\Omega$ by their definition, and they
are assumed to be open with respect to the usual Euclidean
topology. The estimates we need are summarized in \S 2.

\subsection{Some Known Results} In the standard elliptic case when $Q(x) =$
Identity and $\rho(x,y) = |x-y|$ is the ordinary Euclidean metric,
regularity results including Harnack's inequality and local H\"older
continuity for weak solutions of (\ref{eqdiff}) were derived in
\cite{S} and \cite{T} under structural conditions more restrictive than
(\ref{struct}). Obtaining analogues of these results in the degenerate
case is our main concern.

In the degenerate (or subelliptic) case, Harnack's inequality and
H\"older continuity have been studied in \cite[2]{SW1} for {\it
  linear} equations with rough coefficients and nonhomogeneous terms,
and those results are included among the ones we derive
here. Moreover, in the degenerate {\it quasilinear} case, and
under the same structural assumptions as in (\ref{struct}), local
boundedness of weak solutions is proved in \cite{MRW}. In fact, a rich
variety of local boundedness estimates is given there depending on the
strength and type of condition imposed on the coefficients, but still
without any assumption about their differentiability.

In order to describe a known estimate in the degenerate quasilinear
case, we now record (without listing the precise technical data)
a fairly typical form of the local boundedness estimates proved in
\cite{MRW} in case $\gamma = \delta =p$ and $\psi \in [p,
p+1-\sigma^{-1})$ : If $(u,\nabla u) \in W^{1,p}_Q(\Omega)$ is a weak
solution of (\ref{eqdiff}) in a $\rho$-ball $B(y,r)$, then for any $k>0$, there are positive constants $\tau, C,$ and $\bar{Z}$ such that
\begin{equation}\label{typicalbound}
\disp\esup_{x\in B(y,\tau r)} \big(|u(x)| + k\big) \le C
\bar{Z}\,\left(\frac{1}{|B(y,r)|} \int_{B(y,r)} \big(|u(x)|+k\big)^p
\,dx\right)^\frac{1}{p}.
\end{equation}
Here, $\tau$ and $C$ are independent of
$u, k, B(y,r), b, c, d, e, f, g$ and $h$, but $\bar{Z}$
generally depends on all these quantities in very specific ways
described in \cite{MRW} and later in this paper. The richness of
boundedness estimates that we mentioned above results from estimating
$\bar{Z}$ under various assumptions on the coefficients. In fact, the
estimates in Corollaries 1.8--1.11 of \cite{MRW} offer only a sample
of those which are possible. Understanding $\bar{Z}$, removing its
dependence on $u$ and some of the other data, and generalizing the
mean-value estimates which lead to (\ref{typicalbound}) are important
ingredients in deriving the regularity results in this paper, where in the
broad sense we follow the Moser method.

In order to state our results carefully, including
\eqref{typicalbound}, we must describe the technical background,
which is considerable. This is done in the next section.

\section{Technical Background and Hypotheses}

Our principal results are axiomatic in nature and based mainly on the
existence of appropriate Sobolev-Poincar\'e inequalities and Lipschitz
cutoff functions in a space of homogeneous type. In this section, we
describe the setting for our work and list our main assumptions.

\subsection{Homogeneous Spaces}

\vspace{0.25in}

Let $\Omega\subset\mathbb{R}^n$ be an open set and $\rho$ be a
quasimetric defined on $\Omega$ satisfying (\ref{tri}). We will make
two a priori assumptions relating the $\rho$-balls defined in
(\ref{rhoball}) and the Euclidean balls \bea D(x,r)&=&
\;\;\{y\in\Omega : |x-y|<r\}.\nonumber \eea Note that $D(x,r)$ is
the intersection with $\Omega$ of the ordinary Euclidean ball with
center $x$ and radius $r$, and recall that all $\rho$-balls are also
subsets of $\Omega$. As we already mentioned, we will always assume
that every $B(x,r)$ is an open set according to the Euclidean
topology. Second, we will always assume that \bea\label{Cond1}&&
\text{for all } x\in\Omega, \, |x-y|\rightarrow 0 \text{ if }
\rho(x,y)\rightarrow 0.
\end{eqnarray}
As a consequence of (\ref{Cond1}), for every $x\in \Omega$ there
exists $R_0(x)>0$ such that the Euclidean closure $\overline{B(x,r)}$ of
$B(x,r)$ satisfies $\overline{B(x,r)}\subset\Omega$ for all
$0<r<R_0(x)$. See Lemma 2.1 of \cite{MRW} for this result.

\begin{rem}\label{rem2.1} Since $\rho$-balls are assumed to be open
sets, the converse of \eqref{Cond1} automatically holds:
\bea\label{Cond0}&& \text{for all } x,y\in\Omega,\, \rho(x,y)\rightarrow 0
\text{ if } |x-y|\rightarrow 0.
\end{eqnarray}
Furthermore, since $\rho$-balls are open, every $\rho$-ball has positive
Lebesgue measure.
\end{rem}

As is well-known, the triangle inequality (\ref{tri}) implies that
$\rho$-balls have the following swallowing property (see
e.g. \cite[Observation 2.1]{CW1} for the simple proof):

\begin{lem}\label{swalemma} If $x,y\in\Omega$,  $0<t\leq r$
and $B(y,t)\cap B(x,r)\neq \emptyset$, then
\bea\label{swallowing}
B(y,t) \subset  B(x,\gamma^*r)
\eea
where $\gamma^* = \kappa + 2\kappa^2$ with $\kappa$ as in \eqref{tri}.
\end{lem}

\begin{rem}\label{centerswallow} The constant $\gamma^*$ in the
conclusion of Lemma \ref{swalemma} can be decreased if we only require
information about the center of the smaller ball.  Indeed, if
$x,y\in\Omega$, $0<t\leq r$, and $B(y,t)\cap B(x,r)\neq \emptyset$,
then $y\in B(x,2\kappa r)$ by \eqref{tri}.
\end{rem}

\begin{defn} \label{homspace}  We call the triple $(\Omega, \rho,dx)$
a \emph{local homogeneous space} if Lebesgue measure is locally a
doubling measure for $\rho$-balls, i.e., if there are constants $C_0,d_0>0$
and a function $R_1:\Omega\ra (0,\infty)$ such that if $x,y\in\Omega$,
$0<t\leq r <R_1(x)$ and $B(y,t)\cap B(x,r)\neq\emptyset$, then
\bea\label{doubling1} |B(x,r)| \leq  C_0\Big(\frac{r}{t}\Big)^{d_0}
  |B(y,t)|.
\eea
\end{defn}

\noindent This notion generalizes that of a symmetric general homogeneous space
as defined in \cite[p. 71]{SW1}. Also, due to the swallowing property,
\eqref{doubling1} has an equivalent form:  There are constants
$C_0',c'>0$ such that if $x,y\in \Omega$, $0<t\leq r<c'R_1(x)$ and
$B(y,t)\subset B(x,r)$, then
\bea \label{doubling0} |B(x,r)|\leq
  C_0'\Big(\fracc{r}{t}\Big)^{d_0}|B(y,t)|
\eea for the same $d_0$ and $R_1(x)$ as in (\ref{doubling1}).

\begin{rem}\label{rem1-KMR}
By a result of Korobenko-Maldonado-Rios (see \cite{KRM}), the
validity of the local doubling condition \eqref{doubling1} for some
exponent $d_0>0$ and function $R_1(x)>0$ is a consequence of two
conditions that will be introduced below: the local
Sobolev inequality \eqref{sobolev} and the existence of appropriate
sequences of Lipschitz cutoff functions, supported in pseudometric
balls with small radius and adapted to the matrix $Q$, as described in
\eqref{cutoff}.
\end{rem}

We will usually require that $R_1(x)$, as well as similar functions we
will use to restrict sizes of radii, satisfies the local comparability
condition described in the next definition.

\begin{defn}\label{unifrD} Let $E\subset\Omega$. We say that a
function $f:\Omega\ra (0,\infty)$  satisfies a \emph{local uniformity
condition with respect to $\rho$ in $E$} if there is a constant $A_* =
A_*(f,E)\in (0,1)$ such that for all $x\in E$ and all $y\in  B(x,f(x))$,
\bea \label{unifr}
A_*<\frac{f(y)}{f(x)}<\fracc{1}{A_*}.
\eea
\end{defn}

Condition (\ref{unifr}) is automatically true in case $f$ is bounded
above on $E$ and also has a positive lower bound on $E$.  This
condition will be helpful in our proof of the John-Nirenberg estimate
using techniques related to those in \cite{SW1}.  It is not required
in \cite{SW1} since there, $R_0(x),R_1(x)$ (and $R_2(x)$ in \S 2.2
below) are chosen to be the same fixed multiple of the Euclidean
distance $\mbox{dist}(x,\partial\Omega)$ and so (\ref{unifr}) holds
with $f(x) = R_0(x)= R_1(x) = R_2(x)$ on any set $E$ satisfying
$\overline{E}\subset \Omega$.  In some proofs to follow we will choose
$E$ to be a specific quasimetric ball $B(z,r)$.

\subsection{Poincar\'e-Sobolev Estimates and Cutoff Functions}

\vspace{.25in}

Let $p$ and $Q$ be as in \eqref{struct}, and recall that $p \in
(1,\infty)$ and $|Q|\in L^{p/2}_{loc}(\Omega)$.  Before we state the
Sobolev and Poincar\'e estimates that we require, let us make a few
more comments about the Sobolev space $W^{1,p}_{Q}(\Omega)$. A fuller
discussion can be found in \cite{MRW}, \cite{SW2}, and \cite{CRW}. Let
$Lip_{Q,p}(\Omega)$ denote the class of locally Lipschitz functions
with finite $W^{1,p}_Q(\Omega)$ norm; see (\ref{Wnorm}). The space
$W^{1,p}_Q(\Omega)$ is by definition the  Banach space of equivalence
classes of sequences in $Lip_{Q,p}(\Omega)$ which are Cauchy sequences
with respect to the norm (\ref{Wnorm}). Here two Cauchy sequences are
called equivalent if they are equiconvergent in $W^{1,p}_Q(\Omega)$.

To further describe $W^{1,p}_Q(\Omega)$, we consider the form-weighted
space consisting of all (Lebesgue) measurable $\R^n$-valued functions
$\mathbf{f}(x)$ defined in $\Omega$ for which
\begin{equation}\label{calLp}
||\mathbf{f}||_{\mathcal{L}^p(\Omega,Q)} = \left\{\int_\Omega Q(x,
  \mathbf{f}(x))^{\frac{p}{2}} \,dx\right\}^{\frac{1}{p}} < \infty.
\end{equation}
By identifying any two measurable $\R^n$-valued functions $\mathbf{f}$ and
$\mathbf{g}$ with
$||\mathbf{f}-\mathbf{g}||_{\mathcal{L}^p(\Omega,Q)} =0$,
(\ref{calLp}) defines a norm on the resulting Banach space of
equivalence classes.  We denote this Banach
space of equivalence classes by $\mathcal{L}^p(\Omega,Q)$. If
$\{w_k\}=\{w_k\}_{k=1}^\infty \in W^{1,p}_{Q}(\Omega)$, meaning that
$\{w_k\}$ is a Cauchy sequence of $Lip_{Q,p}(\Omega)$ functions with
respect to (\ref{Wnorm}), then there is a unique pair $(w, \mathbf{v})
\in L^p(\Omega) \times \mathcal{L}^p(\Omega,Q)$ such that $w_k
\rightarrow w$ in $L^p(\Omega)$ and $\nabla w_k \rightarrow
\mathbf{v}$ in $\mathcal{L}^p(\Omega,Q)$. The pair $(w,\mathbf{v})$
represents the particular equivalence class in $W^{1,p}_Q(\Omega)$
containing $\{w_k\}$. The space $\mathcal{W}^{1,p}_Q(\Omega)$
is defined to be the collection of all pairs $(w,\mathbf{v})$ that
represent equivalence classes in $W^{1,p}_Q(\Omega)$. Thus,
$\mathcal{W}^{1,p}_Q(\Omega)$ is the image of the isomorphism
$\mathcal{J}:W^{1,p}_Q(\Omega)\ra L^p(\Omega)\times
\mathcal{L}^p(\Omega,Q)$ defined by
\bea \mathcal{J}([\{w_k\}]) = (w,\mathbf{v}),\nonumber
\eea
where $[\{w_k\}]$ denotes the equivalence class in
$W^{1,p}_Q(\Omega)$ containing $\{w_k\}$. Therefore,
$\mathcal{W}^{1,p}_Q(\Omega)$ is a closed subspace of
$L^p(\Omega)\times \mathcal{L}^p(\Omega,Q)$ and hence a Banach space
as well.  Since ${\mathcal{W}}^{1,p}_Q(\Omega)$ and
$W^{1,p}_Q(\Omega)$ are isomorphic, we will often refer to elements
$(w,\mathbf{v})$ of $\mathcal{W}^{1,p}_Q(\Omega)$ as elements of
$W^{1,p}_Q(\Omega)$. Interestingly, $\mathbf{v}$ is generally not
uniquely determined by $w$ for pairs $(w, \mathbf{v})$ in
$\mathcal{W}_Q^{1,p}(\Omega)$, i.e., the projection
\[
P: \mathcal{W}^{1,p}(\Omega) \rightarrow L^p(\Omega)
\]
defined by $P((w,\mathbf{v}))=w$ is not always an injection; see
\cite{FKS} for an example. However, we will generally abuse notation
and denote pairs in $W^{1,p}_Q(\Omega)$ by $(w, \nabla w)$ instead of
$(w, \mathbf{v})$.  \\

$(W^{1,p}_Q)_0(\Omega)$ will denote the space analogous to
$W^{1,p}_Q(\Omega)$ but where the completion is formed by using
Lipschitz functions with compact support in $\Omega$. A typical
element of $(W^{1,p}_Q)_0(\Omega)$ may be thought of as a pair $(w,
\nabla w) \in L^p(\Omega) \times \mathcal{L}^p(\Omega,Q)$ for which
there is a sequence $\{w_k\}\subset Lip_{Q,p}(\Omega) \cap
Lip_0(\Omega)$ such that $w_k \rightarrow w$ in $L^p(\Omega)$ and
$\nabla w_k \rightarrow \nabla w$ in $\mathcal{L}^p(\Omega,Q)$. Here
we again adopt the abuse of notation $\nabla w$ for the second component
$\mathbf{v}$ of a pair $(w, \mathbf{v})$. \\

We can now state the Sobolev-Poincar\'e estimates that we will
assume.  We say that a \emph{local Sobolev inequality} holds in $\Omega$ if there exists a function $R_2:\Omega\ra (0,\infty)$ and constants $C_1>0$ and $\sigma>1$ such that for every $\rho$-ball
$B(y,r)$ with $0<r<R_2(y)$, the inequality
\bea\label{sobolev}
\Big(\frac{1}{|B(y,r)|}\int_{B(y,r)} |w|^{p\sigma}dx\Big)^\frac{1}{p\sigma}
&\leq&C_1\Big[r\Big(\frac{1}{|B(y,r)|}\int_{B(y,r)} |\sqrt{Q}\nabla
  w|^\frac{p}{2}dx\Big)^\frac{1}{p}\\
  &&\;\;\;\;\;\; +\;\Big(\frac{1}{|B(y,r)|}\int_{B(y,r)}
  |w|^pdx\Big)^\frac{1}{p}\Big] \nonumber
\eea
holds for all $(w,\nabla w) \in (W^{1,p}_Q)_0(B(y,r))$.\\

We say that a \emph{local Poincar\'e inequality} holds in $\Omega$
if there are constants $C_2>0$ and $\gb\geq 1$ such that for every
$\rho$-ball $B(y,r)$ with $0<r<R_2(y)$, the inequality
\bea\label{poincare}
\frac{1}{|B(y,r)|}\int_{B(y,r)} |w-w_{B(y,r)}| dx \leq
C_2r\Big(\frac{1}{|B(y,\gb r)|}\int_{B(y,\gb r)} |\sqrt{Q}\nabla
w|^pdx\Big)^{\frac{1}{p}}
\eea
holds for all $(w,\nabla w)\in W^{1,p}_Q(\Omega)$, where $w_{B(y,r)} =
\fracc{1}{|B(y,r)|}\int_{B(y,r)} wdx$.\\

\begin{rem}It is easy to see that \eqref{sobolev} and \eqref{poincare}
hold as stated, that is, for all $(w,\nabla w)$ in
$(W^{1,p}_Q)_0(B(y,r))$ or $W^{1,p}_Q(\Omega)$ respectively, provided
they hold for all $w$ in $Lip_{Q,p}(\Omega) \cap Lip_0(B(y,r))$ or
$Lip_{Q,p}(\Omega)$ respectively.
\end{rem}

As in \cite{MRW}, we ask for two more structural requirements related to our collection of quasimetric $\rho$-balls $\{B(x,r)\}_{r>0;x\in\Omega}$.  The first of these is the existence of appropriate sequences of
Lipschitz cutoff functions (called ``accumulating sequences of Lipschitz
cutoff functions'' in \cite{SW1}).  Specifically, for the function $R_2$ related to the Poincar\'e-Sobolev estimate \eqref{sobolev}, we assume there
are positive constants $s^*, C_{s^*}, \tau$ and $N$, with
$p\sigma^\prime<s^* \le \infty$ and $\tau <1$, such that for
every $\rho$-ball $B(y,r)$ with $0< r< R_2(y)$, there
is a collection of Lipschitz functions $\{\eta_j\}_{j=1}^\infty$ satisfying
\begin{equation}\label{cutoff}
\begin{cases}
\text{supp}\,\eta_1\subset B(y,r)\\
0\leq\eta_j\leq 1\quad \text{ for all }j\\
B(y,\tau r)\subset\{x\in B(y,r):\eta_j(x)=1\}\quad\text{for all $j$}\\
\text{supp}\,\eta_{j+1}\subset\{x\in
B(y,r):\eta_j(x)=1\}\quad\text{for all $j$}\\
\displaystyle{\left(\frac{1}{|B(y,r)|}\int_{B(y,r)}\big|\sqrt{Q}\nabla\eta_j
\big|^{s^*} dx\right)^{1/s^*}\leq
C_{s^*}\frac{N^j}{r}}\quad\text{for all $j$.}
\end{cases}
\end{equation}
This condition is slightly weaker than the corresponding one in
\cite{SW1}; see \cite[p. 149]{MRW} for a fuller discussion. We note
that since $s^*>p\sigma'$, there is a number $s'>\sigma'$ such that
$s^* = ps'$.  The exponent $s=\frac{s^*}{s^*-p}$ dual to $s^\prime$
satisfies $1 \le s <\sigma$ and plays an important role in our
results.

\begin{rem}
As already mentioned in Remark \ref{rem1-KMR}, conditions
\eqref{sobolev} and \eqref{cutoff} imply the validity of the local
doubling condition \eqref{doubling1} for some positive exponent
$d_0$ (see \cite{KRM}). It is important to note that the smaller the
exponent $d_0$ in \eqref{doubling1} can be chosen, the weaker the
required assumptions of local integrability on the coefficients
$b,c,d,e,f,g,h$ in \eqref{struct} will be in the theorems to
follow. See the statements of Proposition \ref{pro1}, of Theorems
\ref{harnackmother}, \ref{holder1}, \ref{harnackmother.2},
\ref{holder1.2}, \ref{harnackmother.3}, \ref{holder1.3} and of
Corollaries \ref{holder2}, \ref{holder2.2}, \ref{holder2.3}.
\end{rem}

Our last requirement is that the following pair of inequalities hold
simultaneously:  There exists $t\in[1,\infty]$ such that for every
$\rho$-ball $B(y,r)$ with $0<r<R_2(y)$, there is a constant
$C_3=C_3(B(y,r))>0$ such that
\begin{equation}\label{E3}
\left(\int_{B(y,r)}|\sqrt{Q}\nabla\eta|^{pt}\,dx\right)^{1/pt}
<\infty \quad\mbox{and}
\end{equation}

\begin{equation}\label{3.6-0}
\bigg(\int_{B(y,r)} | f|^{pt'}dx\bigg)^{1/pt'} \leq
C_3\,||f||_{W^{1,p}_Q(\Omega)} = C_3\bigg(\int_{\Omega}|\sqrt{Q}
\nabla f|^pdx +\int_{\Omega}|f|^pdx\bigg)^{1/p}
\end{equation}
for all $\eta\in \{\eta_j\}$, $\{\eta_j\}$ as in (\ref{cutoff}), and
all $f\in Lip_{loc}(\Omega)$.  As usual, $t'$ denotes the dual
exponent of $t$.  In case $t$ or $t'$ is
infinite, we simply replace the relevant term in (\ref{E3}) or
(\ref{3.6-0}) by an essential supremum.
\begin{rem}  These inequalities are used in \cite{MRW} to derive a
product rule for elements of $W^{1,p}_Q(\Omega)$. They also imply that
functions in $W^{1,p}_Q(\Omega)$,
which are generally not compactly supported, have sufficiently high
local integrability in case the Sobolev inequality (\ref{sobolev})
holds only for compactly supported Lipschitz functions. See
\cite[Section 2, p. 162]{MRW} for these results. It is useful to note
that \eqref{E3} is automatically satisfied for every $t$
with $1\leq t\leq s^*/p$ by \eqref{cutoff}. However, \eqref{E3} may
also hold for larger values of $t$ independently of \eqref{cutoff}.
See \cite[p. 150]{MRW} for details.  In fact, if \eqref{E3} holds with
$t=\infty$ then \eqref{3.6-0} (with $t'=1$) is trivial due to the form
of the $W^{1,p}_Q(\Omega)$ norm (\ref{Wnorm}).
\end{rem}

In order to simplify notation when combining hypotheses, we fix a
single function $r_1:\Omega\ra(0,\infty)$ satisfying
\bea \label{rone}
r_1(x) &\leq& \min\{R_0(x),R_1(x), R_2(x),1\},\quad x\in \Omega,
\eea
where $R_0$ is as described below (\ref{Cond1}), $R_1$ is as in
Definition \ref{homspace} and $R_2$ is as in (\ref{sobolev}), (\ref{poincare}), \eqref{cutoff}, \eqref{E3}, and \eqref{3.6-0}.

\section{Harnack's Inequality}

We begin this section by recalling some notations of \cite{MRW}.
Given a measurable set $E$ and a measurable function $f$ on $E$, we
write \bea \label{n1}|| f ||_{p,E;\overline{dx}} &=&
\left(\frac{1}{|E|}\int_E |f(x)|^p dx\right)^\frac{1}{p}
= \left(\fint_E |f(x)|^p dx\right)^\frac{1}{p},\text{ and}\\
\label{n2}|| f ||_{p,E;dx} &=& \left(\int_E |f(x)|^p dx\right)^\frac{1}{p}.
\eea
In some cases when context is clear, the set $E$ may be dropped from
the left hand side in (\ref{n1}) and (\ref{n2}).

Given a function $u$ and constants $k,
\epsilon_1,\epsilon_2,\epsilon_3$ with $k>0$ and
$\epsilon_1,\epsilon_2,\epsilon_3\in(0,1]$, we denote
\begin{eqnarray}\label{15}
\begin{array}{rclcrcl}
\bar{u}&=&|u|+k,&\qquad&\bar{b}&=&b+k^{1-p}e,\\
\bar{h}&=&h+k^{-p}g,&\qquad&\bar{d}&=&d+k^{1-p}f.
\end{array}
\end{eqnarray}
Here, $b, c, d, e, f, g, h$ denote the coefficients in
\eqref{struct}. Furthermore, for each $\rho$-ball $B(y,r)$, define
\begin{eqnarray} \label{Zetabar}
&&{\bar Z} = {\bar Z}(B(y,r),\bar{u}) =
1+r^{p-1}\avnormb{\bar{b}}{{p' \sigma'}}{B(y,r)}\\&&\nonumber\,\, +
\left(r^p \avnormb{c^\frac{p}{p+1-\psi}
\bar{u}^\frac{p(\psi-p)}{p+1-\psi}}{\frac{p
\sigma'}{p-\epsilon_1}}{B(y,r)}\right)^{\frac{1}{\epsilon_1}}\!
 +\left(r^p\avnormb{\bar{h}}{\frac{p\sigma'}{p-\epsilon_2}}{B(y,r)}
\right)^\frac{1}{\epsilon_2}\!+ \left(r^p
\avnormb{\bar{d}}{\frac{p\sigma'}{p-\epsilon_3}}{B(y,r)}
\right)^\frac{1}{\epsilon_3}\!\!,
\end{eqnarray}
where the exponents $p,\psi, \sigma$ are as usual; see (\ref{struct})
and (\ref{sobolev}).  It is important to note that $\bar{Z}$ is not
monotone in its first argument due to the normalized norms appearing
in its definition. However, if ${\bar Z}(B,{\bar u})<\infty$, then
${\bar Z}(B',{\bar u})<\infty$ whenever $B'\subset B=B(y,r)$ with
$r<R_0(y)$.\\

\subsection{Standing Assumptions}

In order to state our main results efficiently, we list here several
standing assumptions to remain in effect for the rest of this paper.
As above, $\Omega$ will always denote a bounded domain in
$\mathbb{R}^n$, $\rho$ denotes a quasimetric on $\Omega$, and $Q(x)$
denotes a measurable symmetric nonnegative definite matrix defined in
$\Omega$.  We always assume the triple $(\Omega,\rho,dx)$ defines a
local homogeneous space in the sense of Definition \ref{homspace}.
Note that this ensures that the local doubling condition
\eqref{doubling1} is satisfied.  We also assume the validity of the
local Sobolev and Poincar\'e inequalities \eqref{sobolev} and
\eqref{poincare} and the existence of accumulating sequences of
Lipschitz cutoff functions satisfying \eqref{cutoff} for a fixed
$\tau\in(0,1)$ and $s^*>p\sigma'$.  Here $\sigma'$ denotes the dual
exponent to the Sobolev gain factor $\sigma$ of \eqref{sobolev}.
Lastly, we assume that each of \eqref{E3} and \eqref{3.6-0} holds for
some $t\in [1,\infty]$.  We can now state our core Harnack
result. Under certain conditions, it will spawn other versions of
Harnack's inequality that will lead to continuity of weak solutions.

\subsection{Main Results}
\begin{pro}\label{harnackmain} Let $1<p<\infty$ and $|Q(x)|\in
L^{p/2}_{\text{loc}}(\Omega)$. Assume that the functions $A,B$ of
\eqref{eqdiff} satisfy \eqref{struct} with
\begin{equation}\label{ranges1}
\gamma=\delta=p, \quad \psi \in [p, p+1-\sigma^{-1}).
\end{equation}
Fix $y\in \Omega$ and suppose there is a function $r_1(x)$ as in
\eqref{rone} that satisfies a local uniformity condition in
$B(y,r_1(y))$ with constant $A_*=A_*(y,r_1(y))$; see \eqref{unifr}. Let
\begin{equation}\label{5}
C_* =\frac{128\mathfrak{b}\kappa^{10}(\gamma^*)^8}{\tau
A_*^3\min\{A_*^2,(8\kappa^5)^{-1}\}},
\end{equation}
where $\mathfrak{b}$ is from (\ref{poincare}). For $x_0\in
B\big(y,\frac{\tau}{5\kappa}r_1(y)\big)$ and $r\in\big(0,\frac{\tau
A_*}{5\kappa C_*}r_1(y)\big)$, define $$\mathfrak{E} = \{(x,l)\; :\;
B(x,l)\Subset B(x_0,C_* r)\textrm{ and }0<l\leq C_*r\}.$$  Let
$(u,\nabla u)\in W^{1,p}(\Omega)$ be a weak solution of
\eqref{eqdiff}. Assume that
$\epsilon_1,\epsilon_2,\epsilon_3\in(0,1]$ and $k\geq0$ are such
that \bea\label{zbarcond} \disp\sup_{(x,l)\in \mathfrak{E}}
\bar{Z}\big( B(x,l),\bar{u} \big) = M <\infty, \eea where $\bar{Z}$
is defined by \eqref{Zetabar} and $\bar{u}=|u|+k$. If $u\geq 0$ in
$B(x_0,C_*r)$, then the Harnack inequality \bea\label{Harngoal}
\disp\tesup_{z\in B(x_0,\tau r)} \bar{u}(z) &\leq&
C_4\big[C_5\bar{Z}(B(x_0,r),\bar{u})\big]^{C_6M}\disp\teinf_{z\in
B(x_0,\tau r)}\bar{u}(z) \eea holds with
\begin{itemize}
\item[i)] $C_4$ depending on $p,\sigma,s,\epsilon_1,\epsilon_2,\epsilon_3$,
\item[ii)] $C_5$ depending on $a,p,\sigma,s,\epsilon_1,\epsilon_2,\epsilon_3$, $C_1$
in \eqref{sobolev}, $N,C_{s^*}$ in \eqref{cutoff}, on the
pseudometric $\rho$,
\item[iii)] $C_6$ depending on $a,p,\sigma,s,\epsilon_1,\epsilon_2,\epsilon_3$, $C_2,\mathfrak{b}$
in \eqref{poincare}, $\tau,N,C_{s^*}$ in \eqref{cutoff}, $C_0,d_0$
in \eqref{doubling1} and on the pseudometric $\rho$.
\end{itemize}
\end{pro}

\begin{rem}
Since under the hypotheses of Proposition \ref{harnackmain} one has
$(x_0,r)\in\mathfrak{E}$, we obtain \bea\label{Harngoal0}
\disp\tesup_{z\in B(x_0,\tau r)} \bar{u}(z) &\leq&
C_4\big[C_5M\big]^{C_6M}\disp\teinf_{z\in B(x_0,\tau r)}\bar{u}(z)
\eea with $C_4,C_5,C_6$ independent of $(u,\nabla
u),b,c,d,e,f,g,h,k,y,x_0,r,M.$
\end{rem}

A proof of Proposition \ref{harnackmain} is given in \S 6.  The next
proposition provides explicit integrability conditions on structural
coefficients and choices of $\ep_1,\ep_2,\ep_3$ and $k$ that
ensure condition \eqref{zbarcond} is satisfied.  It also provides a
decay condition on $k$ essential for proving H\"older continuity of
weak solutions to \eqref{eqdiff}; see Theorem \ref{holder1} and its proof.


\begin{pro}\label{pro1}
Let $1<p<\infty$, $|Q|\in L_{loc}^{p/2}(\Omega)$, and $\psi \in
[p,p+1-\sigma^{-1})$. Let $(u, \nabla u)\in W_Q^{1,p}(\Omega)$,
fix $y\in \Omega$ and suppose there is a function $r_1(x)$ as
in \eqref{rone} which satisfies a local uniformity condition in
$B(y,r_1(y))$ with constant $A_*=A_*(y,r_1(y))$. Let $C_*$
be defined as in \eqref{5} and assume that
\begin{itemize}
\item[i)] $b,e\in L^\mathcal{B}_\text{loc}(\Omega)$ with
$\mathcal{B}\geq\max\big\{p^\prime\sigma^\prime,\frac{d_0}{p-1}\big\}$;
\item[ii)] $h,g\in L^\mathcal{H}_\text{loc}(\Omega)$ with
$\mathcal{H}\geq\frac{d_0}{p},\,\mathcal{H}>\sigma^\prime$;
\item[iii)] $d,f\in L^\mathcal{D}_\text{loc}(\Omega)$ with
$\mathcal{D}\geq\frac{d_0}{p},\,\mathcal{D}>\sigma^\prime$;
\item[iv)] $c\in L^\mathcal{C}_\text{loc}(\Omega)$ with
$\mathcal{C}\geq\frac{d_0p\sigma}{(p+1-\psi)(p\sigma+d_0)-d_0}>0$ and
$\mathcal{C}>\frac{p\sigma}{\sigma(p+1-\psi)-1}$.
\end{itemize}
For every $x_0\in B\big(y,\frac{\tau}{5\kappa}r_1(y)\big)$ and
$r\in\big(0,\frac{\tau A_*}{5\kappa C_*}r_1(y)\big)$, define
\begin{eqnarray*}
k&=&k(x_0,r)\,=\,k\big(B(x_0,C_*r)\big)\\
&=&\left[(C_*r)^{p-1}\|e\|_{\mathcal{B},B(x_0,C_*r);\overline{dx}}\right]^\frac{
1 } { p-1 }+
\left[(C_*r)^{p}\|f\|_{\mathcal{D},B(x_0,C_*r);\overline{dx}}\right]^\frac{
1 } { p-1 } +
\left[(C_*r)^{p}\|g\|_{\mathcal{H},B(x_0,C_*r);\overline{dx}}\right]^\frac{
1 } { p } .
\end{eqnarray*}
Then
\bea\label{14} k(x_0,r)\leq\Lambda r^\lambda,
\eea
where $\lambda,\Lambda$ are nonnegative numbers independent of $x_0,r$
of the form
\begin{eqnarray*}
\lambda&=&\min\Big\{1-\frac{d_0}{(p-1)\mathcal{B}},\frac{1}{p-1}\Big(p-
\frac{d_0}{\mathcal{D}}\Big), 1-\frac{d_0}{p\mathcal{H}}\Big\},\text{
and}\\
\Lambda&=&C_0^\frac{1}{(p-1)
\mathcal{B}}C_*^{1-\frac{d_0}{(p-1)\mathcal{B}}}r_1(y)^{1-\lambda}
\|e\|_{\mathcal{B},B(y,r_1(y));\overline{dx}}^\frac{1}{p-1}\\
&&+\;C_0^\frac{1}{(p-1)\mathcal{D}}C_*^{\frac{p}{p-1}-\frac{d_0}{(p-1)\mathcal{D}}}
r_1(y)^{\frac{p}{p-1}-\lambda}\|f\|_{\mathcal{D},B(y,r_1(y))
;\overline{dx}}^\frac{1}{p-1}\\
&&+\; C_0^\frac{1}{p\mathcal{H}}C_*^{1-\frac{d_0}{p\mathcal{H}}}r_1(y)^{1-\lambda}
\|g\|_{\mathcal{H},B(y,r_1(y));\overline{dx}}^\frac{1}{p}.
\end{eqnarray*}
Moreover, with $\epsilon_1$, $\epsilon_2$ and $\epsilon_3$ defined
by
\begin{eqnarray}\epsilon_1=\min\Big\{1,\frac{p\sigma(p+1-\psi)-p-(p^2\sigma/
\mathcal{C})}{(\sigma-1)(p+1-\psi)}\Big\},\,\,\,\,
\epsilon_2=\min\Big\{1,p-\frac{p\sigma^\prime}{\mathcal{H}}
\Big\},\,\,\,\,
\epsilon_3=\min\Big\{1,p-\frac{p\sigma'}{\mathcal{D}}\Big\},
\nonumber
\end{eqnarray}
\eqref{zbarcond} is satisfied with
\begin{eqnarray*}
M\!&\!=\!&\!1+C_0^\frac{1}{\mathcal{B}}\left[1+r_1(y)^{p-1}\|b\|_{\mathcal{B},B(y,r_1(y));\overline{dx}}\right]+C_0^\frac{1}{\ep_2\mathcal{H}}\left[1+r_1(y)^{p}\|h\|_{\mathcal{H},B(y,r_1(y));\overline{dx}}
\right]^\frac{1}{\epsilon_2}\\
&&\;\;+\;C_0^\frac{1}{\ep_3\mathcal{D}}\left[1+r_1(y)^{p}\|d\|_{\mathcal{D},B(y,r_1(y));\overline{dx}}
\right]^\frac{1}{\epsilon_3}\\
&&\;\;+C_0^\frac{\psi-p+(p\sigma/\mathcal{C})}{
\epsilon_1(p+1-\psi)\sigma}\Bigg[r_1(y)^{p}\|c\|_{\mathcal{C},
B(y,r_1(y));\overline{dx}}^\frac{p}{p+1-\psi}
\left(\|u\|_{p\sigma,B(y,r_1(y));\overline{dx}} +\Lambda
r_1(y)^\lambda\right)^\frac{p(\psi-p)}{p+1-\psi}\Bigg]^\frac{1}{\epsilon_1},
\end{eqnarray*}
where $C_0$ is as in \eqref{doubling1}.
\end{pro}
Proposition \ref{pro1} is proved in the appendix.

\begin{rem}\label{remA}\begin{enumerate}
\item In part (iv), the assumption that
$\frac{d_0p\sigma}{(p+1-\psi)(p\sigma+d_0)-d_0}>0$ follows from the
condition $\psi\in [p,p+1-\sigma^{-1})$ provided $d_0 \le p\sigma'$;
also, in the classical Euclidean situation, the condition $d_0 \le
p\sigma'$ is true with equality.  If it is the case that $d_0>p\sigma'$, this condition further restricts $\psi\in[p, p+1-\frac{d_0}{d_0+p\sigma}) \subsetneq [p+1-\sigma^{-1})$.
\item The constants $\lambda,\Lambda,M$ in Proposition \ref{pro1} are
independent of $x_0,r$. Moreover $\lambda$ is independent of $y$. The
constant $M$ depends on $u$ only through
$\|u\|_{p\sigma,B(y,r_1(y));dx}$, and it is independent of $u$ when
$\psi=p$.
\item The strict inequalities in $(ii)$, $(iii)$ and $(iv)$
guarantee that $\epsilon_1$, $\epsilon_2$, $\epsilon_3>0$; moreover
$\lambda>0$ if all the inequalities in $(i)$, $(ii)$, $(iii)$ are
strict.
\end{enumerate}
\end{rem}

Combining Propositions \ref{harnackmain} and \ref{pro1} we obtain the
following theorem.


\begin{thm}\label{harnackmother}(Harnack's Inequality, when $\gamma=\delta=p$ and $\psi\geq p$)
Let $1<p<\infty$ and $|Q|\in L_{loc}^{p/2}(\Omega)$.  Let $A,B$ be
functions satisfying \eqref{struct} with $\gamma,\delta, \psi$
restricted to
\begin{equation*}
\gamma=\delta=p, \quad \psi \in [p, p+1-\sigma^{-1}).
\end{equation*}
Fix $y\in \Omega$ and suppose there is a function $r_1(x)$ as in
\eqref{rone} which satisfies a local uniformity condition in
$B(y,r_1(y))$ with constant $A_*=A_*(y,r_1(y))$.  Let $C_*$ be as in
\eqref{5}, $x_0\in B(y,\frac{\tau}{5\kappa}r_1(y))$ and
$r\in(0,\frac{\tau A_*}{5\kappa C_*}r_1(y))$.  Assume that the
structural functions $b,c,d,e,f,g,h$ of \eqref{struct} and
$\epsilon_1,\epsilon_2,\epsilon_3$ and $k=k(x_0,r)$ are as in
Proposition \ref{pro1}. If $(u, \nabla u)\in W_Q^{1,p}(\Omega)$ is a
weak solution of \eqref{eqdiff} in $\Omega$ and $u\geq 0$ in
$B(x_0,C_* r)$, then
\begin{eqnarray}\label{harnack}
\disp\tesup_{z\in B(x_0,\tau r)} \Big(u(z)+k(x_0,r)\Big) &\leq&
C\disp\teinf_{z\in B(x_0,\tau r)}\Big(u(z) + k(x_0,r)\Big),
\end{eqnarray}
with $C=C_4(C_5M)^{C_6M}$, $M$ as in Proposition \ref{pro1} and
$C_4,C_5,C_6$ as in Proposition \ref{harnackmain} with
$\epsilon_1,\epsilon_2,\epsilon_3$ given in Proposition \ref{pro1}.
The constant $C$ depends on $||u||_{p\sigma,B(y,r_1(y));dx}$ only
when $\psi>p$ and only through $M$.
\end{thm}


The proof of Theorem \ref{harnackmother} follows by simply combining
Propositions \ref{harnackmain}, \ref{pro1} and is left to the
reader. Theorem \ref{harnackmother} will allow us to prove H\"older
continuity of weak solutions to \eqref{eqdiff}.  First we recall the
notions of H\"older continuity that we will use.

\begin{defn} Let $w:\Omega\ra \mathbb{R}$ and $S \subset\Omega$. We
say that $w$ is:
\begin{enumerate}
\item {\bf essentially H\"older continuous} with respect to $\rho$ in
$S$ if there are positive constants $C,\mu$ such that
\bea \label{holderconstants2}\tesup_{z,x\in S}
\frac{|w(z)-w(x)|}{\rho(z,x)^\mu} \leq C.
\eea
\item {\bf essentially locally H\"older continuous} with respect to
$\rho$ in $S$ if for every compact set $K\subset S$, there are
positive constants $C,\mu$ such that
\bea \label{holderconstants1}\tesup_{z,x\in K}
\frac{|w(z)-w(x)|}{\rho(z,x)^\mu} \leq C.
\eea
\end{enumerate}
\end{defn}

In these definitions, the notion of H\"older continuity of a function
is relative to the quasimetric $\rho$.  Classical H\"older continuity
with respect to the usual Euclidean metric then follows by imposing
a Fefferman-Phong containment condition on the family of quasimetric
$\rho$-balls. Recall that a Fefferman-Phong condition holds if there
are positive constants $C,\ep$ such $D(x,r) \subset B(x,Cr^\ep)$ for
$x\in \Omega$ and $r>0$ sufficiently small (in terms of $x$).  Several
references impose this condition for such a purpose; see \cite{FP} and
\cite{SW1} for further discussion.\\

Our study of H\"older continuity of weak solutions begins with the
case when the exponents $\gamma, \delta,\psi$ are restricted as in
\eqref{ranges1}.

\begin{thm}\label{holder1}
(H\"{o}lder continuity, when $\gamma=\delta=p$ and $\psi\geq p$) Let
$1<p<\infty$ and $|Q|\in L_{loc}^{p/2}(\Omega)$. Let $(u,\nabla u)$
be a weak solution of \eqref{eqdiff} in $\Omega$ where the functions
$A(x,z,\xi)$ and $B(x,z,\xi)$ satisfy \eqref{struct} with
$\gamma,\delta,\psi$ as in \eqref{ranges1}. Assume that the
coefficient functions of \eqref{struct} satisfy conditions (i)-(iv)
of Proposition \ref{pro1} with strict inequality. Let $y\in\Omega$
and suppose there is a function $r_1(x)$ as in \eqref{rone} which
satisfies a local uniformity condition in $B(y,r_1(y))$ with
constant $A_*=A_*(y,r_1(y))$.  Then $u$ is essentially H\"older
continuous with respect to $\rho$ in
$B(y,\frac{\tau^2}{5\kappa}r_1(y))$. The constants $C$ and $\mu$ in
\eqref{holderconstants1} depend on $y,r_1(y),A_*$, $\kappa$ as in
\eqref{tri}, the Harnack constant $C_4(C_5M)^{C_6M}$ which appears
in Theorem \ref{harnackmother}, $\lambda$ as in Proposition
\ref{pro1}; $C$ depends also on $\|u\|_{p\sigma,B(y,r_1(y));dx}$.
\end{thm}

\begin{rem}
We explicitly note that $\mu$ in the previous Theorem depends on
$\|u\|_{p\sigma,B(y,r_1(y));dx}$ only through $M$, and thus it
depends on $u$ itself only if $\psi>p$.
\end{rem}
Theorem \ref{holder1} is proved in \S 7.  The next result gives
sufficient conditions for essential local H\"older continuity of
solutions in $\Omega$.\\

\begin{cor} \label{holder2}
Let $1<p<\infty$ and $|Q|\in L_{loc}^{p/2}(\Omega)$, and suppose
\eqref{struct} holds with $\gamma,\delta,\psi$ as in \eqref{ranges1}.
Assume also that the coefficient functions of \eqref{struct} satisfy
conditions (i)-(iv) of Proposition \ref{pro1} with strict inequality.
Let $r_1:\Omega\ra (0,\infty)$ be a function satisfying \eqref{rone}
with the property that given any compact $K\subset \Omega$ there is a
positive constant $s_0$ such that $s_0\leq r_1(y)\leq 1$ for
every $y\in K$.  Then if $(u, \nabla u)\in W_Q^{1,p}(\Omega)$ is a
weak solution of \eqref{eqdiff} in $\Omega$, $u$ is essentially
locally H\"older continuous with respect to $\rho$ in $\Omega$.
\end{cor}

\noindent A brief proof of Corollary \ref{holder2} can be found in \S 8.

\subsection{Some consequences}\label{subsec3.3}

The following results are concerned with some of the possible cases
when the exponents $\gamma$, $\delta$, $\psi$ are allowed to vary in
the ranges given in \eqref{ranges}.

In particular, Theorems \ref{harnackmother.2} and \ref{holder1.2}
and Corollary \ref{holder2.2} are devoted to the case when
$\gamma,\delta,\psi< p$. We consider the case when
$\gamma,\delta,\psi> p$ and satisfy \eqref{ranges} in Theorems
\ref{harnackmother.3} and \ref{holder1.3} and in Corollary
\ref{holder2.3}. See \S 9 for their proofs.

Of course, similar results can be obtained for other choices of $\gamma, \delta, \psi$ in the ranges given in \eqref{ranges} but we
won't list them here. Such results can all be derived from Theorems
\ref{harnackmother}, \ref{holder1} and Corollary \ref{holder2}. We
leave the details to the interested reader.


\begin{thm}\label{harnackmother.2}(Harnack's Inequality, when $\gamma,\delta,\psi<p$)
Let $1<p<\infty$ and $|Q|\in L_{loc}^{p/2}(\Omega)$.  Let $A,B$ be
functions satisfying \eqref{struct} with $\gamma,\delta, \psi$
restricted to
\begin{equation}\label{ranges3}
\gamma,\delta,\psi \in (1, p).
\end{equation}
Fix $y\in \Omega$ and suppose there is a function $r_1(x)$ as in
\eqref{rone} which satisfies a local uniformity condition in
$B(y,r_1(y))$ with constant $A_*=A_*(y,r_1(y))$.  Let $C_*$ be as in
\eqref{5}, $x_0\in B(y,\frac{\tau}{5\kappa}r_1(y))$ and
$r\in(0,\frac{\tau A_*}{5\kappa C_*}r_1(y))$.  Assume that the
structural functions $b,d,e,f,g,h$ of \eqref{struct} satisfy
conditions $(i)$, $(ii)$ and $(iii)$ in Proposition \ref{pro1}, that
$c\in L^\mathcal{C}_\text{loc}(\Omega)$ with $\mathcal{C}\geq d_0$
and $\mathcal{C}>p\sigma'$. Let $\epsilon_2,\epsilon_3,\lambda$ be
as in Proposition \ref{pro1} and define
\begin{eqnarray*}
\epsilon_1&=&\min\Big\{1,p-\frac{p^2\sigma'}{\mathcal{C}}\Big\},\\
k_1&=&k_1(x_0,r)\,=\,k_1\big(B(x_0,C_*r)\big)\\
&=&\left[(C_*r)^{p-1}\|b+e\|_{\mathcal{B},B(x_0,C_*r);\overline{dx}}\right]^\frac{
1 } { p-1 }+
\left[(C_*r)^{p}\|c+d+f\|_{\mathcal{D},B(x_0,C_*r);\overline{dx}}\right]^\frac{
1 } { p-1 } \\
&&+\left[(C_*r)^{p}\|g+h\|_{\mathcal{H},B(x_0,C_*r);\overline{dx}}\right]^\frac{
1 } { p },\\
\Lambda_1\!&\!=\!&\!C_0^\frac{1}{(p-1)
\mathcal{B}}C_*^{1-\frac{d_0}{(p-1)\mathcal{B}}}r_1(y)^{1-\lambda}
\|b+e\|_{\mathcal{B},B(y,r_1(y));\overline{dx}}^\frac{1}{p-1}\\
&&+ \;C_0^\frac{1}{(p-1)\mathcal{D}}C_*^{\frac{p\mathcal{D}-d_0}{(p-1)\mathcal{D}}}
r_1(y)^{\frac{p}{p-1}-\lambda}\|c+d+f\|_{\mathcal{D},B(y,r_1(y))
;\overline{dx}}^\frac{1}{p-1}\\
&&+\; C_0^\frac{1}{p\mathcal{H}}C_*^{1-\frac{d_0}{p\mathcal{H}}}r_1(y)^{1-\lambda}
\|g+h\|_{\mathcal{H},B(y,r_1(y));\overline{dx}}^\frac{1}{p},\\
M_1\!&\!=\!&\!1+ C_0^\frac{1}{\mathcal{B}}\left[1+r_1(y)^{p-1}\|b\|_{\mathcal{B},B(y,r_1(y));\overline{dx}}\right]+ C_0^\frac{1}{\ep_2\mathcal{H}}\left[1+r_1(y)^{p}\|h\|_{\mathcal{H},B(y,r_1(y));\overline{dx}}
\right]^\frac{1}{\e_2}\\
&&\;\;+\;C_0^\frac{1}{\ep_3\mathcal{D}}\left[1+r_1(y)^{p}\|d\|_{\mathcal{D},B(y,r_1(y));\overline{dx}}
\right]^\frac{1}{\e_3}+C_0^\frac{p}{\ep_1\mathcal{C}}\left[r_1(y)^{p}\|c\|_{\mathcal{C},
B(y,r_1(y));\overline{dx}}^p\right]^\frac{1}{\e_1}.
\end{eqnarray*}
Then
$$k_1(x_0,r)\leq\Lambda_1 r^\lambda$$ and, if $(u, \nabla u)\in W_Q^{1,p}(\Omega)$ is a weak solution of
\eqref{eqdiff} in $\Omega$ such that $u\geq 0$ in $B(x_0,C_* r)$,
\begin{eqnarray}\label{harnack2}
\disp\tesup_{z\in B(x_0,\tau r)} \Big(u(z)+k_1(x_0,r)\Big) &\leq&
C\disp\teinf_{z\in B(x_0,\tau r)}\Big(u(z) + k_1(x_0,r)\Big),
\end{eqnarray}
where $C=C_4(C_5M_1)^{C_6M_1}$, with $C_4,C_5,C_6$ as in Proposition
\ref{harnackmain}.
\end{thm}

\begin{thm}\label{holder1.2}
(H\"{o}lder continuity, when $\gamma,\delta,\psi<p$) Let
$1<p<\infty$ and $|Q|\in L_{loc}^{p/2}(\Omega)$. Let $(u,\nabla u)$
be a weak solution of \eqref{eqdiff} in $\Omega$ where the functions
$A(x,z,\xi)$ and $B(x,z,\xi)$ satisfy \eqref{struct} with
$\gamma,\delta,\psi$ as in \eqref{ranges3}. Assume that the
coefficient functions of \eqref{struct} satisfy the same conditions
as in Theorem \ref{harnackmother.2} with strict inequality. Let
$y\in\Omega$ and suppose there is a function $r_1(x)$ as in
\eqref{rone} which satisfies a local uniformity condition in
$B(y,r_1(y))$ with constant $A_*=A_*(y,r_1(y))$.  Then $u$ is
essentially H\"older continuous with respect to $\rho$ in
$B(y,\frac{\tau^2}{5\kappa}r_1(y))$. The constants $C$ and $\mu$ in
\eqref{holderconstants1} depend on $y,r_1(y),A_*$, $\kappa$ as in
\eqref{tri}, the Harnack constant $C_4(C_5M_1)^{C_6M_1}$ which
appears in Theorem \ref{harnackmother.2}, $\lambda$ as in
Proposition \ref{pro1}; $C$ depends also on
$\|u\|_{p\sigma,B(y,r_1(y));dx}$, while $\mu$ is independent of $(u,\nabla u)$.
\end{thm}

\begin{cor} \label{holder2.2}
Let $1<p<\infty$ and $|Q|\in L_{loc}^{p/2}(\Omega)$, and suppose
\eqref{struct} holds with $\gamma,\delta,\psi$ as in
\eqref{ranges3}. Assume also that the coefficient functions of
\eqref{struct} satisfy the same conditions as in Theorem
\ref{harnackmother.2} with strict inequality. Let $r_1:\Omega\ra
(0,\infty)$ be a function satisfying \eqref{rone} with the property
that given any compact $K\subset \Omega$ there is a positive
constant $s_0$ such that $s_0\leq r_1(y)\leq 1$ for every
$y\in K$.  Then if $(u, \nabla u)\in W_Q^{1,p}(\Omega)$ is a weak
solution of \eqref{eqdiff} in $\Omega$, $u$ is essentially locally
H\"older continuous with respect to $\rho$ in $\Omega$ with exponent $\mu$ that is independent of the weak solution $(u,\nabla u)$.
\end{cor}


\begin{thm}\label{harnackmother.3}(Harnack's Inequality, when $\gamma,\delta,\psi> p$)
Let $1<p<\infty$ and $|Q|\in L_{loc}^{p/2}(\Omega)$.  Let $A,B$ be
functions satisfying \eqref{struct} with $\gamma,\delta, \psi$
satisfying \eqref{ranges} and restricted to
\begin{equation}\label{ranges4}
\gamma,\delta,\psi>p.
\end{equation}
Fix $y\in \Omega$ and suppose there is a function $r_1(x)$ as in
\eqref{rone} which satisfies a local uniformity condition in
$B(y,r_1(y))$ with constant $A_*=A_*(y,r_1(y))$.  Let $C_*$ be as in
\eqref{5}, $x_0\in B(y,\frac{\tau}{5\kappa}r_1(y))$ and
$r\in(0,\frac{\tau A_*}{5\kappa C_*}r_1(y))$.  Assume that the
structural functions $b,d,e,f,g,h$ of \eqref{struct} satisfy
\begin{itemize}
\item[i)] $b\in L^{\mathcal{B}_0}_\text{loc}(\Omega)$ with
$\mathcal{B}_0\geq\frac{p\sigma}{p\sigma-\sigma-\gamma+1}$ and
$\mathcal{B}_0\geq\frac{d_0p\sigma}{p(p-1)\sigma-d_0(\gamma-p)}>0$;
\item[ii)] $e\in L^{\mathcal{E}}_\text{loc}(\Omega)$ with
$\mathcal{E}\geq\max\big\{p^\prime\sigma^\prime,\frac{d_0}{p-1}\big\}$;
\item[iii)] $h\in L^{\mathcal{H}_0}_\text{loc}(\Omega)$ with
$\mathcal{H}_0>\frac{p\sigma}{p\sigma-\gamma}$ and
$\mathcal{H}_0\geq\frac{d_0p\sigma}{p^2\sigma-d_0(\gamma-p)}>0$;
\item[iv)] $g\in L^\mathcal{G}_\text{loc}(\Omega)$ with
$\mathcal{G}\geq\frac{d_0}{p},\,\mathcal{G}>\sigma^\prime$;
\item[v)] $d\in L^{\mathcal{D}_0}_\text{loc}(\Omega)$ with
$\mathcal{D}_0>\frac{p\sigma}{p\sigma-\delta}$ and
$\mathcal{D}_0\geq\frac{d_0p\sigma}{p^2\sigma-d_0(\delta-p)}>0$;
\item[vi)] $f\in L^{\mathcal{F}}_\text{loc}(\Omega)$ with
$\mathcal{F}\geq\frac{d_0}{p},\,\mathcal{F}>\sigma^\prime$;
\item[vii)] $c\in L^\mathcal{C}_\text{loc}(\Omega)$ with
$\mathcal{C}>\frac{p\sigma}{\sigma(p+1-\psi)-1}$ and
$\mathcal{C}\geq\frac{d_0p\sigma}{(p+1-\psi)(p\sigma+d_0)-d_0}>0$.
\end{itemize}
Define
$$\mathcal{B}=\min\left\{\frac{p\sigma}{\frac{p\sigma}{\mathcal{B}_0}+(\gamma-p)},\mathcal{E}\right\},\,
\mathcal{H}=\min\left\{\frac{p\sigma}{\frac{p\sigma}{\mathcal{H}_0}+(\gamma-p)},\mathcal{G}\right\},\,
\mathcal{D}=\min\left\{\frac{p\sigma}{\frac{p\sigma}{\mathcal{D}_0}+(\delta-p)},\mathcal{F}\right\}.$$
Let $k=k(x_0,r),\epsilon_1,\epsilon_2,\epsilon_3,\lambda,\Lambda$ be
as in Proposition \ref{pro1} and define
\begin{eqnarray}
\nonumber
M_2\!&\!=\!&\!1+C_0^\frac{1}{\mathcal{B}}\Bigg[1+\frac{r_1(y)^{p-1}}{|B\big(y,r_1(y)\big)|^\frac{1}{\mathcal{B}}}\|b\|_{\mathcal{B}_0,B(y,r_1(y));dx}
\|u\|_{p\sigma,B(y,r_1(y));dx}^{\gamma-p}\Bigg]\\
\nonumber &&\;\; +\;C_0^\frac{1}{\ep_2\mathcal{H}}\left[1+\frac{r_1(y)^{p}}{|B\big(y,r_1(y)\big)|^\frac{1}{\mathcal{H}}}\|h\|_{\mathcal{H}_0,B(y,r_1(y));dx}
\|u\|_{p\sigma,B(y,r_1(y));dx}^{\gamma-p}\right]^\frac{1}{\epsilon_2}\\
\label{18}&&\;\;+\;C_0^\frac{1}{\ep_3\mathcal{D}}\left[1+\frac{r_1(y)^{p}}{|B\big(y,
r_1(y)\big)|^\frac{1}{\mathcal{D}}}\|d\|_{\mathcal{D}_0,B(y,r_1(y));dx}
\|u\|_{p\sigma,B(y,r_1(y));dx}^{\delta-p}\right]^\frac{1}{\epsilon_3}
\\
\nonumber&&\;\;+C_0^\frac{\psi-p+(p\sigma/\mathcal{C})}{\epsilon_1
(p+1-\psi)\sigma}\Bigg[r_1(y)^{p}\|c\|_{\mathcal{C},
B(y,r_1(y));\overline{dx}}^\frac{p}{p+1-\psi}\left(\|u\|_{p\sigma,B(y,r_1(y));\overline{dx}}
+\Lambda
r_1(y)^\lambda\right)^\frac{p(\psi-p)}{p+1-\psi}\Bigg]^\frac{1}{\epsilon_1},
\end{eqnarray}
If $(u, \nabla u)\in W_Q^{1,p}(\Omega)$ is a weak solution of
\eqref{eqdiff} in $\Omega$ and $u\geq 0$ in $B(x_0,C_* r)$, then
\begin{eqnarray}\label{harnack3}
\disp\tesup_{z\in B(x_0,\tau r)} \Big(u(z)+k(x_0,r)\Big) &\leq&
C\disp\teinf_{z\in B(x_0,\tau r)}\Big(u(z) + k(x_0,r)\Big),
\end{eqnarray}
where $k=k(x_0,r)$ satisfies \eqref{14} and
$C=C_4(C_5M_2)^{C_6M_2}$, with $C_4,C_5,C_6$ as in Proposition
\ref{harnackmain}.
\end{thm}

\begin{rem}
In parts (i), (iii), (v) and (vii) of the assumptions of Theorem
\ref{harnackmother.3}, the positivity assumptions on
$\frac{d_0p\sigma}{p(p-1)\sigma-d_0(\gamma-p)}$,
$\frac{d_0p\sigma}{p^2\sigma-d_0(\gamma-p)}$,
$\frac{d_0p\sigma}{p^2\sigma-d_0(\delta-p)}$ and
$\frac{d_0p\sigma}{(p+1-\psi)(p\sigma+d_0)-d_0}$ are a consequence
of conditions \eqref{ranges} and \eqref{ranges4} when $d_0 \le
p\sigma'$.  It is also useful to note that $d_0 \le p\sigma'$ is true with equality in the classical Euclidean situation. In case $d_0>p\sigma'$, the positivity conditions of items $(i),(v),$ and $(vii)$ further restrict the ranges of $\gamma, \delta,$ and $\psi$. See also part (1) of
Remark \ref{remA}.
\end{rem}

\begin{thm}\label{holder1.3}
(H\"{o}lder continuity, when $\gamma,\delta,\psi> p$) Let
$1<p<\infty$ and $|Q|\in L_{loc}^{p/2}(\Omega)$. Let $(u,\nabla u)$
be a weak solution of \eqref{eqdiff} in $\Omega$ where the functions
$A(x,z,\xi)$ and $B(x,z,\xi)$ satisfy \eqref{struct} with
$\gamma,\delta,\psi$ as in \eqref{ranges4} and \eqref{ranges}.
Assume that the coefficient functions of \eqref{struct} satisfy the
same conditions as in Theorem \ref{harnackmother.3} with strict
inequality. Let $y\in\Omega$ and suppose there is a function
$r_1(x)$ as in \eqref{rone} which satisfies a local uniformity
condition in $B(y,r_1(y))$ with constant $A_*=A_*(y,r_1(y))$.  Then
$u$ is essentially H\"older continuous with respect to $\rho$ in
$B(y,\frac{\tau^2}{5\kappa}r_1(y))$. The constants $C$ and $\mu$ in
\eqref{holderconstants1} depend on $y,r_1(y),A_*$, $\kappa$ as in
\eqref{tri}, the Harnack constant $C_4(C_5M_2)^{C_6M_2}$ which
appears in Theorem \ref{harnackmother.3}, $\lambda$ as in
Proposition \ref{pro1}; $C$ depends also on
$\|u\|_{p\sigma,B(y,r_1(y));dx}$.
\end{thm}

\begin{cor} \label{holder2.3}
Let $1<p<\infty$ and $|Q|\in L_{loc}^{p/2}(\Omega)$, and suppose
\eqref{struct} holds with $\gamma,\delta,\psi$ as in \eqref{ranges4}
and \eqref{ranges}. Assume also that the coefficient functions of
\eqref{struct} satisfy the same conditions as in Theorem
\ref{harnackmother.3} with strict inequality. Let $r_1:\Omega\ra
(0,\infty)$ be a function satisfying \eqref{rone} with the property
that given any compact $K\subset \Omega$ there is a positive
constant $s_0$ such that $s_0\leq r_1(y)\leq 1$ for every
$y\in K$.  Then if $(u, \nabla u)\in W_Q^{1,p}(\Omega)$ is a weak
solution of \eqref{eqdiff} in $\Omega$, $u$ is essentially locally
H\"older continuous with respect to $\rho$ in $\Omega$.
\end{cor}

We conclude the section with some comments concerning the rate
growth of the Euclidean volume of pseudometric balls $B(x,r)$.

\begin{defn} \label{weakdoub} Let $\Theta\Subset\Omega$ and $r_1:\Omega\ra
(0,\infty)$ be a function satisfying \eqref{rone}. If $q^*$
satisfies $0<q^* <\infty$ and
  there are positive constants $C_7,\alpha$ such that
\begin{equation}\label{D*}
|B(x,r)| \ge C_7 r^{q^*}
\end{equation}
for all $x\in\Theta$ and all $r<\min\{1,\alpha r_1(x)\}$, we will
say that condition weak-$D_{q^*}$ holds on $\Theta$.
\end{defn}

A similar, but slightly stronger, condition called $D_{q^*}$ was
introduced in Definition 1.7 in \cite{MRW} in order to derive some
local boundedness results for weak solutions of equation
\eqref{eqdiff}; see Corollaries 1.8, 1.9 and 1.11  in \cite{MRW}.

Note that by Definition \ref{homspace}, if $(\Omega,\rho, dx)$ is a
local homogeneous space, $\Theta\Subset\Omega$ and $r_1(x)$
satisfies a local uniformity condition in $\Theta$ with constant
$A_*=A_*(\Theta)$ (see \eqref{unifr}), then condition weak-$D_{q^*}$
automatically holds with $q^*=d_0$ on $\Theta$, for some constant
$C_7>0$ and with $\alpha=A_*/2$. See the Appendix for a proof of
this result.

The fact that property \eqref{D*} holds with $q^*=d_0$ for suitable
families of pseudometric balls $B(x,r)$ with small radii is used
repeatedly in the proofs of our results, starting from Proposition
\ref{pro1} (see Steps I and III of the proof in the Appendix) and in
all the theorems and corollaries that follow it.

It is interesting to note that in the proof of Proposition
\ref{pro1}, only condition \eqref{D*} with $q^*=d_0$ is used to
estimate terms involving the structural coefficients $b,c,d,h$,
while the local Doubling Condition \eqref{doubling1} is directly
used to estimate terms involving some local averages of $e,f,g$ (see
Step 6 of the proof in the Appendix).

%




\section{Some Calculus for Degenerate Sobolev Spaces}
\begin{lem}\label{lem1} Let $\Theta\subset\Omega$ be an open set, $(u,\nabla u)\in W^{1,p}_Q(\Omega)$ with $u\in L^\infty(\Theta)$ and let $\e>0$. Let
$$m=\disp\teinf_\Theta\,\,u ,\qquad M=\disp\tesup_\Theta u.$$
Then there exists a sequence
$\{\varphi_j\}_{j\in\N}\subset\text{Lip}_{\text{loc}}(\Omega)\cap
L^\infty(\Omega)$ such that $(\varphi_j,\nabla\varphi_j)\in
W^{1,p}_Q(\Omega)$ and
\begin{itemize}
\item[i)] $(\varphi_j,\nabla\varphi_j)\rightarrow (u,\nabla u)$ in
$W^{1,p}_Q(\Theta)$,
\item[ii)] $\varphi_j(x)\in[m-\e,M+\e]$ for every $x\in\Omega$ and every
$j\in\N$.
\end{itemize}
\end{lem}

\textbf{Proof:} By definition of $W^{1,p}_Q(\Omega)$, there exists a
sequence
$\{\hat\varphi_j\}_{j\in\N}\subset\text{Lip}_\text{loc}(\Omega)$
such that
$(\hat\varphi_j,\nabla\hat\varphi_j)$ converges to $(u,\nabla u)$ in
$W^{1,p}_Q(\Omega)$. By choosing a subsequence, we may assume that
\begin{equation}\label{2}
\begin{array}{rcl}
\displaystyle\hat\varphi_j\!&\!\rightarrow\!&\!u\qquad\,\,\qquad\,\text{in
}L^p(\Omega),\,\text{ in
}\,W^{1,p}_Q(\Omega)\,\text{ and a.e. in }\,\Omega,\\
\displaystyle\sqrt{Q}\nabla\hat\varphi_j\!&\!\rightarrow\!&\!
\sqrt{Q}\nabla u\qquad\text{in }[L^p(\Omega)]^n,\text{ and a.e. in
}\,\Omega.
\end{array}
\end{equation}
Now for every $j\in\N$ and $x\in\Omega$ define
\begin{equation}\label{1}
\varphi_j(x)=\begin{cases}
\hat\varphi_j(x)\qquad\,\,\,\,\,\text{if }\,\,m-\e\leq\hat\varphi_j(x)\leq M+\e,\\
M+\e\qquad\quad\text{if }\,\,\hat\varphi_j(x)>M+\e,\\
m-\e\qquad\quad\text{if }\,\,\hat\varphi_j(x)<m-\e.
\end{cases}
\end{equation}
This immediately yields that
$\varphi_j\in\text{Lip}_{\text{loc}}(\Omega)$ and that
$$m-\e\leq \varphi_j(x)\leq M+\e\qquad$$ for every $j\in\N$ and $x\in\Omega$.  Then $\varphi_j\in L^\infty(\Omega)$ for every $j\in\N$. From \eqref{1}
it follows that
\begin{equation}\label{3}
\nabla\varphi_j(x)=\begin{cases}
\nabla\hat\varphi_j(x)\qquad\text{if }\,\,m-\e<\hat\varphi_j(x)<M+\e,\\
0\qquad\qquad\,\,\,\,\text{otherwise}
\end{cases}
\end{equation}
for each $j\in\N$ and almost every $x\in\Omega$.  Hence, $|\sqrt{Q}\nabla\varphi_j|\leq|\sqrt{Q}\nabla\hat\varphi_j|$
for every $j\in\N$ and a.e. $x \in \Omega$.  We conclude that $(\varphi_j,\nabla\varphi_j)\in W^{1,p}_Q(\Omega)$ for every
$j\in\N$.

Since $u(x)\in[m,M]$ for a.e. $x\in\Theta$ and
$\hat\varphi_j\rightarrow u$ for a.e. $x\in\Omega$ by  \eqref{2}, we
have that $\hat\varphi_j(x)\in(m-\e,M+\e)$ for a.e. $x\in\Theta$ when
$j$ is large enough. It follows from \eqref{1} that one also has
$\varphi_j(x)=\hat\varphi_j(x)$ pointwise a.e. in $\Theta$ when $j$ is
large enough.  Therefore,
$$\varphi_j\rightarrow u\qquad\text{ a.e. in }\Theta.$$
Moreover, by \eqref{3},  $\nabla\varphi_j=\nabla\hat\varphi_j$ a.e. in
$\Theta$  when $j$ is large enough. Hence, by \eqref{2},
$$\sqrt{Q}\nabla\varphi_j\rightarrow\sqrt{Q}\nabla u\qquad\text{ a.e. in }\Theta.$$

Since
$$|u(x)-\varphi_j(x)|^p\leq|u(x)-\hat\varphi_j(x)|^p\leq2^{p-1}
\big[|u(x)|^p+|\hat\varphi_j(x)|^p\big]$$
for a.e. $x\in\Theta$ and $|u|^p+|\hat\varphi_j|^p\rightarrow2|u|^p$
for a.e. $x\in \Omega$ and in $L^1(\Omega)$ by \eqref{2},
Lebesgue's sequentially dominated convergence theorem implies that
$$\varphi_j\rightarrow u\qquad\text{ in }L^p(\Theta).$$
In a similar way, for a.e. $x\in\Theta$ we have
\begin{eqnarray}
\nonumber|\sqrt{Q}\nabla
u(x)-\sqrt{Q}\nabla\varphi_j(x)|^p&\leq&2^{p-1}\big[|\sqrt{Q}\nabla
u(x)|^p+|\sqrt{Q}\nabla\varphi_j(x)|^p\big]\\
\nonumber &\leq&2^{p-1}\big[|\sqrt{Q}\nabla
u(x)|^p+|\sqrt{Q}\nabla\hat\varphi_j(x)|^p\big].
\end{eqnarray}
Further, we have that $|\sqrt{Q}\nabla
u|^p+|\sqrt{Q}\nabla \hat\varphi_j|^p\rightarrow2|\sqrt{Q}\nabla
u|^p$ a.e. in $\Omega$ and in $L^1(\Omega)$ by \eqref{2}.  Lebesgue's
theorem gives
$$\sqrt{Q}\nabla \varphi_j\rightarrow \sqrt{Q}\nabla u\qquad\text{ in }[L^p(\Theta)]^n.$$
We conclude that $$(\varphi_j,\nabla\varphi_j)\rightarrow
(u,\nabla u)\qquad\text{ in }W^{1,p}_Q(\Theta).$$
\begin{flushright}$\Box$
\end{flushright}

\begin{pro}\label{chain}  Let $\Theta\subset\Omega$ be an open set, $(u,\nabla u)\in W^{1,p}_Q(\Omega)$ with
$u\in L^\infty(\Theta)$ and $$m=\disp\teinf_\Theta u
,\qquad M=\disp\tesup_\Theta u.$$ Let $F\in
C^1\big((m-\e_0,M+\e_0)\big)$ for some $\e_0>0$. Then
$\big(F(u),\nabla(F(u))\big)\in W^{1,p}_Q(\Theta)$ with
\begin{equation}\label{4}
\sqrt{Q}\nabla\big(F(u)\big)=F^\prime(u)\sqrt{Q}\nabla u
\end{equation}
almost everywhere in $\Theta.$
\end{pro}

\textbf{Proof:} The proof is a straightforward adaptation of the
techniques used in the proof of Lemma 4.1 in \cite{MRW}. Fix any
$\e\in(0,\e_0)$ and consider the sequence
$\{\varphi_j\}_{j\geq1}\subset\text{Lip}_{\text{loc}}(\Omega)\cap
L^\infty(\Omega)$ provided by Lemma \ref{lem1}. Notice that
$\varphi_j(x)\in[m-\e,M+\e]$ for every $x\in\Omega$ and every
$j$, that $u(x)\in[m-\e,M+\e]$ for a.e. $x\in\Theta$ and that
$$\sup_{t\in[m-\e,M+\e]}|F(t)|<\infty,\qquad \sup_{t\in[m-\e,M+\e]}|F'(t)|<\infty.$$

Arguing as in Lemma 4.1 in \cite{MRW}, it is easy to see that
$\{F(\varphi_j)\}_{j\in\N}\subset\text{Lip}_{\text{loc}}(\Omega)\cap
L^\infty(\Omega)$ and
$\{\big(F(\varphi_j),\nabla(F(\varphi_j))\big)\}_{j\in\N}$ is a Cauchy
sequence in $W^{1,p}_Q(\Theta)$. Thus, $\{\big(F(\varphi_j),
\nabla(F(\varphi_j))\big)\}_{j\in\N}$ defines an element
$\big(F(u),\nabla(F(u))\big)$ of $W^{1,p}_Q(\Theta)$ that satisfies
\eqref{4}.
\begin{flushright}$\Box$
\end{flushright}

\begin{cor}\label{cor1}
Let $\Theta\Subset\Omega$ be an open set and fix a quasimetric ball $B$ with $B\Subset\Theta$. Suppose that for some
$t\in[1,\infty]$, condition \eqref{3.6-0} holds for $B$ and condition
\eqref{E3} holds for a particular function $\eta\in\text{Lip}_0(B)$.
Let $\theta\geq1$, $(u,\nabla u)\in W^{1,p}_Q(\Omega)$ with $u\in
L^\infty(\Theta)$,
$$m=\disp\teinf_\Theta u ,\qquad
M=\disp\tesup_\Theta u,$$ and $F\in
C^1\big((m-\e_0,M+\e_0)\big)$ for some $\e_0>0$. Then
$\big(\eta^\theta F(u),\nabla(\eta^\theta
F(u))\big)\in(W^{1,p}_{Q})_0(B)$ and
$$\sqrt{Q}\nabla\big(\eta^\theta F(u)\big)=\theta\eta^{\theta-1}F(u)\sqrt{Q}\nabla \eta+\eta^\theta F^\prime(u)\sqrt{Q}\nabla u\quad\text{pointwise a.e. in }\Omega.$$
\end{cor}

\textbf{Proof:} This is a simple consequence of Proposition
\ref{chain} together with Proposition 2.2 in \cite{MRW}.
\begin{flushright}$\Box$
\end{flushright}


\begin{rem}\label{diffineqcomp}
Let $(u,\nabla u)\in W^{1,p}_Q(\Omega)$ be such that $u\geq m$ a.e.
in an open set $\Theta\subset \Omega$, and assume that
$F:(m-\ep_0,\infty)\ra \R$ is $C^1$ with
$\disp\sup_{(m-\ep_0,\infty)}|F'|<\infty$ for some $\ep_0>0$. Then
the conclusions of Proposition \ref{chain} and Corollary \ref{cor1}
still hold. We omit the proofs of these facts as they use ideas
similar to those used in the previous proofs.
\end{rem}


\section[2]{The Inequality of John and Nirenberg}\label{sec5}

This section develops a local version of the inequality of John and
Nirenberg adapted to the class $[cR]BMO(E)$ defined in the next paragraph.
The arguments to follow are adaptations of ones in
\cite{SW1}, where $R(x)$ is a small fixed multiple of
$\mbox{dist}(x,\partial\Omega)$. \\

Let $\Omega$ be an open subset in $\mathbb{R}^n$.
Let $\rho$ be a quasimetric in $\Omega$
and fix $R:\Omega\ra (0,\infty)$.  For each $x\in\Omega$ and
$0<c<\infty$, we say that a $\rho$-ball $B(y,t)$ is a $cR(x)$-ball if
$0<t<cR(x)$, $\overline{B(y,\gamma^*t)}\subset\Omega$, and
$B(y,\gamma^*t)\subset B(x,cR(x))$ where $\gamma^*$ is as in Lemma
\ref{swalemma}.  It is useful to note that if $0<c_1<c_2$ then a
$c_1R(x)$-ball $B$ is also a $c_2R(x)$-ball.   Let $E\subset \Omega$,
$E$ open. A function $f\in L^1_{loc}(\Omega)$ is said to belong to the
class $[cR]BMO(E)$ if
\bea\label{Rbmo} ||f||_{[cR]BMO(E)} = \sup_{x\in E}\sup_{B}
\fracc{1}{|B|}\int_B |f(y) - f_B|dy<\infty,
\eea
where the second supremum is taken over all $cR(x)$-balls $B$.

The main result of this section is as follows.

\begin{pro}\label{JN} Let $(\Omega,\rho,dx)$ be a local
homogeneous space as in Definition \ref{homspace}.  Let
$R:\Omega\ra(0,\infty)$ satisfy $R(x)\leq \min\{R_0(x)/(\gamma^*)^2,
R_1(x)/\gamma^*\}$ for all $x$, where $R_0$ is as above Remark
\ref{rem2.1} and $R_1$ is as in Definition \ref{homspace}.  Fix an
open set $E \subset\Omega$ and assume that $R$ satisfies a local
uniformity condition with respect to $\rho$ in $E$ with constant
$A_*=A_*(R,E)$; see \eqref{unifr}.  Then there are positive
constants $\delta_0=\delta_0(R,E), C_8, C_9, c_\rho$ with $\delta_0
<1$ and $c_\rho>1$ such that for all $x\in E$, all $\delta_0
R(x)$-balls $B$, all $f\in [c_\rho R]BMO(E)$ and all $\alpha >0$,
\bea\label{JN1} |\{y\in B\;:\; |f(y) - f_{B}|>\alpha\}| \leq
C_8e^{\frac{-C_9\alpha}{||f||_{[c_\rho R]BMO(E)}}}|B|. \eea
\end{pro}
\begin{rem} The constants $C_8,C_9$ and $c_\rho$ in Proposition
\ref{JN} depend only on the quasimetric $\rho$, while the dependence
of $\delta_0$ on $E$ occurs only through $A_*$. As the proof of
Proposition \ref{JN} shows, $c_\rho=8(\gamma^*)^2\kappa^5$ and
$\delta_0= A_*^2\min\{A_*^2,(8\kappa^5)^{-1}\}/8(\gamma^*)^3\kappa^5$,
where $\kappa$ is the constant in \eqref{tri} and
$\gamma^*=\kappa+2\kappa^2$ as in Lemma \ref{swalemma}.
\end{rem}
The significance of Proposition \ref{JN} is its consequence for a
special class of $A_2$ weights.  Given $0<c<\infty$ and a set
$E\subset \Omega$, a nonnegative function $w \in L^1_{loc}(\Omega)$ is
said to be a $[cR]A_2(E)$ weight if

\bea\label{RA2} ||w||_{[cR]A_2(E)} = \sup_{x\in
E}\sup_B\Big(\frac{1}{|B|}\int_B wdy\Big)\Big(\frac{1}{|B|}\int_B
w^{-1}dy\Big) < \infty,
\eea
where the second supremum is taken over all $cR(x)$-balls $B$.  We
will use the following corollary of Proposition \ref{JN}
in the proof of Proposition \ref{harnackmain}.

\begin{cor}\label{JNcor} Under the hypotheses of Proposition \ref{JN},
there are constants $C_8,C_9>0$ and $c_\rho>1$ such that for any
open set $E\subset\Omega$, there is a $\delta_0 = \delta_0(R,E)>0$
for which \bea\label{JNcor1} ||e^f||_{[\delta_0  R]A_2(E)} \leq
\left(1+ \frac{C_8||f||_{[c_\rho R]BMO(E)}}{C_9 - ||f||_{[c_\rho
R]BMO(E)}}\right)^2 \eea for every $f\in [c_\rho R]BMO(E)$ with
$||f||_{[c_\rho R]BMO(E)}< C_9$. The constants $\delta_0,\; C_8,\;
C_9,\; c_\rho$ are the same as those in Proposition \ref{JN}.
\end{cor}

Except for simple changes, the proof of Corollary \ref{JNcor} is
identical to the proof of \cite[Corollary 61]{SW1}, and we refer the
reader there for its proof.\\

{\bf Proof of Proposition \ref{JN}:} The proof is an adaptation to
$[c_\rho R]BMO(E)$ of the one in \cite[Lemma 60]{SW1}.
We begin by recalling the ``dyadic grids" defined in \cite{SW3}.  Note
by \eqref{Cond0} that the quasimetric space $(\Omega,\rho)$ is
separable since $\Omega$ is separable with respect to Euclidean
distance in $\mathbb{R}^n$ .  Define $\N_{\ell}=\{1,...,\ell\}$ for each
$\ell\in\N$, and let $\N_\infty = \N$.  Set $\lambda=8\kappa^5$ with
$\kappa$ as in (\ref{tri}).  Then for each $m\in\mathbb{Z}$ and every
$k\geq m$, there are points $\{x_j^k\}_{j=1}^{n_k}\subset\Omega$
$(n_k\in\N\cup\{\infty \})$ and Borel sets $\{Q^k_j\}_{j=1}^{n_k}$
satisfying
\bea\label{dyadicstruct}
\label{ds1}& B(x_j^k,\lambda^k)\subset Q^k_j\subset
B(x_j^k,\lambda^{k+1})\mbox{ if } j\in\N_{n_k}, \\
\label{ds2}& \Omega = \cup_{j=1}^{n_k} Q_j^k,\\
\label{ds3}& Q_i^k\cap Q_j^k = \emptyset\mbox{ if } i, j\in
\N_{n_k}\text{ and } i\neq j,\\
\label{ds4}& \mbox{either } Q_j^k\subset Q_i^l\mbox{ or } Q_j^k\cap
Q_i^l = \emptyset \mbox{ if }  k <l, j\in\N_{n_k}\text{ and }i\in\N_{n_l}.
\eea
This dyadic grid of Borel sets depends on the integer $m$, and
there may be different grids for each $m$.  We fix a single grid for
each $m\in\mathbb{Z}$ and denote it by ${\F}_m$:
\bea \label{fset1}{\F}_m = \{Q_j^k\;:\; k,j\in\mathbb{Z}, \;k\geq
m\text{ and }j\in\N_{n_k}\}.
\eea
For fixed $m,k,j\in\mathbb{Z}$ with $k\geq m$ and $j\in\N_{n_k}$ we
will refer to the Borel set $Q_j^k\in\F_m$ as the $j^\text{th}$
``\emph{cube}'' at level $k$. For $0<\delta_0\leq 1$, we will call a
cube $Q^k_j\in\F_m$ ``\emph{$\delta_0$-local}''
if $B(x^k_j,\lambda^{k+1})$ satisfies $\lambda^{k+1} < \delta_0
R(x_j^k)$. For each $m\in\mathbb{Z}$ and $\delta_0\in (0,1]$ we define
\bea
\label{fsets}&\bullet& \F_{m,\delta_0} = \{Q\in\mathcal{F}_m \;:\;
Q\mbox{ is $\delta_0$-local}\}, \mbox{ and}\\
&\bullet& \E_{m,\delta_0} = \{Q=Q_j^k \in \F_{m,\delta_0}\;:\;
B(x_j^k,\lambda^{k+1}) \mbox{ is an $R(x)$-ball for some $x\in
E$}\}.\nonumber
\eea

Set $c_\rho = (\gamma^*)^2\lambda$ and fix $f\in [c_\rho R]BMO(E)$
with $||f||_{[c_\rho R]BMO(E)}=1$.  Let $f^m$  be the discrete
expectation of $f$ on the dyadic grid at level $m$:
\bea\label{expectf} f^m(z) = \sum_{j\in\N_m} \Big(\frac{1}{|Q_j^m|}\int_{Q_j^m} fdy\Big)\chi_{Q_j^m}(z).
\eea
For the moment, we will assume each of the following.
\begin{center}
\begin{itemize}
\item There are positive constants $C_8',C_9'$ and $\delta_1$ with
$\delta_1\leq A_*\lambda^{-1}$ such that for each $m\in\mathbb{Z}$,
$\alpha>0$ and $Q\in \mathcal{E}_{m,\delta_1}$, we have
\bea\label{weaktypefm} |\{y\in Q\;:\; |f^m -f_Q|>\alpha\}| \leq
C_8'e^{-{C_9'\alpha}}|Q|. \eea
Note that $C_8'$, $C_9'$ and $\delta_1$ are independent of $m$, and $C_8',C_9'$ are also independent of $E$. %
\item For almost every $y\in\Omega$,
\bea\label{pwc} f^m(y) \ra f(y) \mbox{ as } m\ra -\infty.
\eea
\end{itemize}
\end{center}

Taking (\ref{weaktypefm}) and (\ref{pwc}) temporarily for granted, let
us now prove Proposition \ref{JN} by using a packing argument and
Fatou's lemma.  To begin, we will use (\ref{weaktypefm}) to derive its
analogue where the $\delta_1$-local cube $Q$ is replaced by any
$\delta_0R(x)$-ball, for any $x\in E$, provided $\delta_0$ is
sufficiently small in terms of $\delta_1$ above.  Indeed, fix $x\in
E$, set $\delta_0 = A_*\delta_1/[(\gamma^*)^3\lambda]$, let
$B=B(z,r)$ be a $\delta_0R(x)$-ball and let $m\in \mathbb{Z}$ with
$\lambda^{m+1}<r$.  Choose $k\in\mathbb{Z}$ with $k>m$ such that
$\lambda^k<r\leq\lambda^{k+1}$.  As $\Omega = \cup_j Q_j^k$, there is
a nonempty collection $\mathcal{G}\subset \N_{n_k}$ such that
\bea Q_j^k \cap B &\neq& \emptyset\mbox{ for all }j\in \mathcal{G} \mbox{, and}\nonumber\\
\label{pack} B \subset \cup_{j\in\mathcal{G}} Q_j^k\subset
\cup_{j\in\mathcal{G}}B(x_j^k,\lambda^{k+1}) &\subset&
B(z,\gamma^*\lambda r)=B^*\subset B(x,(\gamma^*)^2\lambda\delta_0R(x)).
\eea
Here, the third and fourth containments in (\ref{pack}) follow from
(\ref{swallowing}) since $\lambda^{k+1}<\lambda r$.  We now prove
that the set $\mathcal{E}_{m,\delta_1}$ is nonempty.
\begin{lem}\label{oldclaim}  With $x,B,m,k,\delta_1, \mathcal{G}$ and
$\delta_0$ as above, $Q_j^k\in \E_{m,\delta_1}$ for every
$j\in\mathcal{G}$.
\end{lem}

\noindent {\bf Proof of Lemma \ref{oldclaim}}: It is enough to show
that for each $j\in\mathcal{G}$,
\bea \label{packA} &&\lambda^{k+1}<R(x) \mbox{ and } B(x_j^k,\gamma^*\lambda^{k+1})\subset B(x,R(x))\mbox{ (showing that $B(x_j^k,\lambda^{k+1})$ is}\\
&&\mbox{ an $R(x)$-ball), and} \nonumber\\
\label{packB}&& \lambda^{k+1}<\delta_1R(x_j^k) \mbox{ (showing that
$Q_j^k\in \mathcal{F}_{m,\delta_1}$).}
\eea
 Fix $j\in\mathcal{G}$.  To see that (\ref{packA}) holds, we begin by
noting that our choice of $\delta_0$ guarantees that
$(\gamma^*)^3\lambda\delta_0<1$.  Due to our choice of $k$ and using
that $B(z,r)$ is a $\delta_0R(x)$-ball we have
 \bea \lambda^{k+1} < \lambda r < \lambda \delta_0R(x)<R(x).\nonumber
 \eea
Also, since $x_j^k\in B^*$, swallowing gives
$B(x_j^k,\gamma^*\lambda^{k+1})\subset
B(x,(\gamma^*)^3\lambda\delta_0R(x))\subset B(x,R(x))$ establishing
(\ref{packA}). Next, since $R(x)$ satisfies the uniformity condition
(\ref{unifr}) on $E$ with constant $A_*$ and $x_j^k\in B(x,R(x))$,
$A_*R(x)<R(x_j^k)$.  Our choice of $\delta_0$ then guarantees that
 \bea \lambda^{k+1} < \lambda r < \delta_1A_*R(x) <
 \delta_1R(x_j^k),\nonumber
 \eea
 giving (\ref{packB}) and proving the lemma.
\begin{flushright}{$\Box$}
\end{flushright}

\vspace{0.1in}
Next, since $x_j^k\in B^* \cap B(x_j^k,\gamma^*\lambda r)$, we have by
the swallowing lemma, the local doubling property (\ref{doubling1})
and the dyadic structure that
\bea |B^*| \leq |B(x_j^k,(\gamma^*)^2\lambda r)|\leq C_0(\gamma^*\lambda)^{2d_0}|B(x_j^k,\lambda^k)| \leq C'|Q_j^k|\nonumber
\eea
for any $j\in\mathcal{G}$ with $C'=C_0(\gamma^*\lambda)^{2d_0}$.
Therefore, for each $j\in \mathcal{G}$,
\bea |f_B - f_{Q_j^k}| &\leq& |f_{B^*}-f_B| + |f_{B^*} - f_{Q_j^k}| \nonumber\\
&=&\Big|\fracc{1}{|B|}\int_{B} (f-f_{B^*})dy\Big|+\Big|\fracc{1}{|Q_j^k|}\int_{Q_j^k} (f-f_{B^*})dy\Big|\nonumber\\
\label{JN10}&\leq& \fracc{C_0(\gamma^*\lambda)^{d_0}+C'}{|B^*|}\int_{B^*} |f-f_{B^*}|dy\nonumber\\
&\leq& C_0(\gamma^*\lambda)^{d_0}+C' = C''
\eea
since $||f||_{[c_\rho R]BMO(E)} = 1$ and $B^*$ is an $R(x)$-ball (due
to our choice of $\delta_0$) with $x\in E$, and hence $B^*$ is also a $c_\rho
R(x)$-ball.

Consequently, if $y\in B$ and $\alpha>2C'',$ then (\ref{JN10}) and the
standard triangle inequality imply that $|f^m(y)-f_{Q_j^k}| >
\alpha/2$ provided $|f^m(y) - f_B|>\alpha$.  Since
$Q_j^k\in \E_{m,\delta_1}$ for $j\in\mathcal{G}$, the
disjointness in $j$ of the $Q_j^k$ and
(\ref{weaktypefm}) yield
\bea |\{y\in B\;:\;|f^m(y)-f_B|>\alpha\}|&\leq&\sum_{j\in\mathcal{G}} |\{y\in Q_j^k\;:\; |f^m(y)-f_{Q_j^k}|>\alpha/2\}|\nonumber\\
&\leq&\sum_{j\in\mathcal{G}} C_8'e^{-C_9'\alpha/2}|Q_j^k|\nonumber\\
&\leq& C_8'e^{-C_9'\alpha/2}|B^*|, \mbox{ as }\disp\cup_{j\in\mathcal{G}} Q^k_j\subset B^*,\nonumber\\
&\leq& C_8'  e^{-C_9'\alpha/2}C_0(\gamma^*\lambda)^{d_0} |B| \mbox{
by (\ref{doubling1})}.\nonumber \eea In case $\alpha\leq 2C'',$ we
simply use that $|\{y\in B\;:\;|f^m(y)-f_B|>\alpha\}| \leq |B|$ and
replace $C_8'$ with $e^{C_9'C''}$ if necessary.  Hence, there is a
constant $C>0$ independent of $x\in E$ such that for any $\alpha>0$
and any $\delta_0 R(x)$-ball $B$, \bea \label{weaktypeB} |\{y\in
B\;:\; |f^m(y)-f_B|>\alpha\}| \leq Ce^{-C_9'\alpha/2} |B|. \eea Next
we use the pointwise convergence of $f^m$ to $f$ as $m\ra -\infty$
and Fatou's lemma.  Set $E_{m,\alpha} =  \{y\in B\;:\;
|f^m(y)-f_B|>\alpha\}$.  Then
\bea |\{y\in B\;:\; |f(y)-f_B|>\alpha\}| &=& \int \chi_{\{y\in B\;:\; |f(y)-f_B|>\alpha\}}(z)dz\nonumber\\
&\leq& \int \liminf_{m\ra -\infty}\chi_{E_{m,\alpha}}(z)dz\nonumber\\
&\leq& \liminf_{m\ra -\infty}|E_{m,\alpha}|\nonumber\\
&\leq& C  e^{-C_9'\alpha/2} |B|. \eea This proves (\ref{JN1}) with
$C_8=C$ and $C_9=C_9'/2$ in case $||f||_{[c_\rho R]BMO(E)} = 1$ and
$\delta_0 =A_*\delta_1/ [(\gamma^*)^3\lambda]$. The general case
follows by replacing $f$ and $\alpha$ by $f/||f||_{[c_\rho
R]BMO(E)}$ and $\alpha/||f||_{[c_\rho
R]BMO(E)}$ as in \cite{SW1}.\\

The proof now rests on the validity of (\ref{weaktypefm}) and
(\ref{pwc}).  We first prove (\ref{pwc}); the
verification of (\ref{weaktypefm}) is contained in Lemma \ref{310} to
follow.  Given a fixed $x\in\Omega$, the dyadic structure provides a
sequence $\{x_m\}_{m=-1}^{-\infty}\subset\Omega$ such that
\bea (i)&&\mbox{ $x_m = x_{j_m}^m$ for some $j_m\in\N_{n_m}$,
and}\nonumber\\
(ii)&& x\in Q_m = Q_{j_m}^m \subset B(x_m,\lambda^{m+1}) \mbox{ for
each $m$.}\nonumber
\eea
By standard homogeneous space theory (see the proof of Lemma \ref{310}
for further details), almost every point $x\in \Omega$ is a
Lebesgue point of $f$:
\bea\label{diffint}
\lim_{r \to 0} \frac{1}{|B(x,r)|}\int_{B(x,r)} |f(y)-f(x)| dy = 0.
\eea
Fix such an $x$. By Lemma \ref{swalemma}, there exist $m_0$ and $C$ such
that $B(x_m,\lambda^{m+1}) \subset B(x,C\lambda^{m+1})$ if $m\le
m_0$. Thus $Q_m \subset B(x,C\lambda^{m+1})$. Also, by (\ref{ds1}), there
exists $m_1$ such that $|Q_m| \approx |B(x_m,\lambda^{m+1})| \approx
|B(x,C\lambda^{m+1})|$ uniformly in $m$ if $m\le m_1$. By choosing $r
= C\lambda^{m+1}$, we obtain
\[
\lim_{m \to -\infty} \frac{1}{|Q_m|}\int_{Q_m} |f(y) -f(x)|dy = 0.
\]
But $f^m(x) = (1/|Q_m|)\int_{Q_m} f(y) dy$, so $f^m(x) \to f(x)$ as
$m\to -\infty$. This proves (\ref{pwc}). \\

\noindent The next lemma verifies (\ref{weaktypefm}).

\begin{lem}\label{310} Let $(\Omega,\rho,dx)$ be a local homogeneous
space as in Definition \ref{homspace}, and $E$ be an open set in
$\Omega$.  Let $R:\Omega\ra(0,\infty)$ satisfy $R(x)\leq
\min\left\{R_0(x)/(\gamma^*)^2, R_1(x)/\gamma^*\right\}$ where $R_0$
is as above Remark (2.1) and $R_1$ is as in \eqref{doubling1}.
Furthermore, assume that $R(x)$ satisfies a local uniformity
condition on $E$ with constant $A_*$.  Then there are positive
constants $C_8', C_9', \delta_1, c_\rho$ with
$\delta_1\in(0,A_*\lambda^{-1}]$ and $c_\rho>1$ such that for every
$\alpha>0$, $m\in\mathbb{Z}$, and $Q\in\mathcal{E}_{m,\delta_1}$,
\bea\label{JN2} |\{y\in Q\;:\; |f^m(y) - f_{Q}|>\alpha\}| \leq
C_8'e^{-C_9'\alpha}|Q| \eea for all $f\in [c_\rho R]BMO(E)$ with
$||f||_{[c_\rho R]BMO(E)}=1$. The constants $C_8',\; C_9'$ and
$c_\rho$ depend only on $\rho$.

\end{lem}

\noindent {\bf Proof:}  The proof is broken into five steps.\\

\noindent {\bf \emph{I}}:  Recall the dyadic structure described in
(\ref{ds1})--(\ref{ds4}), and set $c_\rho=(\gamma^*)^2\lambda =
8(\gamma^*)^2 \kappa^5$.  Let $f\in [c_\rho R]BMO(E)$ with
$||f||_{[c_\rho R]BMO(E)}=1$. Fix $m\in \mathbb{Z}$ and a cube
$Q_0=Q_{j}^k\in \mathcal{F}_{m,\lambda^{-1}}$; see (\ref{fsets}).
Our first step compares the average of $f$ on $Q_0$ with its average
on the related $\rho$-ball $B(x_j^k,\lambda^{k+1})$.  Indeed,
\bea\label{BMOI0} \left|f_{B(x_j^k,\lambda^{k+1})} - f_{Q_0}\right| &=& \Big|\frac{1}{|Q_0|}\int_{Q_0} (f-f_{B(x_j^k,\lambda^{k+1})})dx\Big|\nonumber\\
&\leq& \frac{1}{|B(x_j^k,\lambda^k)|}\int_{B(x_j^k,\lambda^{k+1})} |f-f_{B(x_j^k,\lambda^{k+1})}|dx.
\eea
Thus, as $Q_0$ is $\lambda^{-1}$-local and $R(x_j^k)\leq R_1(x_j^k)$,
the local doubling condition (\ref{doubling1}) gives
\bea \frac{1}{|Q_0|}\int_{Q_0}|f-f_{Q_0}|dx &\leq& \frac{1}{|Q_0|}
\int_{Q_0} |f-f_{B(x_j^k,\lambda^{k+1})}|dx +
|f_{B(x_j^k,\lambda^{k+1})} - f_{Q_0}|\nonumber\\
&\leq& \frac{2}{|B(x_j^k,\lambda^k)|}\int_{B(x_j^k,\lambda^{k+1})}
|f-f_{B(x_j^k,\lambda^{k+1})}|dx\mbox{ , by (\ref{BMOI0}),}\nonumber\\
&\leq&
\frac{2C_0\lambda^{d_0}}{|B(x_j^k,\lambda^{k+1})|}\int_{B(x_j^k,\lambda^{k+1})}
|f-f_{B(x_j^k,\lambda^{k+1})}|dx.\nonumber
\eea
Therefore, if $B(x_j^k,\lambda^{k+1})$ is also a $c_\rho R(x)$-ball
for some $x\in E$, we may write
\bea\label{BMOI1}
\frac{1}{|Q_0|}\int_{Q_0}|f-f_{Q_0}|dx &\leq& 2C_0\lambda^{d_0} := 2c_0.
\eea

Now, further restrict $Q_0\in\mathcal{E}_{m,\lambda^{-1}}$.  Setting
$h=h(x,Q_0) = (f(x)-f_{Q_0})\chi_{Q_0}(x)$, (\ref{BMOI1}) gives
\bea\label{BMOII0}
\frac{1}{|Q|}\int_Q |h|dx \leq \frac{1}{|Q_0|} \int_{Q_0}
|f-f_{Q_0}|dx \leq 2c_0
\eea
for any dyadic cube $Q$ for which $Q_0\subset Q$.\\

\noindent{\bf \emph{II}}: The dyadic maximal function on local cubes
nearby $E$, acting on $g\in L_{loc}^1(\Omega)$, is defined by \bea
\label{BMOIII0} M^\Delta_E g(x) = \sup \frac{1}{|Q'|}\int_{Q'}
|g(y)|dy, \eea where the supremum is taken over all cubes $Q'\in
\E_{m,A_*\lambda^{-1}}$ such that $x\in Q'$. If $x\in \Omega$ is not
a member of any cube in $\E_{m,A_*\lambda^{-1}}$ then we set
$M^\Delta_Eg(x)=0$.  Next, we note the weak-type $(1,1)$ inequality
for $M^\Delta_E$: \bea\label{BMOIII1} |\{x\in \Omega : M^\Delta_E
g(x)>\alpha\}| \leq \frac{1}{\alpha}\int_\Omega |g(x)|dx. \eea This
is a consequence of the analogous inequality in \cite{SW1} for
a larger dyadic maximal operator.\\

For $\alpha>0$, let $E_\alpha = \{ x\in \Omega : M_E^\Delta h(x) >
\alpha\}$.  If $\alpha \geq 2c_0$, then $E_\alpha \subset Q_0$.
Indeed, if $M_E^\Delta h(x)>\alpha$ there is a cube $Q'\in
\E_{m,A_*\lambda^{-1}}$ containing $x$ such that
\bea \frac{1}{|Q'|} \int_{Q'} |h|dy > 2c_0.\nonumber
\eea
Thus by the definition of $h$, $Q'$ must intersect $Q_0$.  Therefore,
(\ref{BMOII0}) and the dyadic structure give that $Q'\subsetneq
Q_0$, and so $x\in Q_0$.\\

For $\alpha>0$, let $C_\alpha$ be the collection of all cubes $Q
\in\E_{m,A_*\lambda^{-1}}$ for which $|h|_{Q} >\alpha$.  By the above,
if $\alpha\geq 2c_0$ then $Q\subsetneq Q_0$ for each $Q\in
C_\alpha$. Denote the collection of maximal cubes in $C_\alpha$ by
$S_\alpha = \{Q_{\alpha,j}\}$.  Then
\bea (i)&& \mbox{ If $\alpha\geq 2c_0$, the cubes of $S_\alpha$ are
pairwise disjoint.}\nonumber\\
(ii)&&\mbox{ If $\alpha\geq 2c_0$, then  } \cup_{Q\in S_\alpha}Q =
E_\alpha \subset Q_0.\nonumber\\
(iii)&& \mbox{ If $2c_0\leq \alpha < \beta$ and $i,j$ are given then
either $Q_{\beta,j}\subset Q_{\alpha,i}$ or $Q_{\beta,j}\cap
Q_{\alpha,i}=\emptyset$.}\nonumber
\eea
To see (i), note that if two cubes $Q_{\alpha,j}, \; Q_{\alpha,i}$
intersect, the dyadic structure implies that one is contained in
the other, violating maximality.  For (ii), let $x\in E_\alpha$.
Then there is a cube $Q'$ containing $x$ for which $|h|_{Q'}>\alpha$
and so $Q'\subset Q_{\alpha,j}$ for some $j$.  Thus $E_\alpha \subset
\cup_j Q_{\alpha,j}$ and (ii) follows.  (iii)
follows from a similar argument as for (i) using the dyadic structure
and maximality of the cubes in $S_\alpha$.\\

\noindent {\bf \emph{III}:} The local doubling condition
(\ref{doubling1}) translates to a similar property for the Lebesgue
measure of dyadic cubes.  Indeed, fix a $\lambda^{-1}$-local cube
$Q=Q_j^l$ and denote its dyadic predecessor $Q_i^{l+1}$ by $Q_1$.  By
(\ref{ds1}),
\bea\label{IV0} B(x_j^l,\lambda^l)\subset Q\subset
B(x_j^l,\lambda^{l+1}) \mbox{ and}\\
\label{IV1} B(x_i^{l+1},\lambda^{l+1})\subset Q_1\subset
B(x_i^{l+1},\lambda^{l+2})
\eea
where $\lambda^{l+1}< \lambda^{-1}R(x_j^l)\leq
(\gamma^*\lambda)^{-1}R_1(x_j^l)$.  Thus, as $Q\subset Q_1$ and
$\lambda^{l+2} < R_1(x_j^l)/\gamma^*$, Lemma \ref{translate} of
the appendix implies that
\bea\label{IV33} |B(x_i^{l+1},\lambda^{l+2})| \leq
C_0(\gamma^*\lambda^2)^{d_0}|B(x_j^l,\lambda^l)|.
\eea
Therefore, (\ref{IV0}) and (\ref{IV1}) together give
\bea\label{IV2} |Q_1|\leq C_0(\gamma^*\lambda^2)^{d_0}|Q|=c_1 |Q|.
\eea
Next, restrict $Q\in\mathcal{E}_{m,A_*\lambda^{-1}}$.  Then an
inequality similar to (\ref{BMOI1}) holds for $Q_1$, the predecessor
of $Q$.  Indeed, for such $Q$ there is an $x\in E$ such that (keeping
the same labels as in (\ref{IV0}) and (\ref{IV1}))
$\lambda^{l+1}<\min\{R(x),A_*\lambda^{-1}R(x_j^l)\}$ and
$B(x_j^l,\gamma^*\lambda^{l+1})\subset B(x,R(x))$.  Since $x_j^l\in
B(x_i^{l+1},\lambda^{l+2})$, we have that
\bea B(x_i^{l+1},\gamma^*\lambda^{l+2})\subset B(x,(\gamma^*)^2\lambda
R(x)).\nonumber
\eea
Also, since $R(x)$ satisfies a local uniformity condition on $E$ with
constant $A_*$, it follows that $\lambda^{l+1} < A_*\lambda^{-1}R(x_j^l)
< \lambda^{-1}R(x)$, giving $\lambda^{l+2}<R(x)$.
This together with the containment above shows that
$B(x_i^{l+1},\lambda^{l+2})$ is a $c_\rho R(x)$-ball.  Therefore,
\bea\label{IV34} \frac{1}{|Q_1|}\int_{Q_1} |f-f_{Q_1}|dx \leq 2c_1
\eea
using a familiar argument.\\

\noindent{\bf \emph{IV}:}  For each $\alpha\geq 2c_0$, define
$\gamma=\gamma(\alpha) = 1+ (4c_1^2/\alpha)$.  We claim that for all $j$,
\bea\label{V0} |E_{\gamma\alpha}\cap Q_{\alpha,j}|\leq
\frac{1}{2}|Q_{\alpha,j}|.
\eea
Indeed, let $Q$ be the dyadic predecessor of $Q_{\alpha,j}$.  Then,
$Q\subset Q_0$ since $Q_{\alpha,j}$ is a proper subcube of $Q_0$. The
maximality of $Q_{\alpha,j}$ then gives
\bea\label{V1} |h|_Q \leq \alpha.
\eea
Set $g = (h-h_Q)\chi_Q$ and fix $x\in E_{\gamma\alpha}\cap
Q_{\alpha,j}$.  Then, since $\gamma\geq 1$, (ii) and (iii) (see step
{\bf \emph{II}}) give that $x\in Q_{\gamma\alpha,i}$ for some $i$ and
$Q_{\gamma\alpha,i}\subset Q_{\alpha,j}$.  Using this, we have
\bea\label{V2} \gamma\alpha< |h|_{Q_{\gamma\alpha,i}} &\leq&
|g|_{Q_{\gamma\alpha,i}} + \alpha.
\eea
Consequently, for every $x\in E_{\gamma\alpha}\cap Q_{\alpha,j}$,
$M^\Delta_E g(x) > (\gamma-1)\alpha$.  Further, for each $x\in
Q\subset Q_0$ we have
\bea\label{V3} g(x) = h(x)-h_Q &=& f(x) - f_{Q_0}  - (f-f_{Q_0})_Q\nonumber\\
&=& f(x) - f_Q.
\eea
Therefore, by (\ref{BMOIII1}),
\bea\label{V4}
|E_{\gamma\alpha}\cap Q_{\alpha,j}| &\leq& |\{x:M^\Delta_E g(x)>(\gamma-1)\alpha\}|\nonumber\\
&\leq& \frac{1}{(\gamma-1)\alpha}\int_\Omega |g|dx\nonumber\\
&=&\frac{1}{(\gamma-1)\alpha}\int_Q |f-f_Q|dx\nonumber\\
&\leq&   \frac{2c_1}{(\gamma-1)\alpha}|Q|,
\eea
where the last inequality is due to (\ref{IV34}).  Inequality
(\ref{IV2}) combined with our choice of $\gamma$ gives
\bea\label{V5} |E_{\gamma\alpha}\cap Q_{\alpha,j}| \leq
\frac{2c_1^2}{(\gamma-1)\alpha}|Q_{\alpha,j}| \leq
\frac{1}{2}|Q_{\alpha,j}|,
\eea
proving (\ref{V0}).\\

\noindent{\bf \emph{V}:} For $\alpha>0$, define the distribution
function $\omega(\alpha) = |E_\alpha \cap Q_0|$.  We add (\ref{V5})
over $j$ to obtain a useful inequality for $\alpha\geq 2c_0$ (note
that $\omega(\alpha) = |E_\alpha|$ for $\alpha\geq 2c_0$, and that
$\gamma\alpha = \alpha + 4c_1^2$): \bea\label{VI0}
\omega(\alpha+4c_1^2) = \sum_j |E_{\gamma\alpha}\cap
Q_{\alpha,j}|\leq  \frac{1}{2}\sum_j
|Q_{\alpha,j}|=\frac{1}{2}|E_\alpha|= \frac{1}{2}\omega(\alpha).
\eea We now iterate (\ref{VI0}).  Fix $\alpha \geq 2c_0$.  Then
there is a $k\in \N$ such that $\alpha \in [2c_0 + 4(k-1)c_1^2, 2c_0
+ 4kc_1^2]$. Therefore, there is a $\beta\in [2c_0,2c_0+4c_1^2]$ for
which \bea\label{VI1} \omega(\alpha) \leq
\frac{1}{2^{k-1}}\omega(\beta) \leq 2e^{-k\log 2} |Q_0| \eea as
$\omega(s) \leq  |Q_0|$ for $s>0$.  Since \bea k\geq \frac{\alpha -
2c_0}{4c_1^2},\nonumber \eea we obtain \bea\label{VI2}
\omega(\alpha) \leq C_8'e^{-C_9'\alpha}|Q_0| \eea where $C_8'$ and
$C_9'$ depend on $c_0, c_1$. Finally, if $\alpha\in (0,2c_0)$ we use
that $\omega(\alpha)\leq |Q_0|$ to obtain a similar estimate. Hence,
for all $\alpha >0$, \bea\label{VI3} |\{x\in Q_0 :
M^\Delta_E[(f-f_{Q_0})\chi_{Q_0}](x)>\alpha\}| \leq
C_8'e^{-C_9'\alpha}|Q_0|. \eea

The proof will then be complete if we show that
\bea\label{VI4}
|f^m - f_{Q_0}|\chi_{Q_0}  \leq M^\Delta_E \Big[ (f-f_{Q_0})\chi_{Q_0}\Big].
\eea
Using the dyadic structure it is easy to see that (\ref{VI4}) holds
provided $Q_i^m\in\mathcal{E}_{m,A_*\lambda^{-1}}$ whenever
$Q_i^m\subset Q_0$. This proviso is true by further restricting the
size of $\delta_1$.  Set $\delta_1=\min\{A_*^3,A_*\lambda^{-1}\}$ and
suppose $Q_0=Q_j^k\in\mathcal{E}_{m,\delta_1}$.  Omitting as we may
the case when $k=m$, suppose that $Q_i^m\subset Q_0$ and $m<k$.  Recalling that
$\gamma^*=\kappa+2\kappa^2 < 8\kappa^5=\lambda$, we have
\bea\label{VInew1} B(x_i^m,\gamma^*\lambda^{m+1})&\subset&
B(x_i^m,\fracc{\gamma^*\lambda^{k+1}}{\lambda})\nonumber\\
&\subset&B(x_i^m,\lambda^{k+1})\nonumber\\
&\subset&B(x_j^k,\gamma^*\lambda^{k+1}) \subset B(x,R(x))\nonumber
\eea
for some $x\in E$.  Thus $B(x_i^m,\lambda^{m+1})$ is an $R(x)$-ball.
Since $Q_j^k\in\mathcal{F}_{m,\delta_1}$ and $x_i^m,x_j^k\in
B(x,R(x))$, the uniformity condition gives
\bea
\lambda^{m+1}<\fracc{\delta_1}{\lambda}R(x_j^k)\leq
\fracc{\delta_1}{A_*\lambda}R(x)\leq\fracc{\delta_1}{A_*^2\lambda}
R(x_i^m)\leq\fracc{A_*}{\lambda}R(x_i^m),
\eea
and therefore $Q_i^m\in \mathcal{E}_{m,A_*\lambda^{-1}}$.  This
concludes the proof of both Lemma \ref{310} and Proposition \ref{JN}
with $\delta_1=\min\{A_*^3,A_*\lambda^{-1}\}$.
\begin{flushright}
$\Box\;\Box$
\end{flushright}


\section{The Proof of Proposition \ref{harnackmain}}

Proposition \ref{harnackmain}  will be proved using the results of three
lemmas and Corollary \ref{JNcor}. The lemmas give mean-value
estimates for positive and negative powers of weak solutions as well
as a logarithmic estimate. In order to simplify their statements, we
list now some assumptions to remain in force for the rest of the
section. We always assume that $(\Omega,\rho,dx)$ is a local
homogeneous space as in Definition \ref{homspace}, that the Sobolev
inequality \eqref{sobolev} is valid, that \eqref{cutoff} holds for
some $\tau\in(0,1)$ and $s^*>p\sigma^\prime$ with $\sigma$ as in
\eqref{sobolev}, and that (\ref{E3}) and (\ref{3.6-0}) are valid for
some $t\geq1$. Our first lemma concerns positive powers of weak
solutions.

\begin{lem}\label{alpha>0} Noting the assumptions in the paragraph
above, let $(u,\nabla u)$ be a weak solution in $\Omega$ of (1.1),
where \eqref{struct} holds with exponents $\gamma,\delta,\psi$
satisfying \eqref{ranges1}. Let $s= (s^*/p)' \in [1,\sigma')$. Fix
$x_0\in\Omega$, $k>0$, $\epsilon_1,\epsilon_2,\epsilon_3\in(0,1]$, a
$\rho$-ball $B(x_0,r)$ with $0<r<\tau^2 r_1(x_0)$, set
$\bar{u}=|u|+k$, and assume that
${\bar{Z}}(B(x_0,r/\tau),\bar{u})<\infty$. If $u\geq 0$ in
$B(x_0,r)$ then for each $\alpha>0$ there exists $\alpha_1\in
[\alpha \sigma^{-\frac{1}{2}}, \alpha]$ such that \bea \label{qpos}
\disp\tesup_{B(x_0,\tau r)} \bar{u} &\leq&
C_{10}\Big(C_{11}\bar{Z}(B(x_0,r),\bar{u})\Big)^\frac{p\Psi_0}{\alpha_1}||
\bar{u}^{\alpha_1} ||_{s,B(x_0,r);\overline{dx}}^\frac{1}{\alpha_1}.
\eea Here $\Psi_0=\frac{\sigma}{\sigma-s}$, $C_{10}$ depends only on
$p,\sigma,s$ and on $\epsilon_1,\epsilon_2,\epsilon_3$ appearing in
the definition \eqref{Zetabar} of $\bar{Z}$, while $C_{11}$ depends
on $p,\sigma,s,\epsilon_1,\epsilon_2,\epsilon_3,a$ and the constants
$C_1$ from \eqref{sobolev} and $N,C_{s^*}$ in \eqref{cutoff}.
\end{lem}
\begin{rem}\label{remq>0}
\begin{itemize}
\item[i)] $\alpha_1$ is defined by
\begin{equation*}
\alpha_1 = \left\{
\begin{array}{ll}
\alpha & \text{if }
\log_{\frac{\sigma}{s}}\frac{p-1}{\alpha}\leq-\frac{1}{4}\\
\alpha & \text{if } \log_{\frac{\sigma}{s}}\frac{p-1}{\alpha}\in \big[K+\frac{1}{4},K+\frac{3}{4}\big]\,\,\text{ for some }K\in\N\cup\{0\},\\
\alpha\Big(\frac{\sigma}{s}\Big)^{-\frac{1}{2}} & \text{if }
\log_{\frac{\sigma}{s}}\frac{p-1}{\alpha}\in
\big(K-\frac{1}{4},K+\frac{1}{4}\big)\,\,\text{ for some }K\in\N\cup\{0\}.
\end{array} \right.
\end{equation*}
\item[ii)] We explicitly note that the constants $C_{10},C_{11}$ in
\eqref{qpos} are independent of $(u,\nabla u)$,
$k$, $B(x_0,r)$, $b$, $c$, $d$, $e$, $f$, $g$, $h$, and $\alpha$.
\end{itemize}
\end{rem}

\noindent {\bf Proof:} By \cite[Theorem 1.2]{MRW}, the weak solution
$(u,\nabla u)$ satisfies \bea\label{sishere}
||\bar{u}||_{L^\infty(B(x_0, r))} \leq
C \bar{Z}\big(B(x_0, \frac{r}{\tau}),\bar{u}\big)^{\Psi_0}
||\bar{u}||_{sp,B(x_0,\frac{r}{\tau});\overline{dx}}
\eea where $\bar{u} = | u |+k$ and $C>0$ depends only on $p,a$ and
$\psi$. Therefore, as $\bar{Z}\big(B(x_0,r/\tau),\bar{u}\big)<\infty$ by hypothesis,
\cite[Proposition 2.3]{MRW} gives that $\bar{u}$ is bounded on
$B(x_0,r)$.  The proof of the lemma will be completed by following the
proof of \cite[Theorem 1.2]{MRW}, but now using a modified test
function that exploits boundedness of $\bar{u}$. As in \cite{MRW}, we
may assume that $\overline{u}$ satisfies the following modified
structure conditions in terms of the functions $\bar{b},\; \bar{d},$
and $\bar{h}$ (see \cite[(3.1)]{MRW}):
\begin{eqnarray}
\nonumber\xi\cdot A(x,z,\xi)&\geq&
a^{-1}|\sqrt{Q(x)}\;\xi|^p-\bar{h}(x)\bar{z}^p,\\
\label{struct2}\left|\wt{A}(x,z,\xi)\right|&\leq&a|\sqrt{Q(x)}
\;\xi|^{p-1}+\bar{b}(x)\bar{z}^{p-1},\\
\nonumber\left|B(x,z,\xi)\right|&\leq&c|\sqrt{Q(x)}
\;\xi|^{\psi-1}+\bar{d}(x)\bar{z}^{p-1}
\end{eqnarray}
for all $(x,z,\xi)\in \Omega\times\mathbb{R}\times\mathbb{R}^n$ where
$A,\tilde{A}$ and $B$ are as in (\ref{struct}) and $\bar{z} = z+k$.
For simplicity, we will often not indicate the dependence of
$A,\tilde{A},B$, etc. on their variables.

Choose a nonnegative $\eta\in Lip_0(B(x_0,r))$ and
set $v = \eta^p\bar{u}^q$ for $q\in (1-p,0)\cup(0,\infty)$.   By
Corollary \ref{cor1}, $v$ is a feasible nonnegative test function
for any value of $q$ in the indicated range.  Corollary \ref{cor1}
implies that \bea \nabla v\cdot A(x,u,\nabla u) + v B(x,u,\nabla
u)&=&\sqrt{Q}\nabla v\cdot\tilde{A} + vB\nonumber\\
\label{s1}&=&\Big(p\eta^{p-1}\bar{u}^q\sqrt{Q}\nabla\eta
+ q\eta^p\bar{u}^{q-1}\sqrt{Q}\nabla\bar{u}\Big)\cdot \tilde{A} +
\eta^p\bar{u}^qB.
\eea

We now use (\ref{s1}) to derive some pointwise estimates. If $q>0$, we
apply (\ref{struct2}) to (\ref{s1}), giving
\bea \label{q>0.01}\nabla v\cdot A(x,u,\nabla u)+vB(x,u,\nabla u)
&\geq& q\eta^p\bar{u}^{q-1}\Big[a^{-1}|\sqrt{Q}\nabla\bar{u}|^p-
\bar{h}\bar{u}^p\Big]\nonumber\\
&& -p\eta^{p-1}\bar{u}^q|\sqrt{Q}\nabla \eta||\tilde{A}| -
\eta^p\bar{u}^q|B| \nonumber\\
&\geq& a^{-1}q\eta^p\bar{u}^{q-1}|\sqrt{Q}\nabla\bar{u}|^p-
q\bar{h}\eta^p\bar{u}^{q+p-1}\\
&& -ap\eta^{p-1}\bar{u}^q|\sqrt{Q}\nabla\eta||\sqrt{Q}\nabla
\bar{u}|^{p-1}-p\eta^{p-1}\bar{b}\bar{u}^{q+p-1}|\sqrt{Q}\nabla
\eta|\nonumber\\
&&  - c\eta^p|\sqrt{Q}\nabla\bar{u}|^{\psi-1}\bar{u}^q-
\bar{d}\eta^p\bar{u}^{q+p-1}.\nonumber
\eea
If $q<0$ we arrange (\ref{s1}) differently.  For the second term
inside the parentheses on the right side of (\ref{s1}), since $q<0$,
the first estimate of (\ref{struct2}) gives
\bea \label{moredetails1}
q\eta^p\bar{u}^{q-1}\sqrt{Q}\nabla\bar{u}\cdot \tilde{A} &=& q\eta^p\bar{u}^{q-1}\nabla\bar{u}\cdot A\nonumber\\
&\leq& -a^{-1}|q|\eta^p \bar{u}^{q-1}|\sqrt{Q}\nabla{\bar{u}}|^p + |q|\eta^p\bar{u}^{q+p-1}\bar{h}.
\eea
After estimating the other terms of (\ref{s1}) as before, we move the
first term on the right of (\ref{moredetails1}) to the left and obtain
\begin{eqnarray}
\label{q<0.1}&& \nabla v\cdot A(x,u,\nabla u) + v B(x,u,\nabla u) +
a^{-1}|q|\eta^p\bar{u}^{q-1}|\sqrt{Q}\nabla \bar{u}|^p\leq\\
\nonumber&&\qquad\qquad
 |q|\eta^p\bar{u}^{p+q-1}\bar{h}  +
ap\eta^{p-1}\bar{u}^q|\sqrt{Q}\nabla \eta||\sqrt{Q}\nabla
\bar{u}|^{p-1} +p\bar{b}\eta^{p-1}\bar{u}^{q+p-1}|\sqrt{Q}\nabla
\eta|\\
\nonumber&&\qquad\qquad
+c\eta^p\bar{u}^q|\sqrt{Q}\nabla\bar{u}|^{\psi-1}
+\bar{d}\eta^p\bar{u}^{q+p-1}.
\end{eqnarray}

 \noindent Since $u$ is a weak
solution of (\ref{eqdiff}) and $v$ is a feasible test function, we
have that \bea \int_\Omega\nabla v\cdot A+v B=\int_{B(x_0,r)}\nabla
v\cdot A+v B = 0.\nonumber \eea Integrating either (\ref{q>0.01}) or
(\ref{q<0.1}) over $B=B(x_0,r)$, we obtain that for any $q\in
(1-p,0)\cup(0,\infty)$,

\bea \label{q>0.1}|q|\fint_{B} \eta^p\bar{u}^{q-1}|\sqrt{Q}\nabla\bar{u}|^pdx &\leq& C\Big\{ |q|\fint_{B}\bar{h}\eta^p\bar{u}^{p+q-1}dx + \fint_{B} \eta^{p-1}\bar{u}^q|\sqrt{Q}\nabla \eta||\sqrt{Q}\nabla\bar{u}|^{p-1}dx\nonumber\\
&& +\fint_{B}\bar{b}\eta^{p-1}|\sqrt{Q}\nabla\eta|\bar{u}^{p+q-1}dx + \fint_{B} c\eta^p\bar{u}^q|\sqrt{Q}\nabla\bar{u}|^{\psi-1}dx\\
&&+\fint_{B}\bar{d}\eta^p\bar{u}^{p+q-1}dx\Big\},\nonumber \eea
where the constant $C$ in \eqref{q>0.1} depends only on $a,p$. Now
use Young's inequality (\ref{young}) with $\beta=p'$ and $\theta =
p'|q|/(4C)$, where $C$ is as in (\ref{q>0.1}), on the second term of
the right side of (\ref{q>0.1}).  This gives \bea
\label{q>0.2}\fint_{B} \eta^{p-1}\bar{u}^q|\sqrt{Q}\nabla\eta
||\sqrt{Q}\nabla\bar{u}|^{p-1}dx &\leq&
\frac{|q|}{4C}\fint_{B}|\sqrt{Q}\nabla\bar{u}|^p\eta^p\bar{u}^{q-1}dx\\
&& + c_2|q|^{1-p}\fint_{B}|\sqrt{Q}\nabla\eta|^p\bar{u}^{q+p-1}
dx.\nonumber \eea Here $c_2$ depends only on $p,a$. Applying Young's
inequality (\ref{young}) to the fourth term on the right side of
(\ref{q>0.1}) with
$\beta=\frac{p}{\psi-1},\;\beta'=\frac{p}{p+1-\psi}$ and
$\theta=\frac{|q|p}{4(\psi-1)C}$ yields \bea
\label{q>0.3}\fint_{B}c\eta^p\bar{u}^q|\sqrt{Q}\nabla\bar{u}|^{\psi-1}dx
&=&\fint_{B}|\sqrt{Q}\nabla\bar{u}|^{\psi-1}\eta^{\psi-1}
\bar{u}^\frac{(\psi-1)(q-1)}{p}\cdot c\eta^{p+1-\psi}\bar{u}^{q-
\frac{(\psi-1)(q-1)}{p}}dx\nonumber\\  &\leq& \frac{|q|}{4C}
\fint_{B}|\sqrt{Q}\nabla\bar{u}|^p\eta^p\bar{u}^{q-1}dx
+c_3|q|^\frac{1-\psi}{p+1-\psi}\fint_{B}c^\frac{p}{p+1-\psi}\eta^p
\bar{u}^{q+\frac{\psi-1}{p+1-\psi}}dx. \eea Since under our
hypotheses $\psi\in [p,p+1-\sigma^{-1})$, the constant $c_3$ can be
chosen as to depend only on $p,a,\sigma$. Inserting (\ref{q>0.2})
and (\ref{q>0.3}) into (\ref{q>0.1}) and absorbing two terms, we
obtain \bea\label{q>0.4}
|q|\fint_{B}\eta^p\bar{u}^{q-1}|\sqrt{Q}\nabla\bar{u}|^pdx
&\leq&C\Big\{
|q|^{1-p}\fint_{B} |\sqrt{Q}\nabla\eta|^p\bar{u}^{p+q-1}dx\nonumber\\
&&+ \fint_{B}\bar{b}\eta^{p-1}|\sqrt{Q}\nabla\eta|\bar{u}^{p+q-1}dx\nonumber\\
&&+ |q|^\frac{1-\psi}{p+1-\psi}\fint_{B} c^\frac{p}{p+1-\psi}\eta^p\bar{u}^{q+\frac{\psi-1}{p+1-\psi}}dx  \\
&& +
|q|\fint_{B}\bar{h}\eta^p\bar{u}^{p+q-1}dx\nonumber+\fint_{B}\bar{d}\eta^p\bar{u}^{p+q-1}dx\Big\},\nonumber
\eea with $C$ depending only on $a,p,\sigma$. This inequality is
identical to \cite[(3.8)]{MRW} with $\mu=0$. Therefore, we follow
the proof of \cite[Theorem 1.2]{MRW} through steps 5 and 6 with
$Y=p+q-1$, $t=\frac{p}{p+1-\psi}$ and $T=\frac{1-\psi}{p+1-\psi}$.
Note that when dealing with term III in {\emph step 5} of
\cite{MRW}, the exponent $T+p$ may be negative for some values of
$\psi\in [p,p+1-\sigma^{-1})$.  Thus we replace $|q|^{T+p}$ with the
larger term $(|q|+ |q|^{-1})^{|T|+p}$ to arrive at an analogous
inequality to \cite[(3.22)]{MRW}, recalling the notation given in
(\ref{n1}): \bea\label{q>0.5} ||  \eta
\bar{u}^\frac{Y}{p}||_{p\sigma,B;\overline{dx}} &\leq&
C\Big(|q|+\frac{1}{|q|}\Big)^{\tilde{b}_*}\bar{Z}\,\Big\{||\eta
\bar{u}^{\frac{Y}{p}}||_{p,B;\overline{dx}}+r||\bar{u}^{\frac{Y}{p}}
\sqrt{Q}\nabla\eta||_{p,B;\overline{dx}}\Big\} \eea where $B =
B(x_0,r)$, $\bar{Z} = \bar{Z}\big(B(x_0,r),\bar{u}\big)$ and
$\tilde{b}_*\geq b_*>0$ with $b_*$ as in \cite[(3.22)]{MRW}. We
explicitly note that $C$ now depends on $p,a,\sigma$ and on the
constant $C_1$ appearing in \eqref{sobolev}, while $\tilde{b}_*$
depends on $p,\sigma$ and on $\epsilon_1,\epsilon_2,\epsilon_3$ that
appear in the definition \eqref{Zetabar} of $\bar{Z}$.

We now choose $\eta =\eta_j, j \ge 1,$ as in (\ref{cutoff}).  Let
$S_j = supp\;\eta_j$ for $j\ge 1$ and $S_0 = B= B(x_0,r)$. Recall
that $\eta_j=1$ on $S_{j+1}$ and $B(x_0,\tau r)\subset
\disp\cap_{j}S_j$. Since $s^*>p\sigma'$ and $s'p=s^*$, H\"older's
inequality and (\ref{cutoff}) give \bea\label{preiteration}
||\bar{u}^\frac{Y}{p}\chi_{S_{j+1}}||_{p\sigma, B;\overline{dx}}
&\leq&
C\Big(|q|+\frac{1}{|q|}\Big)^{\tilde{b}_*}\bar{Z}N^j||\bar{u}^\frac{Y}{p}
\chi_{S_j}||_{sp, B;\overline{dx}}, \eea which is analogous to
inequality \cite[(3.23)]{MRW} and where $C$ depends on $p,a,\sigma$,
on the constants $C_1$ appearing in \eqref{sobolev} and $N,C_{s^*}$
in \eqref{cutoff}. Now for $\omega\neq0$ and $j\in\N\cup\{0\}$
define \bea \label{PHI} \Phi(j;\omega) = \Big(\fint_{B(x_0,r)}
\bar{u}^\omega\chi_{S_j}dx\Big)^{1/\omega}. \eea By
(\ref{preiteration}), noting that $Y=p+q-1>0$ for all $q\in
(1-p,0)\cup(0,\infty)$, we have \bea\label{PHIit.1}
\Phi(j+1;Y\sigma) \leq
C^\frac{p}{Y}\Big(|q|+\frac{1}{|q|}\Big)^\frac{p\tilde{b}_*}{Y}\bar{Z}^\frac{p}{Y}N^\frac{jp}{Y}\Phi(j;sY).
\eea

Inequality (\ref{PHIit.1}) will be iterated to finish the proof.
Indeed, set ${\X}=\sigma/s >1$ and fix $\alpha_1>0$ as in Remark
\ref{remq>0}.  Set $q_j = \alpha_1{\X}^j+1-p$ and $Y_j =
\alpha_1{\X}^j$ for each $j\in\N\cup\{0\}$.

\emph{Claim:} We claim that $q_j \in (1-p,0) \cup (0,\infty)$ for
$j\in \N\cup\{0\}$ and that \bea\label{uglyq} |q_j|+\frac{1}{|q_j|}
&\leq& \X^j\left[\alpha_1\X^{\frac{1}{2}} + p -1  +
\frac{1}{(p-1)(1-\X^{-\frac{1}{4}})}\right].\eea We start noting
that from $(\X^\frac{1}{8}-\X^{-\frac{1}{8}})^2\geq0$ it follows
that
\begin{equation}\label{19}
\X^\frac{1}{4}-1\geq1-\X^{-\frac{1}{4}}.
\end{equation}
If $\log_{\X}\frac{p-1}{\alpha}\leq-\frac{1}{4}$, then
$\alpha_1=\alpha\geq(p-1)\X^\frac{1}{4}$. Thus for every
$j\in\N\cup\{0\}$ we have
$$q_j=\alpha\X^j+1-p\geq\alpha+1-p>\alpha\X^{-\frac{1}{4}}+1-p\geq0$$
and hence, also using \eqref{19},
\begin{eqnarray}
\nonumber|q_j|+\frac{1}{|q_j|} &\leq&\alpha\X^j +p-1+\frac{1}{\alpha+1-p}\\
\nonumber&\leq&\alpha \X^j + p -1  +
\frac{1}{(p-1)(\X^{\frac{1}{4}}-1)}\,\,\leq\,\,\alpha\X^j
+p-1+\frac{1}{(p-1)(1-X^{-\frac{1}{4}})}.
\end{eqnarray}
Since $\X>1$ and $\alpha=\alpha_1$, \eqref{uglyq} easily follows.

If $\log_{\X}\frac{p-1}{\alpha}>-\frac{1}{4}$, there exists a unique
$K\in\N\cup\{0\}$ such that either
$\log_{\X}\frac{p-1}{\alpha}\in(K-\frac{1}{4},K+\frac{1}{4})$ or
$\log_{\X}\frac{p-1}{\alpha}\in[K+\frac{1}{4},K+\frac{3}{4}]$. In
the first case $\alpha_1=\alpha\X^{-\frac{1}{2}}$ and
$$\alpha\X^{K-\frac{1}{4}}<p-1<\alpha\X^{K+\frac{1}{4}}.$$ Thus if
$j\leq K$ one has
$$q_j=\alpha\X^{j-\frac{1}{2}}+1-p\leq\alpha\X^{K-\frac{1}{2}}+1-p=
\alpha\X^{K-\frac{1}{4}}\X^{-\frac{1}{4}}-(p-1)<-(p-1)(1-\X^{-\frac{1}{4}})<0.$$
On the other hand, if $j\geq K+1$, we have
$$q_j=\alpha\X^{j-\frac{1}{2}}+1-p\geq\alpha\X^{K+\frac{1}{2}}+1-p=
\alpha\X^{K+\frac{1}{4}}\X^{\frac{1}{4}}-(p-1)>(p-1)(\X^{\frac{1}{4}}-1)>0.$$
Thus $q_j\in(1-p,0)\cup(0,\infty)$, and moreover for $j\leq K$
\begin{equation}\label{20}
|q_j|+\frac{1}{|q_j|}\leq
\alpha\X^j+p-1+\frac{1}{(p-1)(1-\X^{-\frac{1}{4}})},
\end{equation}
while for $j\geq K+1$ we have
\begin{equation}\label{21}
|q_j|+\frac{1}{|q_j|}\leq
\alpha\X^j+p-1+\frac{1}{(p-1)(\X^{\frac{1}{4}}-1)}\leq\alpha\X^j+p-1+\frac{1}{(p-1)(1-\X^{-\frac{1}{4}})}.
\end{equation}
From the previous inequalities, since $\X>1$ and
$\alpha=\alpha_1\X^\frac{1}{2}$, \eqref{uglyq} follows.

It remains to consider the case when
$\log_{\X}\frac{p-1}{\alpha}\in[K+\frac{1}{4},K+\frac{3}{4}]$ for
some $K\in\N\cup\{0\}$. Then we have $\alpha_1=\alpha$ and
$$\alpha\X^{K+\frac{1}{4}}\leq p-1\leq\alpha\X^{K+\frac{3}{4}}.$$
Now if $j\leq K$ we have
$$q_j=\alpha\X^j+1-p\leq\alpha\X^K+1-p=\alpha\X^{K+\frac{1}{4}}\X^{-\frac{1}{4}}-(p-1)\leq-(p-1)(1-\X^{-\frac{1}{4}})<0,$$
while if $j\geq K+1$
$$q_j=\alpha\X^j+1-p\geq\alpha\X^{K+1}+1-p=\alpha\X^{K+\frac{3}{4}}\X^\frac{1}{4}-(p-1)\geq(p-1)(\X^{\frac{1}{4}}-1)>0.$$
Hence also in this last case $q_j\in(1-p,0)\cup(0,\infty)$ for every
$j\in\N\cup\{0\}$, and for $j\leq K$ we have \eqref{20}, while for
$j\geq K+1$ we have \eqref{21}. The proof of the \emph{claim} is
complete.

Let $c_4=p -1  + ((p-1)(1-\X^{-\frac{1}{4}}))^{-1}$. By \eqref{PHIit.1} and \eqref{uglyq}, for each
$j\in\N\cup\{0\}$,
\bea \label{PHIit.21} \Phi(j+1;\alpha_1
s\X^{j+1}) &\leq& \Big[C^{\X^{-j}}(\alpha_1\X^\frac{1}{2}+c_4)^
{\tilde{b}_*\X^{-j}}\bar{Z}^{\X^{-j}} \X^{\tilde{b}_*j\X^{-j}}
N^{j\X^{-j}}\Big]^\frac{p}{\alpha_1}\Phi(j;\alpha_1 s\X^j).\nonumber
\eea
Iterating this inequality we see that
\bea \label{PHIit.2}
\Phi(j+1;\alpha_1s\X^{j+1})&\leq& \Big[C^{\Psi_0}(\alpha_1\X^\frac{1}{2}+c_4)^{\tilde{b}_*\Psi_0}
\X^{\tilde{b}_*\Psi_1} N^{\Psi_1}
\bar{Z}^{\Psi_0}\Big]^\frac{p}{\alpha_1}\Phi(0;\alpha_1 s) \eea
for each $j\in \N \cup\{0\}$ where we have set $\Psi_0 = \disp\sum_{j=0}^\infty \X^{-j}$ and $\Psi_1 =
\disp\sum_{j=0}^\infty j\X^{-j}$. Now, since the function
$z^\frac{1}{z}$ achieves its maximum for $z\in(0,\infty)$ at $z=e$,
we have
\begin{eqnarray}
\nonumber(\alpha_1\X^\frac{1}{2}+c_4)^\frac{\tilde{b}_*\Psi_0p}{\alpha_1}&\leq&\big(2\max\{\alpha_1\X^\frac{1}{2},c_4\}\big)^\frac{\tilde{b}_*\Psi_0p}{\alpha_1}
\,\,=\,\,
\left[\big(\max\{2\alpha_1\X^\frac{1}{2},2c_4\}\big)^\frac{1}{2\alpha_1\X^\frac{1}{2}}\right]^{2\tilde{b}_*\X^\frac{1}{2}\Psi_0p}\\
\label{22}&\leq&\left[\max\Big\{e^\frac{1}{e},(2c_4)^\frac{1}{2\alpha_1\X^\frac{1}{2}}\Big\}\right]^{2\tilde{b}_*\X^\frac{1}{2}\Psi_0p}\,\,\leq\,\,c_5c_6^\frac{\Psi_0p}{\alpha_1},
\end{eqnarray}
with $c_5,c_6$ depending on
$p,\sigma,s,\epsilon_1,\epsilon_2,\epsilon_3$. Next, we set
$\Phi(\infty;\infty)
=\disp\limsup_{j\ra\infty}\Phi(j;\alpha_1s\chi^j)$.  Since the right
side of (\ref{PHIit.2}) is independent of $j$, we may allow $j\ra\infty$ in (\ref{PHIit.2}) and use \eqref{22} to obtain
\bea\label{PHIit.3} \Phi(\infty;\infty) &\leq&
C_{10}\Big[C_{11}\bar{Z}\Big]^\frac{\Psi_0p}{\alpha_1}\Phi(0;\alpha_1
s), \eea with $C_{10}$ depending on
$p,\sigma,s,\epsilon_1,\epsilon_2,\epsilon_3$, $C_{11}$ depending on
$p,\sigma,s,\epsilon_1,\epsilon_2,\epsilon_3,a$ and both also depending on the constants
$C_1$ from \eqref{sobolev} and $N,C_{s^*}$ in \eqref{cutoff}. Since
$B(x_0,\tau r)\subset S_j$ for all $j\ge1$ we have
$\disp\esup_{B(x_0,\tau r)} \bar{u} \leq \Phi(\infty;\infty)$, see
for instance \cite{GT}, and therefore we conclude that \bea
\disp\esup_{B(x_0,\tau r)} \bar{u} &\leq& C_{10}\Big[C_{11}\bar{Z}
\Big]^\frac{\Psi_0p}{\alpha_1}||\bar{u}^{\alpha_1}||_{s,
B(x_0,r);\overline{dx}}^\frac{1}{\alpha_1}\, , \eea which completes
the proof of (\ref{qpos}).
\begin{flushright}$\Box$
\end{flushright}

\begin{lem}\label{alpha<0}  Let the assumptions in the opening
paragraph of this section hold, let $s$ be as in Lemma \ref{alpha>0}
and suppose that $(u,\nabla u) \in W^{1,p}_Q(\Omega)$ is a weak
solution in $\Omega$ of \bea\label{supersol1}
\mbox{\emph{div}}(A(x,u,\nabla u)) \leq B(x,u,\nabla u) \eea where
$A,B$ satisfy \eqref{struct} with exponents $\gamma,\delta,\psi$
satisfying \eqref{ranges1}. Fix $x_0\in\Omega$, $k>0$,
$\epsilon_1,\epsilon_2,\epsilon_3\in(0,1]$ and a quasimetric
$\rho$-ball $B(x_0,r)$ with $0<r<\tau^2 r_1(x_0)$, set
$\bar{u}=|u|+k$, and assume that ${\bar
Z}(B(x_0,r/\tau),\bar{u})<\infty$. If $u\geq 0$ in $B(x_0,r)$, then
 \bea\label{alphaneg}
||\bar{u}^{\alpha}||_{s,B(x_0,r);\overline{dx}}^\frac{1}{\alpha}&\leq&C_{10}\Big(C_{11}\bar{Z}(B(x_0,r),\bar{u})
\Big)^\frac{p\Psi_0}{|\alpha|}\disp\teinf_{B(x_0,\tau r)}\bar{u},
\eea where the constants $C_{10},C_{11}$ can be chosen as in
\eqref{qpos}.
\end{lem}

\noindent{\bf Proof:} Since by our assumptions
$\bar{Z}(B(x_0,r/\tau),\bar{u})<+\infty$, $B(x_0,r)\subset
B(x_0,r/\tau)$ and $r/\tau<\tau r_1(x_0)<R_1(x_0)$, we have that
$\bar{Z}:=\bar{Z}(B(x_0,r),\bar{u})$ is also finite. See the comment
following the definition \eqref{Zetabar} of $\bar{Z}$.

Following the same argument as in Lemma \ref{alpha>0} but now with
$q<1-p$ and using that $(u,\nabla u)$ is a weak solution of
(\ref{supersol1}), applying Remark \ref{diffineqcomp} we obtain the
following inequality, similar to (\ref{PHIit.1}): \bea
\Big[\Phi(j+1;Y\sigma)\Big]^{\frac{Y}{p}} \leq
C\Big(|q|+\frac{1}{|q|}\Big)^{\tilde{b}_*}\bar{Z}N^j\Big[\Phi(j;sY)
\Big]^\frac{Y}{p},\nonumber \eea with
$\bar{Z}=\bar{Z}(B(x_0,r),\bar{u})$. Since $Y=p+q-1<0$ for any $q\in
(-\infty,1-p)$, we have \bea\label{PHIit.02} \Phi(j+1;Y\sigma) \geq
C^\frac{p}{Y}\Big(|q|+\frac{1}{|q|}\Big)^\frac{p\tilde{b}_*}{Y}
\bar{Z}^\frac{p}{Y}N^\frac{jp}{Y}\Phi(j;sY). \eea  Let $\alpha<0$,
set $Y_j = \alpha\X^j$ and $q_j=\alpha\X^j+1-p$, where
$\X=\sigma/s>1$ as in Lemma \ref{alpha>0}. Then $Y_j<0$,
$q_j<-(p-1)$ and, with $C$ as in \eqref{PHIit.1}, we have
\bea\label{PHIit.02bis} \Phi(j;s\alpha\X^j) \leq
\Big[C^{\X^{-j}}\Big(|q_j|+\frac{1}{|q_j|}\Big)^{\tilde{b}_*\X^{-j}}
\bar{Z}^{\X^{-j}}N^{j\X^{-j}}\Big]^\frac{p}{|\alpha|}\Phi(j+1;s\alpha\X^{j+1}).
\eea Now note that
\begin{eqnarray*}
|q_j|+\frac{1}{|q_j|}&\leq&|\alpha|\X^j+p-1+\frac{1}{p-1}\,\,<\,\,|\alpha|\X^j+p-1+\frac{1}{(p-1)(1-\X^{-\frac{1}{4}})}\\
&\leq&\X^j\left[|\alpha|\X^\frac{1}{2} + p -1  +
\frac{1}{(p-1)(1-\X^{-\frac{1}{4}})}\right]\,\,=\,\,\X^j\left[|\alpha|\X^\frac{1}{2}
+ c_4\right],
\end{eqnarray*}
with $c_4=p -1  + ((p-1)(1-\X^{-\frac{1}{4}}))^{-1}$ as in Lemma
\ref{alpha>0}. Then from \eqref{PHIit.02bis} we have
\begin{eqnarray}
\nonumber\Phi(j;\alpha s\X^j)
&\leq&\Big[C^{\X^{-j}}(|\alpha|\X^\frac{1}{2}+c_4)^{\tilde{b}_*\X^{-j}}\bar{Z}^{\X^{-j}}
\X^{\tilde{b}_*j\X^{-j}}N^{j\X^{-j}}\Big]^\frac{p}{|\alpha|}
\Phi(j+1;\alpha s\X^{j+1}).
\end{eqnarray}
Iterating the previous inequality we obtain
\begin{eqnarray}
\nonumber\Phi(j;\alpha s\X^j)&\leq&
\Big[C^{\Psi_0}(|\alpha|\X^{\frac{1}{2}}+c_4)^{\tilde{b}_*\Psi_0}
\X^{\tilde{b}_*\Psi_1} N^{\Psi_1}
\bar{Z}^{\Psi_0}\Big]^\frac{p}{|\alpha|}\Phi(\infty;-\infty),
\end{eqnarray}
where $\Phi(\infty;-\infty) = \disp\limsup_{j\ra\infty}\Phi(j;\alpha
s\X^j)$, $\Psi_0=\sum_{j=0}^\infty\X^{-j}$ and
$\Psi_1=\sum_{j=0}^\infty j\X^{-j}$. Also using \eqref{22} we
conclude
\begin{equation*}
\Phi(j;\alpha s\X^j) \leq
C_{10}\Big[C_{11}\bar{Z}\Big]^\frac{\Psi_0p}{|\alpha|}\Phi(\infty,-\infty),
\end{equation*}
with $C_{10},C_{11}$ as in \eqref{PHIit.3}. Since this holds for all
$j\in\N\cup\{0\}$, we obtain \bea \Phi(0,\alpha s)&\leq&C_{10}
\Big[C_{11}\bar{Z}\Big]^\frac{\Psi_0p}{|\alpha|}\Phi(\infty,-\infty).
\eea Since $B(x_0,\tau r)\subset S_j$ for every $j\ge1$ we have
$\disp\text{ess\,inf}_{B(x_0,\tau r)} \bar{u} \geq
\Phi(\infty;-\infty)$, see \cite{GT}; hence \bea
||\bar{u}^{\alpha}||^{\frac{1}{\alpha}}_{s,B(x_0,r);\overline{dx}}
&\leq& C_{10}\Big[C_{11}\bar{Z}\Big]^\frac{\Psi_0p}{|\alpha|}
\disp\einf_{B(x_0,\tau r)}\bar{u}, \eea which proves
(\ref{alphaneg}).
\begin{flushright} $\Box$
\end{flushright}

\begin{lem}\label{logestimate}  Let the assumptions in the opening
paragraph of this section hold, and suppose that $(u,\nabla u) \in
W^{1,p}_Q(\Omega)$ is a weak solution in $\Omega$ of \bea
\label{supersol}\text{\emph{div}}(A(x,u,\nabla u)) \leq B(x,u,\nabla
u) \eea where $A,B$ satisfy \eqref{struct} with exponents
$\delta,\gamma,\psi$ as in \eqref{ranges1}. Furthermore, suppose
that the Poincar\'e inequality \eqref{poincare} holds.  Fix
$\hat{x}\in\Omega$, $k>0$,
$\epsilon_1,\epsilon_2,\epsilon_3\in(0,1]$ and let ${\bar u}=|u|+k$
and $w=\log {\bar u}$. Fix a quasimetric ball $B(\hat{x},\gb
l/\tau)$ for $\gb>1$ as in \eqref{poincare} and $0<l<\tau
r_1(\hat{x})/\gb$. If $u\geq 0$ in $B(\hat{x},\gb l/\tau)$, then
\bea\label{logubar} \fint_{B(\hat{x},l)}|w-w_{B(\hat{x},l)}|\,dx
\leq C_{12}\bar{Z}(B(\hat{x},\gb l/\tau),\bar{u}), \eea  where
${\bar Z}\big(B(\hat{x},\gb l/\tau), {\bar u}\big)$ may be infinite
and where $C_{12}$ depends on $a,p,\sigma$, on
$\epsilon_1,\epsilon_2,\epsilon_3$ in the definition \eqref{Zetabar}
of $\bar{Z}$, on $\mathfrak{b},C_2$ in \eqref{poincare}, on
$d_0,C_0$ in \eqref{doubling1} and  on $C_{s^*},\tau,N$ in
\eqref{cutoff}.
\end{lem}

\noindent {\bf Proof:} We can assume that $\bar{Z}:={\bar
Z}\big(B(\hat{x},\gb l/\tau), {\bar u}\big)$ is finite, otherwise
\eqref{logubar} is trivial. Let $\eta=\eta_1$ be as in
(\ref{cutoff}) relative to $B(\hat{x},\gb l/\tau)$, and set
$v=\eta^p\bar{u}^{1-p}$. Applying Remark \ref{diffineqcomp} on the
quasimetric ball $B(\hat{x},\gb l/\tau)$, we have that $v\in
W^{1,p}_{Q,0}(\Omega)$ with $\text{supp}\,v\subset B(\hat{x},\gb
l/\tau)$ and using \eqref{struct2} we have \bea\label{bridge.1}
\nabla v\cdot A +v B &\leq& a^{-1}(1-p)\eta^p\bar{u}^{-p}|\sqrt{Q}\nabla\bar{u}|^p + (p-1)\bar{h}\eta^p\nonumber\\
&& +ap\eta^{p-1}|\sqrt{Q}\nabla \eta|\bar{u}^{1-p}|\sqrt{Q}\nabla\bar{u}|^{p-1} + p\eta^{p-1}\bar{b}|\sqrt{Q}\nabla \eta|\\
&& +c\eta^p\bar{u}^{1-p}|\sqrt{Q}\nabla\bar{u}|^{\psi-1} +
\eta^p\bar{d}\nonumber \eea a.e. in $\Omega$. We integrate over
$B(\hat{x},\gb l/\tau)$ and use the facts that $(u,\nabla u)$ is a
weak solution of (\ref{supersol}) and that $v$ is a feasible
nonnegative test function, obtaining that the left side of the
resulting inequality is nonnegative. Also, we move the resulting
first term on the right side to the left side and estimate the third
and fifth terms on the right in ways like those used to estimate
similar terms in (\ref{q>0.1}). Then we obtain, as in \eqref{q>0.4}
with $q=1-p$, \bea\label{bridge.2}
\;\;\;\;\fint_{B(\hat{x},\frac{\gb l}{\tau})}
\eta^p\bar{u}^{-p}|\sqrt{Q}\nabla\bar{u}|^pdx &\leq& C\Big\{
\fint_{B(\hat{x},\frac{\gb l}{\tau})}|\sqrt{Q}\nabla\eta|^pdx
+\fint_{B(\hat{x},\frac{\gb l}{\tau})}\bar{h}\eta^pdx\\
&&\,\,+\fint_{B(\hat{x},\frac{\gb l}{\tau})}\eta^{p-1}\bar{b}
|\sqrt{Q}\nabla\eta|\;dx + \fint_{B(\hat{x},\frac{\gb
l}{\tau})}\eta^p\bar{d}\;dx\nonumber\\
&&\,\,+\fint_{B(\hat{x},\frac{\gb l}{\tau})}
c^\frac{p}{p+1-\psi}\eta^p\bar{u}^\frac{p(\psi-p)}{p+1-\psi}dx
\Big\},\nonumber \eea with $C$ depending only on $a,p,\sigma$.
Repeating steps 5 and 6 in the proof of \cite[Theorem 1.2]{MRW} we
obtain \bea\label{bridge.3} \fint_{B(\hat{x},\frac{\gb l}{\tau})}
\eta^p\bar{u}^{-p}|\sqrt{Q}\nabla\bar{u}|^pdx &\leq&
C\bar{Z}^p\Big\{\fint_{B(\hat{x},\frac{\gb l}{\tau})}
|\sqrt{Q}\nabla\eta|^pdx + \frac{1}{l^p}\fint_{B(\hat{x},\frac{\gb l}{\tau})} \eta^p dx\Big\}\\
&\leq& C\bar{Z}^p\Big\{\fint_{B(\hat{x},\frac{\gb
l}{\tau})}|\sqrt{Q}\nabla\eta|^pdx + \frac{1}{l^p}\Big\}\nonumber
\eea analogous to \cite[(3.21)]{MRW} with $Y=0$, noting that $0\leq
\eta\leq 1$ on $B(\hat{x},\frac{\gb l}{\tau})$.  We recall that here
$\bar{Z}=\bar{Z}(B(\hat{x},\gb l/\tau),{\bar u})$ and we note that
$C$ depends on $a,p,\sigma$ and on
$\epsilon_1,\epsilon_2,\epsilon_3$ appearing in the definition
\eqref{Zetabar} of $\bar{Z}$.

Since $\eta$ is the function $\eta_1$ in (\ref{cutoff}) relative to
$B(\hat{x},\gb l/\tau)$, then $\eta\in
Lip_0(B(\hat{x},\gb l/\tau))\cap Lip(\Omega)$ and
$\eta\equiv 1$ on $B(\hat{x},\gb l)$. Recalling that $\gb
l<r_1(y)$, we apply the Poincar\'e
inequality (\ref{poincare}) to $w=\log \bar{u}$ (see Remark
\ref{diffineqcomp}) and get
\bea\label{bridge.4}
\fracc{1}{|B(\hat{x},l)|}\int_{B(\hat{x},l)}|w-w_{B(\hat{x},l)}|dx
&\leq& Cl \Big(\fracc{1}{|B(\hat{x},{\gb l})|}\int_{B(\hat{x},{\gb
l})}|\sqrt{Q}\nabla w|^pdx\Big)^{1/p}\nonumber\\
&\leq& C l \Big(\fracc{1}{|B(\hat{x},{\gb l})|}\int_{B(\hat{x},\gb l/\tau)}\eta^p|\sqrt{Q}\nabla w|^pdx\Big)^{1/p}\nonumber\\
&=& Cl\fracc{|B(\hat{x},\gb l/\tau)|^{1/p}}{|B(\hat{x},{\gb l})|^{1/p}} \Big(\fint_{B(\hat{x},\gb l/\tau)}\eta^p\bar{u}^{-p}
|\sqrt{Q}\nabla \bar{u}|^pdx\Big)^{1/p}\nonumber\\
&\leq& Cl\Big(\fracc{|B(\hat{x},\gb l/\tau)|}{|B(\hat{x},{\gb
l})|}\Big)^{\frac{1}{p}}\bar{Z}\Big[\Big(\fint_{B(\hat{x},\gb
l/\tau)}|\sqrt{Q}\nabla \eta|^pdx\Big)^\frac{1}{p}
+\frac{1}{l}\Big]
\eea
where the last line is obtained using \eqref{bridge.3}.  Also, $C$ in \eqref{bridge.4} depends on $a,p,\sigma,\epsilon_1,\epsilon_2,\epsilon_3$ and on
the constants $\mathfrak{b},C_2$ appearing in \eqref{poincare}. In
(\ref{bridge.4}), use H\"older's inequality with exponents
$\frac{s^*}{p},\frac{s^*}{s^*-p}$ together with (\ref{doubling1})
and (\ref{cutoff}) to obtain
\bea\label{bridge.5}
\fracc{1}{|B(\hat{x},l)|}\int_{B(\hat{x},l)}|w-w_{B(\hat{x},l)}|dx
&\leq& Cl \bar{Z}\Big(\fracc{|B(\hat{x},\gb
l/\tau)|}{|B(\hat{x},{\gb l})|}\Big)^{1/p}
\Big[\Big(\fint_{B(\hat{x},\gb l/\tau)}|\sqrt{Q}\nabla\eta|^{s^*}dx\Big)^\frac{1}{s^*} + \frac{1}{l}\Big]\nonumber\\
&\leq& C\bar{Z}\tau^{-\frac{d_0}{p}}\Big[\fracc{\tau N}{\gb
}+1\Big]\nonumber\\ &=& C_{12}\bar{Z},
\eea where $C_{12}$ depends
on $a,p,\sigma,\epsilon_1,\epsilon_2,\epsilon_3,\mathfrak{b},C_2$,
on the constants $d_0,C_0$ in \eqref{doubling1} and $C_{s^*},\tau,N$
in \eqref{cutoff}.
\begin{flushright}
$\Box$
\end{flushright}


\noindent {\bf Proof of Proposition \ref{harnackmain}:} We will use the
notation and assumptions of Proposition \ref{harnackmain} and divide the
proof into steps.

\emph{Step 1.} We have $B(x_0,C_*r)\Subset
B(y,\frac{\tau}{2}r_1(y))$. Indeed if $\xi\in B(x_0,C_*r)$, then
$$\rho(\xi,y)\leq\kappa(\rho(\xi,x_0)+\rho(x_0,y))<\kappa\Big(C_*r+\frac{\tau}{5\kappa}r_1(y)\Big)
<\Big(\frac{\tau
A_*}{5}+\frac{\tau}{5}\Big)r_1(y)<\frac{\tau}{2}r_1(y).$$

\emph{Step 2.} By using Lemmas \ref{alpha>0} and \ref{alpha<0}, let
us show that for every $\alpha>0$ there exists
$\alpha_1\in[\alpha\sigma^{-1/2},\alpha]$ such that
\begin{equation}\label{6}
\esup_{B(x_0,\tau r)}\bar{u}\leq
C_{10}\big[C_{11}\bar{Z}\big(B(x_0,r),\bar{u}\big)\big]^\frac{p\psi_0}{\alpha_1}\|\bar{u}^{\alpha_1}\|_{s,B(x_0,r);\overline{dx}}^\frac{1}{\alpha_1}
\end{equation}
and that for every $\alpha_2<0$
\begin{equation}\label{23}
\|\bar{u}^{\alpha_2}\|_{s,B(x_0,r);\overline{dx}}^{\frac{1}{\alpha_2}}\leq
C_{10}\big[C_{11}\bar{Z}\big(B(x_0,r),\bar{u}\big)\big]^\frac{p\psi_0}{|\alpha_2|}\einf_{B(x_0,\tau
r)}\bar{u},
\end{equation}
with $C_{10},C_{11}$ and $\psi_0$ independent of $(u,\nabla
u),k,B(x_0,r),y,b,c,d,e,f,g,h,\alpha,\alpha_1,\alpha_2$.

Indeed, by our assumptions, $r_1$ satisfies a local uniformity
condition with respect to $\rho$ on $B=B(y,r_1(y))$ with constant
$A_*$. Since $x_0\in B(y,r_1(y))$ we have $r_1(x_0)>A_*r_1(y)$, so
that $r< \frac{\tau A_*}{5\kappa C_*}r_1(y)<\tau^2r_1(x_0)$.
Moreover $\tau^{-1}<C_*$, so that $r/\tau<C_*r$ and $B(x_0,
r/\tau)\subset B(x_0,C_*r)$. Thus by (\ref{zbarcond}) we conclude
that $\bar{Z}(B(x_0,r/\tau),\bar{u})\leq M<+\infty$, and hence all
the assumptions of Lemmas \ref{alpha>0} and \ref{alpha<0} are
satisfied.

\emph{Step 3.} We start implementing the ideas of Section \ref{sec5}. Let
$$R(\xi)=\min\left\{\frac{16(\gamma^*)^4\kappa^5}{A_*^2\min\{A_*^2,
(8\kappa^5)^{-1}\}}r,\frac{\tau^2A_*}{40\kappa^6(\gamma^*)^4
\mathfrak{b}}r_1(\xi)\right\}\,\,\,\text{for }\xi\in\Omega, \quad
\text{and let }B=B(x_0,r).$$ Let us show that if $\xi\in
B(y,r_1(y))$ then
$R(\xi)=\frac{16(\gamma^*)^4\kappa^5}{A_*^2\min\{A_*^2,
(8\kappa^5)^{-1}\}}r$. Indeed for every $\xi\in B(y,r_1(y))$ we have
$A_*r_1(y)<r_1(\xi)$ by our assumptions on $r_1$. Since
$r\in(0,\frac{\tau A_*}{5\kappa C_*}r_1(y))$ with $C_*$ as in
\eqref{5}, we obtain
$$\frac{16(\gamma^*)^4\kappa^5}{A_*^2\min\{A_*^2,
(8\kappa^5)^{-1}\}}r<\frac{16\tau(\gamma^*)^4\kappa^4}{5C_*A_*\min\{A_*^2,
(8\kappa^5)^{-1}\}}r_1(y)= \frac{\tau^2A_*^2}{40\kappa^6(\gamma^*)^4
\mathfrak{b}}r_1(y)<\frac{\tau^2A_*}{40\kappa^6(\gamma^*)^4
\mathfrak{b}}r_1(\xi).$$ Note that by the second restriction above
on $R(\xi)$, we have $R(\xi)<\frac{r_1(\xi)}{(\gamma^*)^4} <
\frac{r_1(\xi)}{(\gamma^*)^2}$ for every  $\xi\in\Omega$, which
meets some of the requirements of Proposition \ref{JN}. Moreover,
since $R(\xi)$ is constant on $B$, it satisfies a local uniformity
condition on $B$ with respect to $\rho$ for any constant in $(0,1]$,
in particular for $A_*$. Thus $R(\xi)$ satisfies the requirements in
the statement of Proposition \ref{JN} and consequently can be used
in Corollary \ref{JNcor}. Hence there are constants $C_8,C_9$,
$c_\rho=8\kappa^5(\gamma^*)^2$,
$\delta_0=\frac{A_*^2\min\{A_*^2,(8\kappa^5)^{-1}\}}{8\kappa^5(\gamma^*)^3}$
(cf. Remark 5.2 for the values of $c_\rho$ and $\delta_0$) such that
\eqref{JNcor1} holds for every function $f\in [c_\rho R]BMO(B)$ with
$\|f\|_{[c_\rho R]BMO(B)}<C_9$.

\emph{Step 4.} We claim that if $\xi\in B$ and $B(z,t_0)$ is a
$c_\rho R(\xi)$--ball then
\begin{equation}\label{7}
\begin{array}{l}
\displaystyle i)\quad B(z,t_0)\subset
B\Big(z,\frac{\mathfrak{b}}{\tau}t_0\Big)\subset
B(x_0,C_*r),\\
\displaystyle ii)\quad
0<t_0<\frac{\mathfrak{b}}{\tau}t_0<C_*r<\frac{\tau
A_*}{5\kappa}r_1(y)<r_1(z)
\end{array}
\end{equation}
Since $\xi\in B$ and $B(z,t_0)$ is a $c_\rho R(\xi)$--ball (see the
definition above (\ref{Rbmo})), we have
\begin{equation}\label{8}
t_0<c_\rho
R(\xi)=\frac{128(\gamma^*)^6\kappa^{10}}{A_*^2\min\{A_*^2,(8\kappa^5)^{-1}\}}r=\frac{2\gamma^*c_\rho}{\delta_0}r
\end{equation}
and $B(z,\gamma^*t_0)\subset B(\xi,c_\rho
R(\xi))=B\Big(\xi,\frac{2\gamma^*c_\rho}{\delta_0}r\Big)$. Thus,
also using the swallowing property of the pseudometric balls as
described in Lemma \ref{swalemma},
\begin{eqnarray*}
B(z,t_0)\!&\!\subset\!&\!
B\Big(z,\frac{\mathfrak{b}}{\tau}t_0\Big)\subset
B\Big(z,\frac{2\gamma^*c_\rho\mathfrak{b}}{\delta_0\tau}r\Big)\subset
B\Big(\xi,\frac{2(\gamma^*)^2c_\rho\mathfrak{b}}{\delta_0\tau}r\Big)\\
\!&\!\subset\!&\!B\Big(x_0,\frac{2(\gamma^*)^3
c_\rho\mathfrak{b}}{\delta_0\tau}r\Big)=B\Big(x_0,C_*A_*r\Big)\subset
B(x_0,C_*r).
\end{eqnarray*}
Since by Step 1 we have $B(x_0,C_*r)\subset B(y,r_1(y))$ and since
$r_1$ satisfies a local uniformity condition with respect to $\rho$
with constant $A_*$ on that set, we conclude that
$A_*r_1(y)<r_1(z)$. Hence, also using \eqref{8},
$$0<t_0<\frac{\mathfrak{b}}{\tau}t_0<\frac{2\mathfrak{b}\gamma^*c_\rho}{\tau\delta_0}r=\frac{A_*}{(\gamma^*)^2}C_*r<C_*r<\frac{\tau A_*}{5\kappa}r_1(y)<r_1(z).$$
Here the next-to-last inequality is due to the relation between $r$
and $r_1(y)$ that we noted earlier. The proof of our claim is complete.

\emph{Step 5.} Now we show that $w=\log\bar{u}$ is a function in
$[c_\rho R]BMO(B)$, where as above $B= B(x_0,r)$. Let $\xi\in B$ and
consider a $c_\rho R(\xi)$--ball $B(z,t_0)$. By Step 4, Lemma
\ref{logestimate} and condition \eqref{zbarcond} we conclude that
$$\fint_{B(z,t_0)} |w(\zeta)-w_{B(z,t_0)}|\,d\zeta\,\,\leq\,\,
C_{12}\overline{Z}\Big(B\Big(z,\frac{\mathfrak{b}}{\tau}t_0\Big),
\bar{u}\Big)\,\,\leq\,\, C_{12}M,$$ and thus, by the definition
given in \eqref{Rbmo}, we have $\|w\|_{[c_\rho R]BMO(B)}\leq
C_{12}M$.

Now we choose $\alpha=\frac{C_9}{2sC_{12}M}$, where $C_9$ is as in
Proposition \ref{JN} and where $s$ is as in Lemmas \ref{alpha>0} and
\ref{alpha<0} and also \eqref{6}. Then the corresponding $\alpha_1$
from inequality \eqref{6} satisfies
$\alpha_1\in[\alpha\sigma^{-1/2},\alpha]$.
Then we have $$\|s\alpha_1w\|_{[c_\rho R]BMO(E)}\leq
s\alpha_1C_{12}M\leq s\alpha C_{12}M\leq \frac{C_9}{2},$$ and by
Step 3 we can use Corollary \ref{JNcor} to conclude that
\begin{equation}\label{9}
\|e^{s\alpha_1w}\|_{[\delta_0R]A_2(B)}\leq(1+C_8)^2.
\end{equation}

\emph{Step 6.} We notice here that $x_0\in B(x_0,r)$ and that
$B(x_0,r)$ is a $\delta_0R(x_0)$--ball, and use this fact in
conjunction with \eqref{9}. We start by recalling that since $x_0\in
B(y,r_1(y))$, Step 3 shows that
$R(x_0)=\frac{16(\gamma^*)^4\kappa^5}{A_*^2\min\{A_*^2,(8\kappa^5)^{-1}\}}r$.
Now a simple calculation gives
$0<r<\gamma^*r<2\gamma^*r=\delta_0R(x_0)$, $B(x_0,\gamma^*r)\subset
B(x_0,2\gamma^*r)=B(x_0,\delta_0R(x_0))$ and
$\overline{B(x_0,\gamma^*r)}\subset
\overline{B(x_0,r_1(x_0))}\subset\Omega$.

Since $B(x_0,r)$ is a $\delta_0R(x_0)$--ball with $x_0\in B(x_0,r)$, by
\eqref{9} and definition \eqref{RA2} we have
\begin{equation*}
\left(\fint_{B(x_0,r)}e^{s\alpha_1w}\,d\zeta\right)\left(\fint_{B(x_0,r)}e^{-s\alpha_1w}\,d\zeta\right)\leq(1+C_8)^2,
\end{equation*}
and thus we conclude
\begin{equation}\label{10}
\|\bar{u}^{\alpha_1}\|^\frac{1}{\alpha_1}_{s,B(x_0,r);\overline{dx}}\leq(1+C_8)^\frac{2}{s\alpha_1}\|\bar{u}^{-\alpha_1}\|^{-\frac{1}{\alpha_1}}_{s,B(x_0,r);\overline{dx}}
\end{equation}

\emph{Step 7.} Now we use \eqref{6}, \eqref{23} with
$\alpha_2=-\alpha_1<0$ and \eqref{10} to finish the proof:
\begin{eqnarray*}
\text{ess}\!\!\!\sup_{\!\!\!\!B(x_0,\tau
r)}\bar{u}&\leq&C_{10}\big[C_{11}\bar{Z}\big(B(x_0,r),\bar{u}\big)\big]^\frac{p\psi_0}{\alpha_1}\|\bar{u}^{\alpha_1}\|_{s,B(x_0,r);\overline{dx}}^\frac{1}{\alpha_1}\\
&\leq&C_{10}(1+C_8)^\frac{2}{s\alpha_1}\big[C_{11}\bar{Z}\big(B(x_0,r),\bar{u}\big)\big]^\frac{p\psi_0}{\alpha_1}\|\bar{u}^{-\alpha_1}\|_{s,B(x_0,r);\overline{dx}}^{-\frac{1}{\alpha_1}}\\
&\leq&C_{10}^2(1+C_8)^\frac{2}{s\alpha_1}\big[C_{11}\bar{Z}\big(B(x_0,r),\bar{u}\big)\big]^\frac{2p\psi_0}{\alpha_1}\text{ess}\!\!\!
\inf_{\!\!\!\!\!B(x_0,\tau r)}\bar{u}.
\end{eqnarray*}
Since $C_8+1\geq 1$, $\bar{Z}\big(B(x_0,r),\bar{u}\big)\geq1$ and
$\alpha_1\geq\alpha\sigma^{-1/2}$, we see that
$$\text{ess}\!\!\!\sup_{\!\!\!\!B(x_0,\tau
r)}\bar{u}\leq
C_{10}^2\big[(1+C_8)^\frac{1}{sp\Psi_0}C_{11}\bar{Z}\big(B(x_0,r),\bar{u}\big)\big]^\frac{2p\psi_0\sqrt{\sigma}}{\alpha}\text{ess}\!\!\!
\inf_{\!\!\!\!\!B(x_0,\tau r)}\bar{u}$$
 which, recalling the definition of $\alpha$ given in Step 5, is inequality
\eqref{Harngoal}, with $C_4=C_{10}^2$,
$C_5=(1+C_8)^\frac{1}{sp\Psi_0}C_{11}$ and $C_6={4\sqrt\sigma
p\Psi_0sC_{12}}/{C_9}$.
\begin{flushright}$\Box$
\end{flushright}


\section{The Proof of Theorem \ref{holder1}}

Let $(u,\nabla u)$ be a weak solution of \eqref{eqdiff} in $\Omega$
and let $x_0\in B=B(y,\frac{\tau^2}{5\kappa }r_1(y))$. For $r>0$
(sufficiently small so that $\overline{B(x_0,r)}\subset\Omega$), define
\bea M(r) = \disp\esup_{B(x_0,r)} u,\;\;  m(r) = \disp\einf_{
B(x_0,r)} u, \text{ and } \omega_{x_0}(r)=M(r)-m(r).\nonumber
\eea
We will refer to $\omega_{x_0}(r)$ as the oscillation of $u$ in
$B(x_0,r)$.  Now, let $r\in (0,\frac{\tau^2 A_*}{5\kappa C_*}r_1(y))$,
where $C_*$ is as in (\ref{5}), and set $M_0=M(C_*r)$, $m_0=m(C_*r)$,
noting that $M_0$ and $m_0$ are finite by \cite[Theorem 1.2]{MRW} and
Proposition \ref{pro1} as $B(x_0,C_*r)\subset B(y,\tau^2
r_1(y))$. Denote
\bea (u_1,\nabla u_1)=(M_0-u,-\nabla u) &\text{ and }& (u_2,\nabla
u_2)=(u-m_0,\nabla u).\nonumber
\eea
Clearly $(u_1,\nabla u_1), (u_2,\nabla u_2) \in W^{1,p}_Q(\Omega)$
and $u_1, u_2 \ge 0$ in $B(x_0,C_*r)$.  For $(x,z,\xi)\in \Omega\times
\mathbb{R}\times\mathbb{R}^n$, let
\bea \label{alternates}
{A_1}(x,z,\xi)=-A(x,M_0-z,-\xi),&& {A_2}(x,z,\xi) =
A(x,z+m_0,\xi),\nonumber\\
\tilde{{A_1}}(x,z,\xi) = -\tilde{A}(x,M_0-z,-\xi),&&
\tilde{{A_2}}(x,z,\xi)=\tilde{A}(x,z+m_0,\xi),\text{ and}\nonumber\\
{B_1}(x,z,\xi)=-B(x,M_0-z,-\xi),&& {B_2}(x,z,\xi) =
B(x,z+m_0,\xi).\nonumber
\eea
It is not difficult to see that $(u_1,\nabla u_1)$ and $(u_2,\nabla
u_2)$ are respectively weak solutions in $\Omega$ of
\bea\label{altdiffeqs}
\text{div}\Big({A_1}(x,u,\nabla u)\Big) &=& {B_1}(x,u,\nabla u),\text{ and}\\
\text{div}\Big({A_2}(x,u,\nabla u)\Big) &=& {B_2}(x,u,\nabla u).\nonumber
\eea

We now check that these equations satisfy $\eqref{struct}$ for
coefficients which satisfy $(i)-(iv)$ of Proposition \ref{pro1}.  The
calculations are simple and we only provide an example.  Let us show
that ${A_1}$ satisfies item (ii) of \eqref{struct} for appropriately
modified definitions of $g, h$.  Indeed, since $A$ satisfies
\eqref{struct}, we have
\bea \xi\cdot {A_1}(x,z,\xi) &=&-\xi \cdot A(x,M_0-z,-\xi)\nonumber\\
&\geq& a^{-1}|\sqrt{Q(x)}\xi|^p-h(x)|M_0-z|^p-g(x)\nonumber\\
&\geq& a^{-1}|\sqrt{Q(x)}\xi|^p -2^{p-1}h(x)|z|^p
-(g(x)+2^{p-1}h(x)|M_0|^p).\nonumber
\eea
Setting ${h_1}(x)=2^{p-1}h(x)$ and ${g_1}(x)=g(x)+2^{p-1}|M_0|^ph(x)$
it follows that ${A_1}(x,z,\xi)$ satisfies \eqref{struct}(ii) with
$h,g$ there replaced by $h_1,g_1$ respectively.  Furthermore,
$h_1, g_1 \in L^\mathcal{H}_{loc}(\Omega)$ since $|M_0|<\infty$ and
both $h, g \in L^\mathcal{H}_{loc}(\Omega)$ by hypothesis.  Other
verifications are similar using the modified functions
\bea {h_2}(x)=2^{p-1}h(x),&&
{g_2}(x)=g(x)+2^{p-1}|m_0|^ph(x),\nonumber\\
{e_1}(x)=e(x)+2^{p-1}|M_0|^{p-1}b(x),&& {e_2}(x)=e(x)+2^{p-1}|m_0|^{p-1}b(x),\text{ and}\nonumber\\
{f_1}(x)=f(x)+2^{p-1}|M_0|^{p-1}d(x),&& {f_2}(x)=f(x)+2^{p-1}|m_0|^{p-1}d(x)\nonumber
\eea
with
\bea
{b_1}(x)\;\;=&b_2(x)&=\;\;2^{p-1}b(x),\nonumber\\
{c_1}(x)\;\;=&{c_2}(x)&=\;\;c(x),\text{ and}\nonumber\\
{d_1}(x)\;\;=&{d_2}(x)&=\;\;2^{p-1}d(x).\nonumber
\eea

Therefore both $(u_1,\nabla u_1)$ and $(u_2,\nabla u_2)$ are weak
solutions of equations satisfying the hypotheses of Theorem
\ref{harnackmother}.  As a consequence, $u_1,u_2$ satisfy
\bea\label{harnackA1} \disp\esup_{z\in B(x_0,\tau r)}u_1(z) +
k_1(x_0,r) &\leq& \hat{C}_1\Big[\disp\einf_{z\in B(x_0,\tau
r)}u_1(z)+k_1(x_0,r)\Big],\text{
and}\\
\label{harnackA2} \disp\esup_{z\in B(x_0,\tau r)}u_2(z) + k_2(x_0,r)
&\leq& \hat{C}_2\Big[\disp\einf_{z\in B(x_0,\tau
r)}u_2(z)+k_2(x_0,r)\Big]. \eea Here $k_1=k_1(x_0,r)$ and
$k_2=k_2(x_0,r)$ are defined as $k$ in Proposition \ref{pro1} using
the structural coefficient functions
${b_1},{c_1},{d_1},{e_1},{f_1},{g_1},{h_1}$ and
${b_2},{c_2},{d_2},{e_2},{f_2},{g_2},{h_2}$ respectively.  By
Proposition \ref{pro1}, \bea\label{kjbound} k_j(x_0,r) \leq
\Lambda_j r^\lambda, \quad j=1,2, \eea with $\lambda$ exactly as in
Proposition \ref{pro1} and $\Lambda_1,\Lambda_2$ defined as
$\Lambda$ in Proposition \ref{pro1} using instead the structural
coefficients related to ${A_1},{A_2},{B_1},{B_2}$. Each of
$\hat{C}_1,\hat{C}_2, \Lambda_1,\Lambda_2$ depends on $p$, $\psi$,
$M_0$, $m_0$, $||u||_{p\sigma,\tilde{B};dx}$,
$||b||_{\mathcal{B},\tilde{B};dx}$,
$||c||_{\mathcal{C},\tilde{B};dx}$,
$||d||_{\mathcal{D},\tilde{B};dx}$,
$||e||_{\mathcal{B},\tilde{B};dx}$,
$||f||_{\mathcal{D},\tilde{B};dx}$,
$||g||_{\mathcal{H},\tilde{B};dx}$,
$||h||_{\mathcal{H},\tilde{B};dx}$, $C_0$, $d_0$, $s$, $a$, and $N$,
where $\tilde{B}=B(y,r_1(y))$.  It is important to also note that
$\lambda$ is independent of $u$, and that when
$\psi\in(p,p+1-\sigma^{-1})$ the dependence of
$\hat{C}_1,\hat{C}_2,\Lambda_1,\Lambda_2$ on $||u||_{p\sigma,B;dx}$
occurs through $M_1,M_2$ of Proposition \ref{pro1} and through
$M_0$, $m_0$, see \cite[Theorem 1.2]{MRW}. Moreover
$\hat{C}_1,\hat{C}_2,M_1,M_2$ are independent of
$||u||_{p\sigma,B;dx},M_0,m_0$ when $\psi=p$.

Setting $C=\max\{\hat{C}_1,\hat{C}_2\}$ and rewriting
\eqref{harnackA1} and \eqref{harnackA2} in terms of $M(r)$ and
$m(r)$ gives \bea \label{hars}M_0-m(\tau r) +k_1(x_0,r) &\leq& C
\Big(M_0- M(\tau
r) + k_1(x_0,r)\Big),\text{ and}\\
M(\tau r) - m_0 +k_2(x_0,r) &\leq& C\Big( m(\tau r) - m_0
+k_2(x_0,r)\Big).\nonumber
\eea
Adding the inequalities in \eqref{hars}, rearranging and inserting the
oscillation $\omega_{x_0}$, we obtain
\bea \omega_{x_0}(\tau r)(C+1) \leq (C-1)\Big(\omega_{x_0}(C_*r) +
(k_1+k_2)\Big)\nonumber
\eea
and so
\bea \label{GT823soon}
\omega_{x_0}(\tau r) \leq \frac{C-1}{C+1}\Big(\omega_{x_0}(C_*r) +
\hat{\Lambda}r^\lambda\Big),\quad 0<r < \frac{\tau A_*}{5\kappa
C_*}r_1(y),
\eea
where we have used estimate \eqref{kjbound} and set $\hat{\Lambda} =
\Lambda_1+\Lambda_2$.  Define $R=C_*r$, $R_0 = \frac{\tau^2
A_*}{6\kappa}r_1(y)$ and $\hat{\Lambda}_0=
\hat{\Lambda}C_*^{-\lambda}$.  Recall that $\tau < 1 <C_*$. Then
\eqref{GT823soon} and the monotonicity of $\omega_{x_0}$ imply that
for every $\nu\leq \tau/C_*$, one has
\bea\label{GT823next} \omega_{x_0}(\nu R) \leq
\frac{C-1}{C+1}\Big(\omega_{x_0}(R) + \hat{\Lambda}_0R^\lambda\Big)
\quad\text{for all $R\in(0,R_0]$}.
\eea

We now iterate \eqref{GT823next} using powers of $\nu$ to obtain
essential H\"older continuity of $u$.  Indeed, for any $\nu\leq
\tau/C_*$ and $j\geq 1$, we have \bea
\label{it1}\omega_{x_0}(\nu^jR_0) \leq
\Big(\frac{C-1}{C+1}\Big)^j\Big\{\omega_{x_0}(R_0) +
\hat{\Lambda}_0R_0^\lambda \disp\sum_{i=0}^{j-1}
\Big[\frac{C+1}{C-1}\, \nu^\lambda \Big]^i\Big\}. \eea Now choose
$\nu\leq \frac{\tau}{C_*}$ such that $\frac{C+1}{C-1}\nu^{\lambda}
\leq \frac{1}{2}$ to obtain $\disp\sum_{i=0}^\infty
\Big[\frac{C+1}{C-1} \, \nu^{\lambda}\Big]^{i} \leq 2$.  Then
\eqref{it1} gives \bea\label{it2} \omega_{x_0}(\nu^jR_0) \leq
\Big(\frac{C-1}{C+1}\Big)^j\Big[\omega_{x_0}(R_0) +
2\hat{\Lambda}_0R_0^\lambda \Big]. \eea Now let $0<R\leq R_0\nu$ and
choose $j\in\mathbb{N}$ such that $\nu^{j+1}R_0<R\leq \nu^jR_0$.
This choice implies that \bea\label{jsize}
j+1>\frac{\ln\Big(\frac{R}{R_0}\Big)}{\ln \nu}. \eea Combining
\eqref{jsize} with \eqref{it2} we obtain \bea\label{HC10}
\omega_{x_0}(R) \leq \frac{C+1}{C-1}\Big(\frac{R}{R_0}\Big)^\mu
\Big(\omega_{x_0}(R_0)+2\hat{\Lambda}_0R_0^\lambda\Big), \eea where
$\mu=\frac{\ln\frac{C-1}{C+1}}{\ln\nu}>0$. Thus there are positive
constants $c_7,\;\mu$ independent of $x_0$ such that \bea
\label{HC11}\omega_{x_0}(R)\leq c_7R^\mu \quad\text{if $0<R\le
R_0\nu$.} \eea As a consequence of \eqref{HC11}, $u$ is essentially
H\"older continuous with respect to $\rho$ in $B=B(y,\frac{\tau^2
}{5\kappa }r_1(y))$.  To see this, first note that since $u\in
L^\infty(B)$, there is a set $E_y\subset B$ with $|E_y|=0$ such that
\bea  |u(x)| \leq ||u||_{L^\infty(B)} \eea for all $x\in B\setminus
E_y$. Choosing $x,w\in B\setminus E_y$, there
are two cases to consider.\\

\noindent \emph{Case I:} $\rho(x,w)<\frac{\nu R_0}{2}$. Applying
\eqref{HC11} in the ball $B(x,2\rho(x,w))$ we obtain \bea
|u(x)-u(w)| \leq \omega_x(2\rho(x,w)) \leq c_7 2^\mu\rho(x,w)^\mu.
\eea

\noindent \emph{Case II:} $\rho(x,w)\geq \frac{\nu R_0}{2}$. Then
\bea |u(x)-u(z)|\leq 2||u||_{L^\infty(B)} \leq
\frac{2^{\mu+1}||u||_{L^\infty(B)}}{\nu^\mu R_0^\mu}\rho(x,w)^\mu.
\eea

Setting $c_8=\max\{c_7 2^\mu,
\frac{2^{\mu+1}\|u\|_{L^\infty(B)}}{\nu^\mu R_0^\mu}\}$ and
combining estimates, it follows that $u$ is essentially H\"older
continuous with respect to $\rho$ in $B$, which completes the proof.
\begin{flushright}$\Box$
\end{flushright}

\section{Proof of Corollary \ref{holder2}}

Fix a compact set $K\subset \Omega$.  By hypothesis, there is a positive
constants $s_0$ depending only on $K$ such that $s_0\leq
r_1(y)\leq 1$ for all $y\in K$.  As a result, the constants
$\Lambda$ and $M$ of Proposition \ref{pro1} can be chosen larger so
that they depend only on $K$ and $S=\disp\bigcup_{y\in K} B(y,r_1(y))
\Subset\Omega$.  More precisely, this is achieved by replacing in
those definitions all instances of $|B(y,r_1(y))|$ with
$\disp\inf_{y\in K}|B(y,r_1(y))|>0$, expanding all norms calculated on
$B(y,r_1(y))$ so that they are calculated over $S$, and replacing
$r_1(y)$ itself by $s_0$ and $1$ as appropriate.  Moreover, since
$\bar{S}$ is a compact subset of $\Omega$, $r_1$ satisfies a local uniformity
condition on every ball $B(y,r_1(y))$ with $y\in K$, with a uniform
constant $A_*=A_*(S)$.  In fact, one can choose $A_*=s_0'$ where
$s_0'$ satisfies  $0<s_0'\leq r_1(y)\leq 1$ for
all $y\in S$.  Combining these observations with Theorem
\ref{harnackmother}, it follows that any weak solution $(u,\nabla u)$
of \eqref{eqdiff} in $\Omega$ satisfies the Harnack inequality
\eqref{harnack} when $x_0\in K$ and where the constant $C_1$ there is
chosen to depend only on $K$ and $S$ via the norms of structural
coefficients and the $L^{p\sigma}(S)$ norm of $u$.\\

Fix now a weak solution $(u,\nabla u)$ of \eqref{eqdiff} in
$\Omega$. Working through the proof of Theorem \ref{holder1} with
the observations above, one sees that $\nu$ (see \eqref{it1}) can
now be chosen to depend only on $K$ and $S$ as $C$ there depends
only on these quantities, and $\tau/C_*$ depends only on $K, S$
through $C_*=C_*(A_*)$. As a result, for every $y\in K$ we have the
estimate \bea\label{globholder} \disp\sup_{z,w\in B(y,r)\setminus
E_y} \frac{|u(z)-u(w)|}{\rho(z,w)^\mu} \leq c_9, \eea where
$r=\frac{\tau^2}{5\kappa} s_0$ and the constants $c_9,\mu$ are
independent of $y$ and $E_y$ is as in the proof Theorem
\ref{holder1}. We now cover $K$ with a finite collection of
$\rho$-balls of the form $B(y_j,\frac{r}{2\kappa})$ and set
$E=\bigcup E_{y_j}$.  Then $|E|=0$,
and for $x,z\in K\setminus E$ there are two cases to consider:\\

\noindent \emph{Case I:} $\rho(x,z) < \frac{r}{2\kappa}$.  We claim
that there exists $y_j\in K$ such that both $x,z\in B(y_j,r)$.  Indeed,
choose $y_j$ such that $x\in B(y_j,\frac{r}{2\kappa})$.  Then
\begin{equation*} \rho(z,y_j)\leq \kappa(\rho(z,x)+\rho(x,y_j)) < r.
\end{equation*}
Since $x,z\in B(y_j,r)\setminus E_{y_j}$, we may apply
\eqref{globholder} to obtain
\begin{equation*} |u(x)-u(z)| \leq c_9\rho(x,z)^\mu.
\end{equation*}

\noindent \emph{Case II:} $\rho(x,z)\geq \frac{r}{2\kappa}$.  Arguing
as at the end of the proof of Theorem \ref{holder1}, we have
\begin{equation*} |u(x)-u(z)| \leq 2||u||_{L^\infty(S)} \leq \frac{2^{\mu+1}\kappa^\mu}{r^\mu}||u||_{L^\infty(S)}\rho(x,z)^\mu.
\end{equation*}

Combining both cases, it follows that $u$ is essentially H\"older
continuous in $K$ and, therefore, essentially locally H\"older
continuous in $\Omega$.

\begin{flushright}
$\Box$
\end{flushright}

\section{Proofs of results in Subsection \ref{subsec3.3}}

\textbf{Proof of Theorems \ref{harnackmother.2}, \ref{holder1.2} and
of Corollary \ref{holder2.2}.} For every $w\in[0,\infty)$ and every $\alpha\in(0,p]$ we have
\begin{equation}\label{12}
w^\alpha\leq1+w^p.
\end{equation}
Now, if $A(x,z,\xi)$ and $B(x,z,\xi)$ satisfy the structural
assumptions \eqref{struct} with $\gamma,\delta,\psi$ satisfying
\eqref{ranges3}, then by \eqref{12} they also satisfy the modified
structural conditions
\begin{equation}\label{struct1}
\begin{cases}
(i)\quad A(x,z,\xi)=\sqrt{Q(x)}{\tilde A}(x,z,\xi),\\
(ii)\quad\xi \cdot A(x,z,\xi) \geq a^{-1}\Big|\sqrt{Q(x)}\,
\xi\Big|^p
- h(x)|z|^p - (g(x)+h(x)),\\
(iii)\quad\Big|\tilde{A}(x,z,\xi)\Big| \leq a\Big|\sqrt{Q(x)}\,
\xi\Big|^{p-1} + b(x)|z|^{p -1} + (b(x)+ e(x)),\\
(iv)\quad\Big|B(x,z,\xi)\Big|\leq c(x)\Big|\sqrt{Q(x)}\,
\xi\Big|^{p-1} +d(x)|z|^{p-1}+(c(x)+d(x)+f(x))
\end{cases}
\end{equation}
for every $x\in\Omega$, every $z\in\R$ and every $\xi\in\R^n$. Thus,
in order to conclude, it is sufficient to apply Theorems
\ref{harnackmother}, \ref{holder1} and Corollary \ref{holder2} using
the new structural conditions \eqref{struct1}, i.e. with
$\gamma=\delta=\psi=p$ and with $e$, $f$, $g$ replaced respectively
by $b+e$, $c+d+f$ and $g+h$.
\begin{flushright}
$\Box$
\end{flushright}

\textbf{Proof of Theorems \ref{harnackmother.3}, \ref{holder1.3} and
of Corollary \ref{holder2.3}.} As is clear from their proofs, in
order to obtain Theorems \ref{harnackmother}, \ref{holder1} and
Corollary \ref{holder2}, one needs the structural assumptions
\eqref{struct} to hold with $\gamma=\delta=p$ and
$\psi\in[p,p+1-\sigma^{-1})$ not for every
$(x,z,\xi)\in\Omega\times\R\times\R^n$, but only for
$(x,z,\xi)=(x,u(x),\nabla u(x))$ for almost every $x\in\Omega$,
where $u$ is the weak solution of equation \eqref{eqdiff} under
consideration.

If $A(x,z,\xi)$ and $B(x,z,\xi)$ satisfy the structural assumptions
\eqref{struct} with $\gamma,\delta,\psi>p$, then we can write
\begin{equation}\label{struct3}
\begin{cases}
(i)\quad A(x,z,\xi)=\sqrt{Q(x)}{\tilde A}(x,z,\xi),\\
(ii)\quad\xi \cdot A(x,z,\xi) \geq a^{-1}\Big|\sqrt{Q(x)}\,
\xi\Big|^p
- \big(h(x)|z|^{\gamma-p}\big)|z|^p - g(x),\\
(iii)\quad\Big|\tilde{A}(x,z,\xi)\Big| \leq a\Big|\sqrt{Q(x)}\,
\xi\Big|^{p-1} + \big(b(x)|z|^{\gamma -p}\big)|z|^{p-1} + e(x),\\
(iv)\quad\Big|B(x,z,\xi)\Big|\leq c(x)\Big|\sqrt{Q(x)}\,
\xi\Big|^{\psi-1} +\big(d(x)|z|^{\delta-p}\big)|z|^{p-1}+f(x).
\end{cases}
\end{equation}
Now, by replacing $z$ with $u(x)$ in \eqref{struct3} we can conclude the proof through the application of Theorems \ref{harnackmother}, \ref{holder1} and Corollary
\ref{holder2}.  This is done using the modified structural conditions
\eqref{struct3} that correspond to \eqref{struct} with
$\gamma=\delta=p$ and with $h$, $b$, $d$ replaced by
$h_1=h|u|^{\gamma-p}$, $b_1=b|u|^{\gamma -p}$ and
$d_1=d|u|^{\delta-p}$ respectively.  Indeed, note that the map
$\mathcal{B}_0\mapsto\frac{p\sigma\mathcal{B}_0}{p\sigma+(\gamma-p)\mathcal{B}_0}$
is increasing and hence
$\frac{p\sigma\mathcal{B}_0}{p\sigma+(\gamma-p)\mathcal{B}_0}\geq\max\{p'\sigma',\frac{d_0}{p-1}\}$ since
$\mathcal{B}_0\geq\max\{\frac{p\sigma}{p\sigma-\sigma-\gamma+1},\frac{d_0p\sigma}{p(p-1)\sigma-d_0(\gamma-p)}\}$.
Thus we conclude that
$$\mathcal{B}=\min\Big\{\frac{p\sigma\mathcal{B}_0}{p\sigma+(\gamma-p)\mathcal{B}_0},\mathcal{E}\Big\}\geq\max\Big\{p'\sigma',\frac{d_0}{p-1}\Big\}$$
and, since $\mathcal{B}\leq\mathcal{E}$, we have $e\in
L^\mathcal{B}_{\text{loc}}(\Omega)$. Moreover for any compact subset
$\Theta\subset\Omega$ with positive measure
\begin{eqnarray*}
\|b_1\|_{\mathcal{B},\Theta,\overline{dx}}\!&\!\leq\!&\!\|b_1\|_{\frac{p\sigma\mathcal{B}_0}{p\sigma+(\gamma-p)\mathcal{B}_0},\Theta,\overline{dx}}\,=\,\left(\fint_\Theta
b^\mathcal{B}|u|^{(\gamma
-p)\mathcal{B}}\,dx\right)^\frac{1}{\mathcal{B}}\\
\!&\!\leq\!&\!\left(\fint_\Theta
b^\frac{p\sigma\mathcal{B}_0}{p\sigma+(\gamma-p)\mathcal{B}_0}|u|^\frac{(\gamma
-p)p\sigma\mathcal{B}_0}{p\sigma+(\gamma-p)\mathcal{B}_0}\,dx\right)^\frac{p\sigma+(\gamma-p)\mathcal{B}_0}{p\sigma\mathcal{B}_0}.
\end{eqnarray*}
By using H\"{o}lder's inequality with conjugate exponents
$q=\frac{p\sigma+(\gamma-p)\mathcal{B}_0}{p\sigma} $
and $q'=\frac{p\sigma+(\gamma-p)\mathcal{B}_0}{(\gamma-p)\mathcal{B}_0}$ we obtain
\begin{equation}\label{13}
\|b_1\|_{\mathcal{B},\Theta,\overline{dx}}\,\,\leq\,\,\left(\fint_\Theta
b^{\mathcal{B}_0}\,dx\right)^\frac{1}{\mathcal{B}_0}\left(\fint_\Theta
|u|^{p\sigma}\,dx\right)^\frac{\gamma-p}{p\sigma}\,\,=\,\,\|b\|_{\mathcal{B}_0,\Theta,\overline{dx}}\|u\|_{p\sigma,\Theta,\overline{dx}}^{\gamma-p}\,\,<\,\,+\infty,
\end{equation}
and hence $b_1=b|u|^{\gamma-p}\in
L^\mathcal{B}_{\text{loc}}(\Omega)$.

For $\mathcal{H}$ and $\mathcal{D}$ as defined in Theorem \ref{harnackmother.3}, one can prove in a similar way that
$\mathcal{H},\mathcal{D}\geq\frac{d_0}{p}$, that
$\mathcal{H},\mathcal{D}>\sigma'$, that $g,h_1=h|u|^{\gamma-p}\in
L^\mathcal{H}_{\text{loc}}(\Omega)$, and that $f,d_1=d|u|^{\delta-p}\in
L^\mathcal{D}_{\text{loc}}(\Omega)$ with
\begin{equation}\label{17}
\|h_1\|_{\mathcal{H},\Theta,\overline{dx}}\,\,\leq\,\,\|h\|_{\mathcal{H}_0,\Theta,\overline{dx}}\|u\|_{p\sigma,\Theta,\overline{dx}}^{\gamma-p}\qquad\text{
and }\qquad
\|d_1\|_{\mathcal{D},\Theta,\overline{dx}}\,\,\leq\,\,\|d\|_{\mathcal{D}_0,\Theta,\overline{dx}}\|u\|_{p\sigma,\Theta,\overline{dx}}^{\delta-p}.
\end{equation}
Finally, if $M$ is defined as in the statement of Proposition
\ref{pro1} (with $b,d,h$ replaced by $b_1,d_1,h_1$ respectively) and
$M_2$ is as in \eqref{18}, then \eqref{13} and \eqref{17}
imply that $M\leq M_2$.

Thus we can conclude by applying Theorems \ref{harnackmother},
\ref{holder1} and Corollary \ref{holder2}.
\begin{flushright}
$\Box$
\end{flushright}


\section{Appendix}

\begin{thm}\label{youngsinequality} (Young's Inequality) Let
$a_1,a_2,\theta >0$ and $\beta,\beta'\geq 1$ satisfy
$\frac{1}{\beta}+\frac{1}{\beta'} =1$.  Then
\bea\label{young} a_1\, a_2&\leq& \theta\,\fracc{a_1^{\beta}}{\beta} +
\fracc{1}{\theta^{\beta'/\beta}}\,\fracc{a_2^{\beta'}}{\beta'}.
\eea
\end{thm}

\begin{lem}\label{translate} Let $(\Omega,\rho,dx)$ be a local
homogeneous space and $\gamma^*$ be as in \eqref{swallowing}.  Fix
$x,y\in\Omega$, $\lambda\geq 1$, $t>0$, $l\in\Z$, and
$k\in\N\cup\{0\}$. Then if $t\leq \lambda^{l+k} < R_1(x)/\gamma^*$ and
$B(y,t)\cap B(x,\lambda^{l+k})\neq \emptyset$, we have
\bea\label{translate1}
|B(y,t)| \leq C_0(\gamma^*\lambda^k)^{d_0}|B(x,\lambda^l)|.
\eea
\end{lem}

{\bf Proof:} The swallowing property (\ref{swallowing}) gives $B(y,t)
\subset B(x,\gamma^*\lambda^{l+k})$, and since
$\gamma^*\lambda^{l+k} < R_1(x)$, we have that
\bea |B(y,t)| \leq |B(x,\gamma^*\lambda^{l+k})| \leq C_0\Big(\frac{\gamma^*\lambda^{l+k}}{\lambda^l}\Big)^{d_0}|B(x,\lambda^l)|= C_0(\gamma^* \lambda^k)^{d_0}|B(x,\lambda^l)|.\nonumber
\eea
\begin{flushright}
$\Box$
\end{flushright}

\begin{pro}
Let $(\Omega,\rho, dx)$ be a local homogeneous space, see Definition
\ref{homspace}, let $\Theta\Subset\Omega$ and assume $r_1(x)$ is a
function as in \eqref{rone} that satisfies a local uniformity
condition in $\Theta$ with constant $A_*=A_*(\Theta)$, see
\eqref{unifr}. Then condition weak-$D_{q^*}$, see Definition
\ref{weakdoub}, holds with $q^*=d_0$ on $\Theta$, for some constant
$C_7>0$ and with $\alpha=A_*/2$.
\end{pro}
\textbf{Proof:} Since $\Theta\Subset\Omega$ is compact, we can cover
it with a finite number of pseudometric balls
$B(y_1,r_1(y_1)),\ldots,B(y_P,r_1(y_P))$ with
$y_1,\ldots,y_P\in\Theta$. Let $x\in\Theta$,
$r\in\big(0,\frac{A_*}{2}r_1(x)\big)$ and choose
$y_k\in\{y_1,\ldots,y_P\}$ such that $x\in B(y_k,r_1(y_k))$. Then,
conditions \eqref{unifr} and \eqref{rone} imply that
$0<r<\frac{A_*}{2}r_1(x)<\frac{r_1(y_k)}{2}< R_1(y_k)$.  Using \eqref{doubling1} we conclude that
$$|B(y_k,r_1(y_k))|\leq C_0\left(\frac{r_1(y_k)}{r}\right)^{d_0}|B(x,r)|.$$
It now follows that condition weak-$D_{q^*}$ holds with $q^*=d_0$,
$\alpha=A_*/2$ and
$$C_7=\frac{1}{C_0}\min_{k=1,\ldots,N}\left\{\frac{|B(y_k,r_1(y_k))|}{(r_1(y_k))^{d_0}}\right\}.$$
\begin{flushright}
$\Box$
\end{flushright}

\textbf{Proof of Proposition \ref{pro1}.} \emph{Step 1.}  We start by
recalling that if $x_0\in B\big(y,\frac{\tau}{5\kappa}r_1(y)\big)$,
$r\in\big(0,\frac{\tau A_*}{5\kappa C_*}r_1(y)\big)$ and $C_*$ is as
in \eqref{5}, then $B(x_0,C_*r)\Subset B(y, r_1(y))$; see Step 1 of
the proof of Proposition \ref{harnackmain}. Since by the definition of
$C_*$ we have $0<C_*r<r_1(y)\leq R_1(y)$, Definition
\ref{homspace} gives
\begin{equation*}
|B(y,r_1(y))|\leq
C_0\left(\frac{r_1(y)}{C_*r}\right)^{d_0}|B(x_0,C_*r)|.
\end{equation*}

\emph{Step 2.} We now prove \eqref{14}. By Step 1 and the definition
of $k(x_0,r)$,
\begin{eqnarray*}
k(x_0,r)\!&\!\!=\!\!&\!\left[\frac{(C_*r)^{p-1}}{|B(x_0,C_*r)|^\frac{1}{\mathcal{B}}}\|e\|_{\mathcal{B},B(x_0,C_*r);dx}\right]^\frac{1}{p-1}+
\left[\frac{(C_*r)^{p}}{|B(x_0,C_*r)|^\frac{1}{\mathcal{D}}}\|f\|_{\mathcal{D},B(x_0,C_*r);dx}\right]^\frac{1}{p-1}\\
\!&\!\!\!\!&\!\hspace{1cm}+
\left[\frac{(C_*r)^{p}}{|B(x_0,C_*r)|^\frac{1}{\mathcal{H}}}\|g\|_{\mathcal{H},B(x_0,C_*r);dx}\right]^\frac{1}{p}\\
\!&\!\!\leq\!\!&\!\left[\left(\frac{C_0r_1(y)^{d_0}}{|B(y,r_1(y))|}\right)^\frac{1}{\mathcal{B}}(C_*r)^{p-1-\frac{d_0}{\mathcal{B}}}
\|e\|_{\mathcal{B},B(y,r_1(y));dx}\right]^\frac{1}{p-1}\\
\!&\!\!\!\!&\!\hspace{1cm}+\left[\left(\frac{C_0r_1(y)^{d_0}}{|B(y,r_1(y))|}\right)^\frac{1}{\mathcal{D}}(C_*r)^{p-\frac{d_0}{\mathcal{D}}}\|f\|_{\mathcal{D},B(y,r_1(y));dx}\right]^\frac{1}{p-1}\\
\!&\!\!\!\!&\!\hspace{2cm}+\left[\left(\frac{C_0r_1(y)^{d_0}}{|B(y,r_1(y))|}\right)^\frac{1}{\mathcal{H}}(C_*r)^{p-\frac{d_0}{\mathcal{H}}}\|g\|_{\mathcal{H},B(y,r_1(y));dx}\right]^\frac{1}{p}.
\end{eqnarray*}
Thus, by the definitions of $\lambda$, $\Lambda$ and the fact that
$C_*r<r_1(y)$,
\begin{eqnarray*}
k(x_0,r)\!&\!\!\leq\!\!&\!\left(\frac{C_0}{|B(y,r_1(y))|}\right)^\frac{1}{(p-1)\mathcal{B}}
r_1(y)^{1-\lambda}(C_*r)^\lambda\|e\|_{\mathcal{B},B(y,r_1(y));dx}^\frac{1}{p-1}\\
\!&\!\!\!\!&\!\hspace{1cm}+\left(\frac{C_0}{|B(y,r_1(y))|}\right)^\frac{1}{(p-1)\mathcal{D}}r_1(y)^{\frac{p}{p-1}-\lambda}(C_*r)^\lambda\|f\|_{\mathcal{D},B(y,r_1(y));dx}^\frac{1}{p-1}\\
\!&\!\!\!\!&\!\hspace{2cm}+\left(\frac{C_0}{|B(y,r_1(y))|}\right)^\frac{1}{p\mathcal{H}}r_1(y)^{1-\lambda}(C_*r)^\lambda\|g\|_{\mathcal{H},B(y,r_1(y));dx}^\frac{1}{p}\,\,\,=\,\,\,\Lambda
r^\lambda.
\end{eqnarray*}

\emph{Step 3.} If $B(x,l)\subseteq B(x_0,C_*r)$ and $0<l\leq C_*r$,
then $B(x,l)\Subset B(y,r_1(y))$ and $0<l< r_1(y)\leq R_1(y)$. As in
Step 1 we conclude that
\begin{equation*}
|B(y,r_1(y))|\leq C_0\left(\frac{r_1(y)}{l}\right)^{d_0}|B(x,l)|.
\end{equation*}
Also note that since $r_1$ satisfies a local uniformity
condition on $B(y,r_1(y))$ with respect to $\rho$ with constant
$A_*$ and since $x_0\in B(y,r_1(y))$, then $A_*r_1(y)<r_1(x_0)$.
Thus, from $$0<l\leq C_*r<A_*r_1(y)<r_1(x_0)\leq R_1(x_0)$$ and
$B(x,l)\subset B(x_0,C_*r)$, we deduce by Definition \ref{homspace}
that
\begin{equation*}
|B(x_0,C_*r)|\leq C_0\left(\frac{C_*r}{l}\right)^{d_0}|B(x,l)|.
\end{equation*}
\emph{Step 4.} We are now going to prove that
$\bar{Z}(B(x,l),\bar{u})\leq M$, with $M$ as in the conclusion of
Proposition \ref{pro1}. Using the definition of $\bar{Z}$ given in
\eqref{Zetabar} and equations \eqref{15}, we have
\begin{eqnarray}
\label{16}\bar{Z}(B(x,l),\bar{u})\!\!&\!\!\leq\!\!&\!\!1+l^{p-1}\|b\|_{p'\sigma',B(x,l);\overline{dx}}+\frac{l^{p-1}}{k^{p-1}}\|e\|_{p'\sigma',B(x,l);\overline{dx}}\\
\nonumber\!\!&\!\!\!\!&\!\!+\left(l^{p}\|c^\frac{p}{p+1-\psi}\bar{u}^\frac{p(\psi-p)}{p+1-\psi}\|_{\frac{p\sigma'}{p-\epsilon_1},B(x,l);\overline{dx}}\right)^\frac{1}{\epsilon_1}\\
\nonumber\!\!&\!\!\!\!&\!\!+\left(l^{p}\|h\|_{\frac{p\sigma'}{p-\epsilon_2},B(x,l);\overline{dx}}+\frac{l^{p}}{k^{p}}\|g\|_{\frac{p\sigma'}{p-\epsilon_2},B(x,l);\overline{dx}}\right)^\frac{1}{\epsilon_2}\\
\nonumber\!\!&\!\!\!\!&\!\!+\left(l^{p}\|d\|_{\frac{p\sigma'}{p-\epsilon_3},B(x,l);\overline{dx}}+\frac{l^{p}}{k^{p-1}}\|f\|_{\frac{p\sigma'}{p-\epsilon_3},B(x,l);\overline{dx}}\right)^\frac{1}{\epsilon_3}.
\end{eqnarray}
\emph{Step 5.} Recalling the conditions on $\mathcal{B}$ and by Step
3, we have by H\"{o}lder's inequality that
\begin{eqnarray*}
 l^{p-1}\|b\|_{p'\sigma',B(x,l);\overline{dx}}&\leq&l^{p-1}\|b\|_{\mathcal{B},B(x,l);\overline{dx}}\,\,=\,\,\frac{l^{p-1}}{|B(x,l)|^\frac{1}{\mathcal{B}}}\|b\|_{\mathcal{B},B(x,l);dx}\\
 &\leq&\frac{C_0^\frac{1}{\mathcal{B}}r_1(y)^\frac{d_0}{\mathcal{B}}l^{p-1-\frac{d_0}{\mathcal{B}}}}{|B(y,r_1(y))|^\frac{1}{\mathcal{B}}}\|b\|_{\mathcal{B},B(y,r_1(y));dx}\\
 &\leq&C_0^\frac{1}{\mathcal{B}}r_1(y)^{p-1}\|b\|_{\mathcal{B},B(y,r_1(y));\overline{dx}}.
\end{eqnarray*}
The terms including norms of $h$ and $d$ are treated in a similar
way, also recalling the definitions of $\epsilon_2,\epsilon_3$. Thus
we obtain
\begin{eqnarray*}
 l^p\|h\|_{\frac{p\sigma'}{p-\epsilon_2},B(x,l);\overline{dx}}&\leq&C_0^\frac{1}{\mathcal{H}}r_1(y)^p\|h\|_{\mathcal{H},B(y,r_1(y));\overline{dx}},\\
 l^p\|d\|_{\frac{p\sigma'}{p-\epsilon_3},B(x,l);\overline{dx}}&\leq&C_0^\frac{1}{\mathcal{D}}r_1(y)^p\|d\|_{\mathcal{D},B(y,r_1(y));\overline{dx}}.
\end{eqnarray*}
\emph{Step 6.} Again using the conditions on $\mathcal{B}$, Step 3 and the
definition of $k=k(x_0,r)$, we have
\begin{eqnarray*}
 \frac{l^{p-1}}{k^{p-1}}\|e\|_{p'\sigma',B(x,l);\overline{dx}}&\leq&\frac{l^{p-1}}{(C_*r)^{p-1}\|e\|_{\mathcal{B},B(x_0,C_*r);\overline{dx}}}\|e\|_{\mathcal{B},B(x,l);\overline{dx}}\\
 &=&\frac{l^{p-1}}{(C_*r)^{p-1}}\,\frac{|B(x_0,C_*r)|^\frac{1}{\mathcal{B}}}{\|e\|_{\mathcal{B},B(x_0,C_*r);dx}}\,\frac{\|e\|_{\mathcal{B},B(x,l);dx}}{|B(x,l)|^\frac{1}{\mathcal{B}}}\\
 &\leq&\frac{l^{p-1}}{(C_*r)^{p-1}}\,\frac{|B(x_0,C_*r)|^\frac{1}{\mathcal{B}}}{|B(x,l)|^\frac{1}{\mathcal{B}}}
 \,\,\leq\,\,C_0^\frac{1}{\mathcal{B}}\left(\frac{l}{C_*r}\right)^{p-1-\frac{d_0}{\mathcal{B}}}\,\,\leq\,\,C_0^\frac{1}{\mathcal{B}}.
\end{eqnarray*}
The terms involving $g$ and $f$ can be estimated in a similar way,
giving
$$\frac{l^p}{k^p}\|g\|_{\frac{p\sigma'}{p-\epsilon_2},B(x,l);\overline{dx}}\leq C_0^\frac{1}{\mathcal{H}}\qquad\text{ and }\qquad
 \frac{l^p}{k^{p-1}}\|f\|_{\frac{p\sigma'}{p-\epsilon_3},B(x,l);\overline{dx}}\leq C_0^\frac{1}{\mathcal{D}}.$$
\emph{Step 7.} We estimate the remaining term
$$I:=l^{p}\|c^\frac{p}{p+1-\psi}\bar{u}^\frac{p(\psi-p)}{p+1-\psi}\|_{\frac{p\sigma'}{p-\epsilon_1},B(x,l);\overline{dx}}
=l^p\left(\fint_{B(x,l)}c^\frac{p^2\sigma'}{(p+1-\psi)(p-\epsilon_1)}\,\bar{u}^\frac{p^2\sigma'(\psi-p)}{(p+1-\psi)(p-\epsilon_1)}dx\right)^\frac{p-\epsilon_1}{p\sigma'}$$
starting with an application of H\"{o}lder inequality with conjugate
exponents
$$q=\frac{(\sigma-1)(p+1-\psi)(p-\epsilon_1)}{p(\psi-p)}>1,\qquad
q'=\frac{(\sigma-1)(p+1-\psi)(p-\epsilon_1)}{(\sigma-1)(p+1-\psi)(p-\epsilon_1)-p(\psi-p)},$$
where we will associate $q$ with $\bar{u}$ and $q'$ with $c$; note
also that $q>1$ due to the definition of $\epsilon_1$ (see
Proposition \ref{pro1}). Thus, also recalling the conditions on $c$,
we obtain
\begin{eqnarray*}
 I&\leq&l^p
  \left(\fint_{B(x,l)}c^\frac{p^2\sigma}{(\sigma-1)(p+1-\psi)(p-\epsilon_1)-p(\psi-p)}dx\right)^\frac{(\sigma-1)(p+1-\psi)(p-\epsilon_1)-p(\psi-p)}{p\sigma(p+1-\psi)}
  \left(\fint_{B(x,l)}\bar{u}^{\sigma p}dx\right)^\frac{\psi-p}{(p+1-\psi)\sigma}\\
 &=&l^p\|c\|^\frac{p}{p+1-\psi}_{\frac{p^2\sigma}{(\sigma-1)(p+1-\psi)(p-\epsilon_1)-p(\psi-p)},B(x,l);\overline{dx}}
  \|\bar{u}\|^\frac{p(\psi-p)}{(p+1-\psi)}_{{p\sigma},B(x,l);\overline{dx}}\\
 &\leq&l^p\|c\|^\frac{p}{p+1-\psi}_{\mathcal{C},B(x,l);\overline{dx}}\|\bar{u}\|^\frac{p(\psi-p)}{(p+1-\psi)}_{{p\sigma},B(x,l);\overline{dx}}\,\,=\,\,
  \frac{l^p}{|B(x,l)|^{\frac{p}{(p+1-\psi)\mathcal{C}}+\frac{\psi-p}{\sigma(p+1-\psi)}}}\|c\|^\frac{p}{p+1-\psi}_{\mathcal{C},B(x,l);dx}\|\bar{u}\|^\frac{p(\psi-p)}{(p+1-\psi)}_{{p\sigma},B(x,l);dx},
\end{eqnarray*}
where the last inequality follows from the second part of the
minimum in the definition of $\epsilon_1$. Thus, by the first
display in Step 3,
\begin{eqnarray*}
 I&\leq&\frac{l^{p-\frac{d_0p}{(p+1-\psi)\mathcal{C}}-\frac{d_0(\psi-p)}{\sigma(p+1-\psi)}}(C_0r_1(y)^{d_0})^{\frac{p}{(p+1-\psi)\mathcal{C}}+\frac{\psi-p}{\sigma(p+1-\psi)}}}{|B(y,r_1(y))|^{\frac{p}{(p+1-\psi)\mathcal{C}}+\frac{\psi-p}{\sigma(p+1-\psi)}}}
   \|c\|^\frac{p}{p+1-\psi}_{\mathcal{C},B(y,r_1(y));dx}\|\bar{u}\|^\frac{p(\psi-p)}{(p+1-\psi)}_{{p\sigma},B(y,r_1(y));dx}\\
 &\leq&\frac{r_1(y)^pC_0^\frac{p\sigma+(\psi-p)\mathcal{C}}{(p+1-\psi)\sigma\mathcal{C}}}{|B(y,r_1(y))|^\frac{p\sigma+(\psi-p)\mathcal{C}}{(p+1-\psi)\sigma\mathcal{C}}}
   \|c\|^\frac{p}{p+1-\psi}_{\mathcal{C},B(y,r_1(y));dx}\left[\|u\|_{{p\sigma},B(y,r_1(y));dx}+k(x_0,r)|B(y,r_1(y))|^\frac{1}{p\sigma}\right]^\frac{p(\psi-p)}{(p+1-\psi)},
\end{eqnarray*}
where we have used the facts that $l<r_1(y)$ and
$p-\frac{d_0p}{(p+1-\psi)\mathcal{C}}-\frac{d_0(\psi-p)}{\sigma(p+1-\psi)}\geq
0$, due to the first condition on $\mathcal{C}$ in item (iv) of
Proposition \ref{pro1}. Finally, since $C_*r<r_1(y)$, Step 2 applies
to $k(x_0,r)$ and we have
$$I\leq r_1(y)^pC_0^\frac{p\sigma+(\psi-p)\mathcal{C}}{(p+1-\psi)\sigma\mathcal{C}}
   \|c\|^\frac{p}{p+1-\psi}_{\mathcal{C},B(y,r_1(y));\overline{dx}}\left[\|u\|_{{p\sigma},B(y,r_1(y));\overline{dx}}+\Lambda r_1(y)^\lambda\right]^\frac{p(\psi-p)}{(p+1-\psi)}.$$
\emph{Step 8.} It is now sufficient to insert the estimates from
Steps 5,6,7 into inequality \eqref{16} to conclude the proof.
\begin{flushright}
$\Box$
\end{flushright}

\bibliographystyle{plain}

\end{document}